\documentclass[11pt,draftcls,onecolumn]{IEEEtran}

\usepackage{cite}
\usepackage{amsmath} 
\usepackage{amssymb} 
\usepackage{dsfont}
\usepackage[mathcal]{euscript}
\usepackage{empheq}
\usepackage{graphicx}
\usepackage{graphics}
\usepackage{bm}
\usepackage{psfrag}
\usepackage{mathrsfs}
\usepackage{bbm}
\usepackage{color}
\usepackage{cancel}
\input{epsf}
\newtheorem{theorem}{Theorem}
\newtheorem{lemma}{Lemma}
\newtheorem{corollary}{Corollary}

\newtheorem{remark}{Remark}
\newtheorem{definition}{Definition}
\newtheorem{assumption}{Assumption}

\renewcommand{\P}{\mathbb{P}}
\newcommand{\E}{\mathbb{E}}

\newcommand{\beq}{\begin{equation}}
\newcommand{\eeq}{\end{equation}}
\newcommand{\beqa}{\begin{eqnarray}}
\newcommand{\eeqa}{\end{eqnarray}}
\newcommand{\dfz}{\triangleq}

\begin{document}

\title{
Graph Learning over Partially Observed Diffusion Networks: Role of Degree Concentration
}
\author{Vincenzo~Matta, Augusto~Santos, and Ali~H.~Sayed
\thanks{V.~Matta is with DIEM, University of Salerno,
via Giovanni Paolo II, I-84084, Fisciano (SA), Italy (e-mail: vmatta@unisa.it).

A. Santos was with the Adaptive Systems Laboratory, EPFL, CH-1015 Lausanne, Switzerland (e-mail: augusto.pt@gmail.com).

A.~H.~Sayed is with the \'Ecole Polytechnique F\'ed\'erale de Lausanne EPFL, School of Engineering, CH-1015 Lausanne, Switzerland (e-mail: ali.sayed@epfl.ch).
}
\thanks{
This paper was presented in part at the 2018 Asilomar Conference on Signals, Systems, and Computers~\cite{MattaSantosSayedAsilomar2018}, and at the 2019 IEEE International Symposium on Information Theory~\cite{MattaSantosSayedISIT2019}. 
}
\thanks{
The work of A.~H.~Sayed was supported in part by NSF grant CCF-1524250 and by grant 
205121-184999 from the Swiss National Science Foundation.
}
}

\maketitle

\begin{abstract}
This work examines the problem of graph learning over a diffusion network when data can be collected from a limited portion of the network ({\em partial observability}). 
The main question is to establish technical guarantees of consistent recovery of the subgraph of probed network nodes, $i)$ despite the presence of unobserved nodes; and $ii)$ under different connectivity regimes, including the dense regime where the probed nodes are influenced by many connections coming from the unobserved ones.
We ascertain that suitable estimators of the combination matrix (i.e., the matrix that quantifies the pairwise interaction between nodes) possess an {\em identifiability gap} that enables the discrimination between connected and disconnected nodes within the observable subnetwork.
Fundamental conditions are established under which the subgraph of monitored nodes can be recovered, with high probability as the network size increases, through {\em universal} clustering algorithms.
This claim is proved for three matrix estimators: $i)$ the Granger estimator that adapts to the partial observability setting the solution that is exact under full observability ; $ii)$ the one-lag correlation matrix; and $iii)$ the residual estimator based on the difference between two consecutive time samples. 
A detailed characterization of the asymptotic behavior (as the network size increases) of these matrix estimators is established in terms of an error bias and of the identifiability gap, through revealing closed-form solutions that embody the dependence on the main system features, such as level of observability and network connectivity. 
These estimators are then examined in terms of their {\em sample complexity}, i.e., in terms of how the number of samples grows with the network size to achieve consistent learning.
Comparison among the estimators is performed through illustrative examples that show how estimators that are not optimal in the full observability regime can outperform the Granger estimator in the partial observability regime.
The analysis reveals that the fundamental property enabling consistent graph learning is the {\em statistical concentration} of node degrees.
\end{abstract}


\section{Introduction}
\IEEEPARstart{L}{earning} the graph structure that governs the evolution of a networked dynamical system from data collected at some accessible nodes is a challenging inverse problem with applications across many domains. 
The objective of such inferential problems is to discover the interaction profile among the network nodes since the topology has a critical effect on system behavior~\cite{Barrat,liggett,Queues,porterdynamical}.
Graph learning plays a central role in many applications including, among other possibilities: estimating the longevity or the source of an epidemics~\cite{topoepidemics,Pedrosource}; revealing commonalities and agent influence over social networks~\cite{KrimetalTSIPN2016,EarnestSocialTSIPN2016,SalamiYingSayedTSIPN2017}; discovering the routes of clandestine information flows~\cite{VenkitasubramaniamHeTongTIT2008,EmbeddingCapacityTIT2013}; identifying defective elements~\cite{SaligramaetalTIT2012}; addressing the fundamental issue in neuroscience that links brain functional connectivity (i.e., a ``functional'' topology estimated from blood-oxygenation-level-dependent signals) to brain structural connectivity (i.e., the anatomical topology of neuron interconnections)~\cite{Honeyetal2009,bratislav}.

Depending on the particular context, the aforementioned class of problems can be referred to in different ways, including {\em topology inference}~\cite{mateos}, {\em network tomography}~\cite{tomo, tomo_journal}, {\em structure learning}, or {\em graph learning}~\cite{AnandkumarValluvanAOS2013}. We adopt these terminologies almost interchangeably throughout our treatment.

This article addresses the graph learning problem under the framework of {\em partial observability}, i.e., when only a fraction of the network nodes can be probed. Providing technical guarantees of graph learning under such a demanding setting is critical in large-scale networked systems, where it is not feasible to gather data from all nodes comprising the network. 
Furthermore, we solve the learning problem for distinct regimes of network wiring, including {\em densely-connected} networks, a case often overlooked in the literature and challenging to be addressed in general. 

We establish that, over a diffusion network, and under certain assumptions on the generative model --- i.e., on the entries of the combination matrix and its underlying support graph --- the problem of topology inference becomes {\em local}, i.e., all the information about the underlying subnetwork connecting the observed agents is contained in the samples of the observed agents, asymptotically as the network size gets large, and irrespective of the density of connections.

\section{Overview of the Learning Problem}
\label{sec:idhasa}
\begin{figure} [t]
\begin{center}
\includegraphics[scale= 0.27]{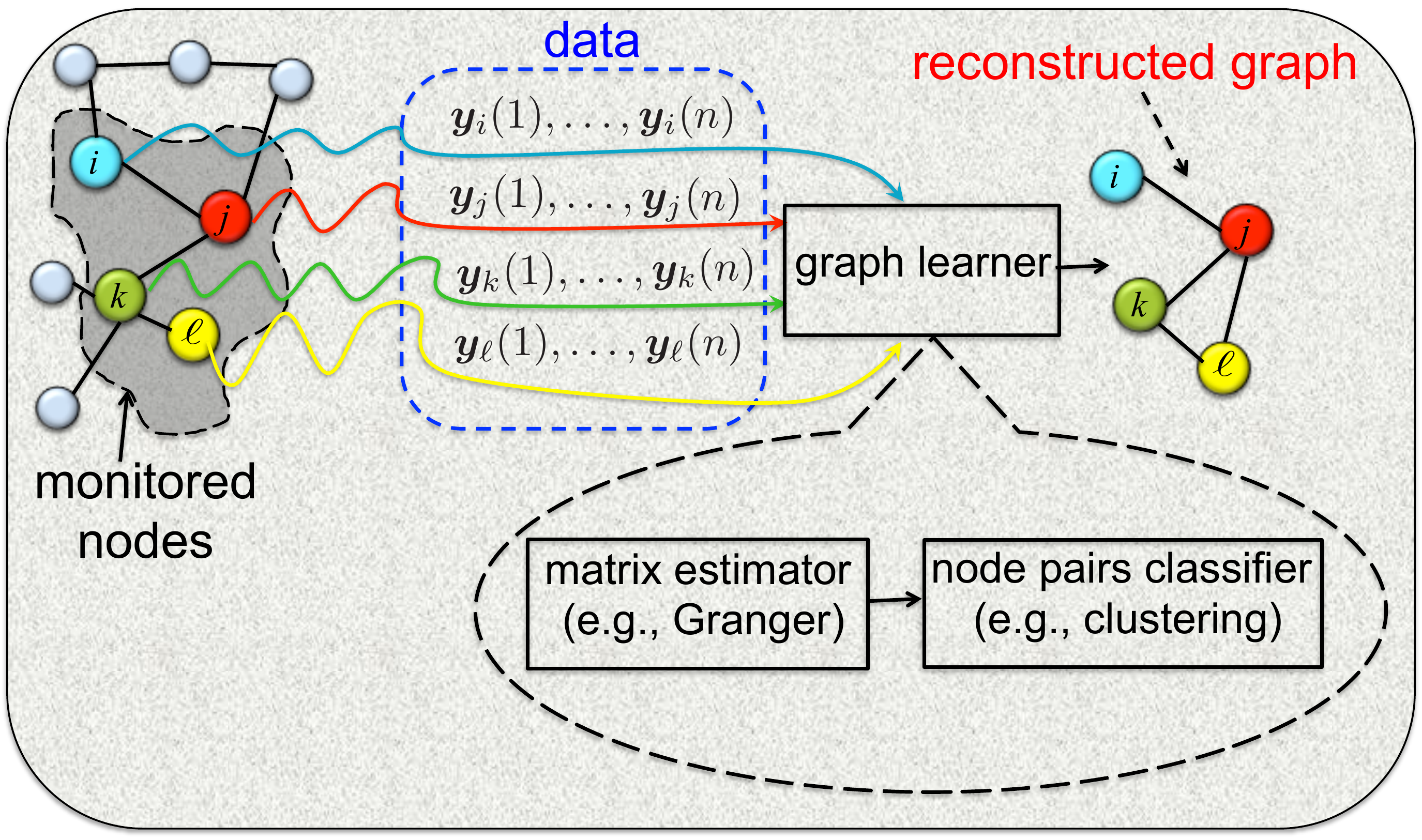}
\caption{Illustration of the graph learning problem considered in this work. 
}
\label{fig:scheme1}
\end{center}
\end{figure}

The structure of the learning problem addressed in this work can be summarized as follows. 
Streams of data originating from a certain subnetwork are collected, and the goal is to estimate the (unknown) topology linking the nodes of this subnetwork from the collected data. 
The graph learning protocol will involve two main steps: an estimation step, where a combination matrix (i.e., a matrix quantifying the strength of the connections among the network nodes) is estimated; and a thresholding step, where node pairs linked by a strong edge (i.e., node pairs whose corresponding combination matrix entry lies above some threshold) are deemed connected. 
A structurally-consistent estimator is one that ends up assigning strong ties to interacting pairs and weak ties to non-interacting pairs. In this way, at the thresholding stage, one can correctly classify the pairs as interacting and non-interacting. 
Figure~\ref{fig:scheme1} gives a graphical summary of the aforementioned procedure.

It is useful to illustrate the learning problem in relation to some popular networked systems. 
To start with, let us neglect some practical limitations, in particular: $i)$ assume that all nodes can be monitored (full observability); $ii)$ that there are no limitations in terms of computational power; and $iii)$ there are no limitations on the available time samples. 
Then, the first inferential stage consists of finding a matrix estimator to quantify the strength of pairwise interactions in the network. 
One notable estimator relies on computing the (spatial) correlation matrix $R_0=\lim_{n\rightarrow \infty}\mathbb{E}\left[\bm{y}_n\bm{y}_n^{\top}\right]$, where $\bm{y}_n$ denotes the vector collecting the data from all nodes at time $n$, and where assumption $iii)$ above justifies the limit and (under an ergodic assumption) implies that the statistical average can be learned from the data. 
When the matrix $R_0$ provides a consistent estimator for the connection strengths, we talk of a {\em correlation network}~\cite{mateos}. 
For these networks, interactions between two nodes are {\em direct} and they are directly captured from pairwise correlations. One example of a correlation network is the ferromagnetic Ising model~\cite{Bento2009} with independent and identically distributed (i.i.d.) time samples, and under certain constraints of sparsity on the network and of regularity on the weights of interaction. 

Another classic model for graph learning is a {\em Gaussian graphical model}. In this case, $R_0$ is no longer the proper estimator, but its inverse $R_0^{-1}$ (which is often referred to as the precision or concentration matrix) is a consistent estimator, in that its support coincides with the underlying graph of interactions. 
Over Gaussian graphical models, the pairwise interaction between adjacent nodes is affected by other nodes, and this latent influence is the reason why spatial correlation between measurements is no longer sufficient to capture the network structure. 

For most standard graphical models, pairwise interactions are described through dependent random variables defined on the network nodes. It is usually assumed that i.i.d. samples of these variables are available for the learning process. 
In other words, over graphical models the data samples do not arise from a dynamical process governing the time evolution of the nodes' outputs.
In contrast, there is the more interesting class of graph models that correspond to dynamical graph systems. 
In this case, the signal time-evolution at any node is influenced by the signal evolutions at neighboring nodes.
One relevant example is the diffusion or first-order Vector Autoregressive (VAR) system described by~(\ref{eq:VARmodel}) further ahead. 
For such graphs, the proper estimator for graph connectivity turns out to be the Granger estimator, $R_1 R_0^{-1}$, which combines in a suitable way information contained in the correlation matrix, $R_0$, and in the {\em one-lag} correlation matrix, $R_1$. 

In this article, we will be dealing with dynamic graphs, where signals evolve at the nodes and are affected by the evolution of the signals at neighboring agents as well.

\subsection{Structural-Consistency, Hardness and Sample-Complexity}
\label{subsec:schsc}
We are now ready to introduce three concepts that play an important role in graph learning problems.

\vspace*{5pt}
\noindent
{\bf {\em Structural consistency}}.
If it is possible to discover the graph structure from a statistical descriptor related to the measurements (in the previous examples, $R_0$, $R_0^{-1}$, or $R_1 R_0^{-1}$), we shall say that the graph learning problem is {\em identifiable}. 
When a graph learning problem is identifiable, the goal is to find a matrix estimator that is able to guess the correct graph. 
Given a matrix estimator that allows identifying (at least asymptotically as the network size goes to infinity) the graph structure, we shall say that this matrix estimator achieves {\em structural consistency}. 
We see that identifiability is a property of a given graph learning problem, whereas consistency is a property of a given matrix estimator applied to an identifiable problem.

In this work we focus on the case in which the measurements are available from only a limited subset of nodes, with identifiability referring here to the subgraph connecting these monitored nodes. 
The identifiability issue becomes particularly critical under this demanding setting\footnote{It should be remarked, however, that certain issues may arise also in the full observability case. For instance, as the network size increases, the entries of the matrix estimators might become smaller and smaller, and, hence, it is critical to identify whether a stable threshold can be found to classify connected/disconnected node pairs.}, giving rise to many interesting as well as challenging questions, such as: {\em Does partial observability impair identifiability of the monitored subnetwork? If not, how can we design structurally-consistent estimators?}

We remark that the concepts of identifiability and structural consistency disregard complexity issues, since they assume that the necessary statistical quantities are available. 
For example, we assume that $R_1 R_0^{-1}$ can be computed exactly, which means that we have sufficient time to learn and that matrix inversion is possible whatever the size of the network.
    
\vspace*{5pt}
\noindent
{\bf {\em Hardness}}.
How much computational complexity is required to evaluate the matrix estimator necessary for a particular graph learning problem? 
For example, if one is interested in the precision matrix, $R_0^{-1}$, the hardness is related to the complexity of the matrix inversion. 
Many works attempt to reduce the complexity by leveraging particular constraints such as smoothness or sparsity of the graph signals. 

Notably, there are relevant topology inference problems that are computationally intractable (e.g., NP-hard)~\cite{graphlearnhard1,SHIMONY,Breslerhardness,complexity_markov}.
It is therefore critical to identify meaningful graph classes where this does not happen~\cite{BreslerEfficiently}. 
The present work sheds light on a class of networked systems whereby structure learning with {\em non-independent} observations and under {\em partial} observability has affordable computational complexity, even over the hardly explored and nonetheless of critical importance {\em dense} (thus, loopy and non-tree like) networks.
However, we remark that the concept of hardness disregards the complexity associated with the {\em empirical} estimation of the pertinent matrices from the available data samples. This fundamental element of complexity is usually referred to as {\em sample} complexity.

\vspace*{5pt}
\noindent
{\bf {\em Sample complexity}}.
In practice, only a finite amount of data is available and, therefore, only approximate versions of the aforementioned matrix estimators can be computed. 
There exist several results about sample complexity in the context of high-dimensional graphical models, where the number of samples necessary to get some prescribed accuracy is related to the system parameters, e.g., to the network size and to the density of connections. 
Results relative to sample complexity are comparably less mature over {\em dynamical} graph systems~\cite{mateos}. Over these systems, the dependence among the time samples induced by the dynamical model complicates the theoretical analysis of the convergence of the empirical estimators, and, hence, the analysis of sample complexity. 
Useful results about the sample complexity of vector autoregressive models like the one addressed in this work are available in~\cite{BentoIbrahimiMontanari,HanLuLiuJMLR,LohWainwrightAOS2012,RaoKipnisJavidiEldarGoldsmithCDC2016}. These results do not consider the {\em partial-observations} setting and, hence, they do not apply here.\footnote{In order to avoid misunderstandings, we remark that the terminology ``partial samples'' or ``missing data''  used in~\cite{LohWainwrightAOS2012,RaoKipnisJavidiEldarGoldsmithCDC2016} refers to a different problem. In our setting samples from only a subset of the network are available, while in~\cite{RaoKipnisJavidiEldarGoldsmithCDC2016} samples from the whole network are available. However, given an overall amount of data generated at each node, a certain fraction per node is randomly lost and/or corrupted.}
Nonetheless, in~\cite{BentoIbrahimiMontanari,HanLuLiuJMLR,LohWainwrightAOS2012,RaoKipnisJavidiEldarGoldsmithCDC2016} we can find useful techniques to bound the errors associated to the empirical correlation matrices over vector autoregressive models. These types of bounds will be exploited in the proof of Theorem~\ref{theor:samplecomplexity}, where we establish the sample complexity of the estimators proposed in this work.

\section{Related Work}
The learning problem considered in this work lies in the broad field of signal processing over graphs~\cite{Sayed,SayedProcIEEE2014,MattaSayedCoopGraphSP2018,SayedZhaoCoopGraphSP2018,MouraDSP1,MouraDSP2,KovacevicSP2015,OrtegaSPmag,Ortegaetal2016,TsitsiveroBarbarossaDiLorenzo2016,MarquesSegarraLeusRibeiroSP2017,ChepuriLeusTSIPN2017,PerraudinVandergheynstSP2017}, dealing in particular with the identification of an unknown network topology from measurements gathered at the network nodes~\cite{mateos}.
These types of inferential problems can be addressed under different settings, including the case where measurements from all network nodes are available (full observability), and the case where only a fraction of nodes is accessible (partial observability).
Even if we focus on the partial observability setting, we deem it useful to start with some results pertaining to the full observability setting.

\subsection{Graph Learning under Full Observability}
The majority of works on graph learning over networks focuses on linear system dynamics, with nonlinear dynamics typically being tackled by variational characterizations under a small-noise assumption~\cite{Ching2017ReconstructingLI, Napoletani, Noise_Jien_Ren}, or by increasing the dimensionality of the observable space~\cite{ScienceRobustNetInference, Koopman_Gon}.

Topology inference for a general class of linear stochastic dynamical systems (e.g., VAR models of arbitrary order, or even non-causal linear models) is addressed in~\cite{MaterassiSalapakaTAC2012}. 
An approach based on Wiener filtering is proposed to infer the topology, which provides exact reconstruction for self-kin networks or, in general, guarantees reconstruction of the smallest self-kin network embracing the true network.

There exist works dealing with more general dynamical systems and graphs. For example, in~\cite{Kiyavash1} the concept of directed information graphs is advocated to discover dependencies in networks of interacting processes linked by causal dynamics. And a metric to learn causal relationships by means of {\em functional dependencies} is proposed in~\cite{Kiyavash2} to track the graph structure over a possibly nonlinear network.

Moving closer to our setting, among linear (or linearized) systems, special attention is devoted to autoregressive diffusion models~\cite{Moneta, pasdeloup, SantiagoTopo, MeiMoura}. 
For instance, in~\cite{MeiMoura} causal graph processes are exploited to devise a computationally tractable algorithm for graph structure recovery with performance guarantees.
Recent works exploit optimization and graph signal processing techniques to feed the graph learning algorithm with proper structural constraints. In~\cite{pasdeloup, SantiagoTopo}, it is shown how to capitalize on the fact that the weighting matrix and the correlation matrix share the same eigenvectors, and how to solve the topology inverse problem through optimization methods under sparsity constraints. An account of the methods for the full observability regime can be found in~\cite{mateos}. 

We stress that most of the aforementioned methods work in the {\em graph spectral domain}. This has to be contrasted with the methods proposed in this paper, which rely instead on the graph edge domain. 
Working in the edge domain allows us to obtain a transparent relationship for the matrix estimators, which is critical to establish identifiability under the challenging partial observability setting. 

In summary, while in certain cases (e.g., general linear models and/or nonlinear models) identifiability can be an issue, most of the works on diffusion models focus on reducing complexity by exploiting proper structural constraints (e.g., smoothness or sparsity of the signals defined on the graph)~\cite{mateos}.
However, all the aforementioned results pertain to the case where node measurements from the entire network are available. 
We focus instead on the case in which only partial observation of the network is permitted.

\subsection{Graph Learning under Partial Observability}
A fundamental challenge of our work is performing structure learning when only {\em partial} observation is allowed. 
Under this setting, results for retrieving particular network graphs (polytrees) are available in~\cite{MaterassiSalapakaCDC2012} and~\cite{KiyavashPolytrees}.
Considering instead the case of general topologies, and with focus on VAR/diffusion models like the one considered in this work, references~\cite{Geigeretal15, MaterassiSalapakaCDC2015} establish technical conditions for exact or partial topology identifiability.
However, these identifiability conditions act at a very ``microscopic'' level (they are formulated in terms of some precise details of the local graph structure and/or of the statistical model), and are therefore impractical over large-scale networks. 
In contrast, in this work we pursue a statistical asymptotic approach that is genuinely tailored to the large-scale setting: an asymptotic (i.e., large-network) framework is considered, where a thermodynamic limit of large graphs is afforded by using {\em random} graphs, and the conditions on the network connection topology are summarized at a {\em macroscopic} level through average descriptive indicators (e.g., probability of drawing an edge). 
In a similar vein, detailed asymptotic analysis with performance guarantees are available for graphical models with latent variables. For example, in~\cite{ChandrasekaranParriloWillskyAOS2012} it is shown that, under certain conditions concerning the interactions between the observed and the unobserved network nodes, the ``{\em sparsity+low-rank}'' framework can be exploited to estimate the amount of latent variables~\cite{ChandrasekaranParriloWillskyAOS2012}, and to reconstruct the topology of the observable subnetwork.
Likewise, in~\cite{AnandkumarValluvanAOS2013} the graph learning problem is tackled in the context of locally-tree graphs, whereas in~\cite{AnandkumarTanHuangWillskyJMLR2012} a local separation criterion is imposed to deal with Gaussian graphical models. 
Still in the framework of learning graphical models with latent variables, in~\cite{BreslerBoltzmann} an influence-maximization metric is proposed, to show that ferromagnetic restricted Boltzmann machines with bounded degree are an instance of graphical models that can be efficiently learned.

However, and as already explained in Sec.\ref{sec:idhasa}, graphical models do not match the networked dynamical models considered in this work. 
For these models, results for graph learning under partial observation have been recently obtained in~\cite{tomo, tomo_journal, tomo_icassp, tomo_dsw, tomo_isit}.
More specifically, $i)$ in~\cite{tomo} the whole network graph is assumed to follow an Erd\H{o}s-R\'enyi construction and the number of observable nodes grows with the overall network size $N$; whereas $ii)$ in~\cite{tomo_journal} the number of monitored nodes is held fixed (and, hence, the fraction of observable nodes vanishes in the limit of large $N$), the graph of the monitored nodes is left arbitrary, and the unobserved component continues to obey an Erd\H{o}s-R\'enyi model. 
The present work focuses on the former model\footnote{Even if the machinery used to prove our results can be applied to the latter model as well, we deem it useful to focus on a single model, in order to make the exposition more organic and to convey better the main message of the work.}.   
It is therefore necessary to explain clearly why the present work constitutes a significant progress in the context of local tomography over diffusion networks, in comparison to~\cite{tomo}.

\subsection{Main Advances}
In this work we consider the same setting adopted in~\cite{tomo}, except for a slightly more restrictive assumption on the class of combination matrices, which is summarized in Assumption~\ref{assum:regdiffmat} further ahead. 
These matrices, referred to as {\em regular diffusion matrices}, arise quite naturally when the nodes' combination weights are collected into a (scaled, i.e., stable version of a) symmetric doubly-stochastic matrix, whose support graph must match the underlying graph of connections. It is important to remark that all the matrices considered in the examples of the aforementioned previous works, as well as the most typical combination matrices, satisfy automatically such an assumption. The key contributions here in relation to~\cite{tomo} are as follows: 

{---} One first advance regards the regime of connectivity. 
Reference~\cite{tomo} addresses only the case that the network is sparsely connected, which means that the connection probability is allowed to vanish with $N$, but in a way that preserves network connectedness. 
In this work we are able to establish that consistent local tomography is possible also under the {\em dense} regime, which means that the connection probability is {\em not} vanishing. 
Our results are established in terms of a strong notion of consistency (referred to as {\em universal local structural consistency}), which will be formally introduced in Sec.~\ref{sec:ident}. 
Moreover, the analysis will reveal that the key-property to prove this consistency result is the {\em statistical concentration} of node degrees.

{---} We advance also with respect to the results currently available under the sparse regime. 
More specifically, in~\cite{tomo} a consistency result is proved, for all sparsely connected networks, in terms of an  {\em average} fraction of misclassified node pairs. 
In the present work we are able to strengthen significantly this consistency result for the subclass of sparsely connected networks whose node degrees fulfill a statistical concentration property. 
For these networks, we ascertain that universal local structural consistency holds. 
Using such stronger notion of consistency answers also the following important question (posed, and only partially answered in~\cite{tomo}, where the answer was obtained under a certain {\em approximation} of independence): {\em does the fraction of mistakes scale properly with $N$?} 

{---} We are able to offer a rigorous proof that the connected and non-connected agent pairs can be recovered through some {\em universal} clustering algorithm. In particular, we propose a variant of the $k$-means algorithm that is shown to be asymptotically consistent.

{---} The third advance regards the topology inference algorithms. 
Under the full-observability setting, it can be relatively easy, for a given dynamical system, to retrieve the right matrix estimator that contains information about the underlying graph. The three matrix estimators in Sec.~\ref{sec:idhasa} (i.e., $R_0$, $R_0^{-1}$, or $R_1 R_0^{-1}$) constitute examples of how to manage this aspect in traditional graph learning models. 
On one hand, the matrix estimators available in the full-observability setting might orient the choice of a matrix estimator for the partial observability setting. For example, one can replace the Granger estimator with a version that considers only the subnetwork of observed measurements. This choice is widely adopted in causal inference from time series (when one neglects the existence of latent components), and has been adopted, e.g., in~\cite{tomo} and~\cite{tomo_journal}.
On the other hand, there is in principle no reason why the best solution to the partial observability setting is obtained by mimicking one matrix estimator that is good for the full observability setting. 
This concept is put forth in the present work, where, exploiting novel analytical methods to characterize the graph learning errors in terms of powers of the combination matrix, we are able to construct and examine different matrix estimators. 
In particular we will consider, along with the Granger estimator, two other matrix estimators, namely, the {\em one-lag correlation matrix} and the {\em correlation matrix between the residuals} (i.e., difference between subsequent time samples). 
First, we will be able to show that all the three estimators are structurally consistent. 
Then, we will show by numerical simulations that the Granger estimator can be outperformed, which confirms the hypothesis that what is good under full observability might not be the best option to use under partial observability.

{---} Last but not least, we achieve the significant advance of ascertaining the {\em sample complexity} of the proposed estimators, i.e., we establish how the number of samples must scale with the network size to let these estimators learn the graph of probed nodes faithfully.

The bottom line of the novel results presented in this work is that the {\em fortunate coupling between the Erd\H{o}s-R\'enyi model and the regular diffusion matrices renders the problem of topology inference a local problem}.  
Moreover, the analysis will reveal that sparsity is not necessarily the key-enabler for structural consistency under partial observability. 
Instead, the main enabling feature will be seen to be the error concentration induced by the aforementioned coupling.

\subsection{Outline}
The article is organized as follows. 
Section~\ref{sec:dyngraph} introduces the dynamical system considered in this work, which encodes the relationships between the output measurements and the network topology. 
In Sec.~\ref{sec:ERgraphs} we describe the random model used for the graph, we present useful {\em statistical concentration} results that are critical for our treatment, and we introduce the {\em partial observation} setting. We remark that adopting a statistical approach is key to address the {\em large-scale} setting, since a {\em random} graph enables an asymptotic study of the thermodynamic behavior of the system as the network size grows. Our achievability results will be indeed established in a precise statistical sense, by showing that the estimated graph converges with high probability to the true graph. Following these lines, Sec.~\ref{sec:ident} introduces rigorous notions of consistency, and fundamental concepts such as margins, identifiability gap and bias, which are the building blocks to characterize in a precise statistical sense the limiting behavior of the inferential procedure. 
In relation to the concept of {\em universal} consistency, we prove the existence of  a clustering algorithm that enables universal (i.e., fully data-driven) recovery of the true topology.
In Sec.~\ref{sec:ME} we propose three matrix estimators, whose universal consistency is established in Sec.~\ref{sec:MainTh}. 
Finally, Sec.~\ref{sec:illustexam} collects some examples aimed at illustrating the theoretical results. 

\vspace*{5pt}
{\em Notation}. 
We use boldface letters to denote random variables, and normal font letters for their realizations. 
Matrices are denoted by capital letters, and vectors by small letters. 
This convention can be occasionally violated, for example, the total number of network nodes is denoted by $N$. 
The symbol $\stackrel{\textnormal{p}}{\longrightarrow}$ denotes convergence in probability as $N\rightarrow\infty$. 

Sets and events are denoted by upper-case calligraphic letters, whereas the corresponding normal font letter will denote the set cardinality. For example, the cardinality of ${\cal S}$ is $S$. The complement of $\mathcal{S}$ is denoted by $\mathcal{S}'$. 

For a $K\times K$ matrix $Z$, the submatrix spanning the rows of $Z$ indexed by set $\mathcal{S}\subseteq\{1,2,\ldots, K\}$ and the columns indexed by set $\mathcal{T}\subseteq\{1,2,\ldots,K\}$, is denoted by $Z_{\mathcal{S} \mathcal{T}}$, or alternatively by $[Z]_{\mathcal{S} \mathcal{T}}$. When $\mathcal{S}=\mathcal{T}$, the submatrix $Z_{\mathcal{S} \mathcal{T}}$ is abbreviated as $Z_{\mathcal{S}}$. 
Moreover, in the indexing of the submatrix we keep the index set of the original matrix. For example, if $\mathcal{S}=\{2,3\}$ and $\mathcal{T}=\{2,4,5\}$, the submatrix $M=Z_{\mathcal{S} \mathcal{T}}$ is a $2\times 3$ matrix, indexed as follows:
\beq
M=
\left(
\begin{array}{llll}
z_{22}&z_{24}&z_{25}\\
z_{32}&z_{34}&z_{35}
\end{array}
\right)
=\left(
\begin{array}{llll}
m_{22}&m_{24}&m_{25}\\
m_{32}&m_{34}&m_{35}
\end{array}
\right).
\eeq
For a matrix $M$, the symbols $\|M \|_1$, $\| M\|_2$ and $\|M\|_{\infty}$ denote the vector-induced $\ell_1$, $\ell_2$ and $\ell_{\infty}$ norms of $M$, respectively. 
The symbol $\| M \|_{\max}$ denotes instead the maximum absolute entry of $M$.
The symbol $\log$ denotes the natural logarithm.

\section{Dynamical Graph System}
\label{sec:dyngraph}
Let $\bm{y}_i(n)$ be the output measurement produced by node $i$ at time $n$. 
Likewise, let $\bm{x}_i(n)$ be the input source (e.g., streaming data or noise) exciting node $i$ at time $n$.
The random variables $\bm{x}_i(n)$ have zero mean and unit variance, and are independent and identically distributed (i.i.d.), both spatially (i.e., w.r.t. to index $i$) and temporally (i.e., w.r.t. to index $n$). 
It is convenient to stack the input and output variables, respectively, into the vectors:
\beqa
\bm{x}_n&=&[\bm{x}_1(n),\bm{x}_2(n),\ldots,\bm{x}_N(n)]^{\top}, \nonumber\\
\bm{y}_n&=&[\bm{y}_1(n),\bm{y}_2(n),\ldots,\bm{y}_N(n)]^{\top}.
\eeqa
The stochastic dynamical system considered in the present work is given by the following {\em network diffusion} process (a.k.a. first-order Vector AutoRegressive (VAR) model):
\beq
\boxed{
\bm{y}_n = A\,\bm{y}_{n-1}+ \sigma \, \bm{x}_n
}
\label{eq:VARmodel}
\eeq
where $A$ is some stable $N\times N$ matrix with nonnegative entries, and $\sigma^2$ is a variance factor.
By rewriting~(\ref{eq:VARmodel}) on an entrywise basis:
\beq
\bm{y}_{i}(n)=\sum_{\ell=1}^N a_{i \ell} \, \bm{y}_{\ell}(n - 1) + \sigma\,\bm{x}_i(n),
\label{eq:1storder}
\eeq
we readily see that the support-graph of $A$ reflects the connections among the network nodes. 
Indeed, Eq.~(\ref{eq:1storder}) shows that, at time $n$, the output of node $i$ is updated by {\em combining} the outputs of other nodes from time $n-1$. 
In particular, node $i$ scales the output of node $\ell$ by using a combination weight $a_{i\ell}$, which implies that the output of agent $\ell$ is {\em effectively used} by node $i$ if, and only if $a_{i\ell}\neq 0$. 
After the combination step, the output measurement $\bm{y}_i(n)$ is adjusted by incorporating the streaming-source value, $\bm{x}_i(n)$, which is locally available at node $i$ at current time $n$. 

In our {\em partial observation} setting, only a subset of nodes can be probed: for each node $i$ belonging to the subset of probed nodes, $\mathcal{S}$, a stream of $n$ measurements, $\bm{y}_i(1), \bm{y}_i(2),\ldots, \bm{y}_i(n)$ is acquired. 
The learning task is to reconstruct the graph of interconnections corresponding to the combination (sub)matrix $A_{\mathcal{S}}$.

Formulations like the one in~(\ref{eq:VARmodel}) arise naturally across application domains, e.g., in economics~\cite{Moneta}, in the variational characterization of nonlinear stochastic dynamical systems~\cite{PhysRevE.lai}, or in distributed network processing applications where several useful strategies such as consensus~\cite{XiaoBoydSCL2004,TsitsiklisBertsekasAthansTAC1986,BoydGhoshPrabhakarShahTIT2006,DimakisKarMouraRabbatScaglioneProcIEEE2010,MouraetalTSP2011,MouraetalTSP2012,KarMouraJSTSP2011,BracaMaranoMattaTSP2008,BracaMaranoMattaWillettTSP2010,SayedProcIEEE2014,ChenSayedTIT2015part1,ChenSayedTIT2015part2} and diffusion~\cite{LopesSayedTSP2008,CattivelliSayedTSP2010,CattivelliSayedTSP2011,SayedTuChenZhaoTowficSPmag2013,MattaBracaMaranoSayedTIT2016,MattaBracaMaranoSayedTSIPN2016} lead to data models of the form in~(\ref{eq:1storder}). 
The forthcoming section provides an example belonging to the latter category.

\subsection{One Example: Adaptive Diffusion Networks}
\label{subsec:motivation}
A team of $N$ agents is engaged in solving some inferential task (e.g., an estimation task or a detection task~\cite{MattaSayedCoopGraphSP2018,SayedZhaoCoopGraphSP2018}) regarding a physical phenomenon of interest. 
This phenomenon is modeled through a spatially and temporally i.i.d. sequence of zero-mean and unit-variance streaming data, $\{\bm{x}_i(n)\}$, with $n=1,2,\ldots$, and $i=1,2,\ldots,N$.
In order to accomplish the inferential task, the agents organize themselves into a network, where each individual agent can rely on {\em local cooperation} with its neighbors.
Moreover, the network agents are required to be responsive to drifts in the observed phenomenon, and, to this end, they implement an {\em adaptive diffusion} strategy~\cite{Sayed, MattaBracaMaranoSayedTIT2016}. 
One useful implementation is the so-called combine-then-adapt (CTA) strategy, which has been studied in some revealing detail in the aforementioned references. 

The CTA rule consists of a combination step followed by an adaptation step.
During the combination step, agent $i$ {\em aggregates} data from its neighbors, by scaling these data with convex (i.e., nonnegative and adding up to one) combination weights $w_{i\ell}$, for $\ell=1,2,\dots,N$. The combination step produces the intermediate variable:
\beq
\bm{v}_i(n-1)=\sum_{\ell=1}^N w_{i\ell}\,\bm{y}_{\ell}(n-1).
\label{eq:combine}
\eeq
Subsequently, during the {\em adaptation} step, agent $i$ adjusts its output variable by comparing $\bm{v}_i(n-1)$ against the fresh streaming data $\bm{x}_i(n)$, and updating its state by using a small step-size $\mu\in(0,1)$:
\beq
\bm{y}_i(n)=\bm{v}_i(n-1) + \mu[\bm{x}_i(n) - \bm{v}_i(n-1)].
\label{eq:adapt}
\eeq
Merging~(\ref{eq:combine}) and~(\ref{eq:adapt}) into a single step yields:
\beq
\bm{y}_i(n)=(1-\mu)\sum_{\ell=1}^N w_{i\ell}\, \bm{y}_{\ell}(n-1) + \mu\,\bm{x}_i(n),
\label{eq:CTA1}
\eeq
which is equivalent to~(\ref{eq:1storder}), provided that one introduces a (scaled) combination matrix, $A$, with entries:
\beq
a_{ij}\dfz(1-\mu)w_{ij},
\label{eq:aweightsmatdef}
\eeq
and that the standard deviation $\sigma$ multiplying $\bm{x}_i(n)$ in~(\ref{eq:1storder}) is set equal to $\mu$.
For later use, it is useful to remark that, under this adaptive diffusion framework, the matrix $A$ is naturally nonnegative and its normalized counterpart, $A/(1-\mu)$, is a {\em right-stochastic matrix} (i.e., its rows add up to one). 

\section{Random Graph Model}
\label{sec:ERgraphs}
In the following, we shall denote by $\bm{G}$ the adjacency matrix of the network graph, whose entry $\bm{g}_{ij}$ is equal to one if nodes $i$ and $j$ are connected, and is equal to zero otherwise. The bold notation is used because we deal with {\em random} graphs. 

In this article we address the useful case where the network graph is generated according to the Erd\H{o}s-R\'enyi {\em random graph} model, namely, an undirected graph whose edges are drawn, one independently from the other, through a sequence of Bernoulli experiments with identical success (i.e., connection) probability~\cite{erd,BollobasRandom}.
In particular, the notation $\bm{G}\sim\mathscr{G}(N,p_N)$ will signify that the off-diagonal entries of the adjacency matrix $\bm{G}$ originate from an Erd\H{o}s-R\'enyi graph over $N$ nodes, and with connection probability $p_N$. 
Accordingly, the variables $\bm{g}_{ij}$, for $i=1,2,\dots, N$ and $j>i$, are independent Bernoulli random variables with $\P[\bm{g}_{ij}=1]=p_N$, and the matrix $\bm{G}$ is symmetric.
As it will be clear soon, the explicit dependence of the connection probability upon $N$ will be critical to examine the evolution of random graphs in the thermodynamic limit of large $N$. 

As one fundamental graph descriptor, in this work we use the {\em degree} of a node. The degree of node $i$ is defined as: 
\beq
\bm{d}_i=1 + \sum_{\ell\neq i} \bm{g}_{i\ell}=1 + \sum_{\ell\neq i} \bm{g}_{\ell i},
\label{eq:degdef}
\eeq
namely, the cardinality of the $i$-th node neighborhood (including $i$ itself). 
In particular, we shall use the minimal and maximal degrees that are defined as, respectively:
\beq
\bm{d}_{\min}\dfz \displaystyle{\min_{i=1,2,\dots N} \bm{d}_i},\qquad
\bm{d}_{\max}\dfz \displaystyle{\max_{i=1,2,\dots N} \bm{d}_i}.
\eeq

\subsection{Thermodynamic Limit of Random Graphs}
One meaningful (and classic) way to characterize the behavior of random graphs is to examine their thermodynamic limit as the network size goes to infinity.
Such an asymptotic characterization is useful because it captures average behavior that emerges with high probability over large networks. 
 
In examining the thermodynamic behavior of random graphs, the connection probability $p_N$ is generally allowed to scale with $N$. This degree of freedom allows representing different types of asymptotic graph behavior. For example, recalling that the average number of neighbors over an Erd\H{o}s-R\'enyi graph scales as $N p_N$, different graph evolutions can be obtained with different choices of $p_N$. For example, a constant $p_N$ will let the number of neighbors grow linearly with $N$. In comparison, a $p_N$ scaling as $(\log N)/N$ would correspond to a number of neighbors growing logarithmically with $N$.
In summary, different limiting regimes are determined by the way the connection probability evolves with $N$. 
It is useful for our purposes to list briefly the main regimes that are of interest for the forthcoming treatment.

--- {\em Connected regime}. In this work we focus on the regime where the graph is connected with high probability. This regime prescribes that the pairwise connection probability scales as~\cite{erd,BollobasRandom}:
\beq
p_N=\frac{\log N + c_N}{N}, \quad c_N\stackrel{N\rightarrow\infty}{\longrightarrow} \infty.
\label{eq:pNconn}
\eeq

--- {\em Sparse (connected) regime}. 
The connected regime can be obtained also when the pairwise connection probability, $p_N$, vanishes as $N$ gets large. In particular, we shall refer to this scenario as the {\em sparse} connected regime:
\beq
p_N\stackrel{N\rightarrow\infty}{\longrightarrow} 0\quad \textnormal{under (\ref{eq:pNconn})}
\quad
\textnormal{[Sparse connection regime]}.
\eeq

--- {\em Dense regime}. 
We call {\em dense} the regime where the pairwise connection probability converges to a nonzero quantity, namely $p_N\rightarrow p>0$. 

The aforementioned taxonomy basically focuses on the concepts of connectedness and sparsity. 
These concepts have been advocated in previous works related to topology inference with partial observations, and, in particular, some useful structural consistency results have been proved under the sparse (connected) regime. 

One essential element of novelty in our analysis is exploiting a different feature, namely, the {\em concentration} of graph degrees. 
We wish to avoid confusion here: the term ``concentration'' does {\em not} refer to the number of node connections.  
Instead, the concept of concentration is borrowed from a common terminology in statistics, which is used to refer to statistical quantities that concentrate around some (deterministic) value as $N\rightarrow\infty$~\cite{ConcIneqBook}. 
In particular, we will focus our attention on the {\em uniform concentration properties of the minimal and maximal degrees of random graphs}.

\begin{table}[t]
  \begin{center}
    \begin{tabular}{l|c|c} 
      Connection probability & Concentration & Sparsity\\
      \hline
      \hline
      &&
      \\
      $p_N=\omega_N\displaystyle{\frac{\log N}{N}}\rightarrow p>0$ & Uniform & Dense\\
      &&
      \\
      \hline
      &&
      \\
      $p_N=\omega_N\displaystyle{\frac{\log N}{N}}\rightarrow 0$, 	& Uniform & Sparse\\
      &&
      \\
      \hline
    \end{tabular}
    \vspace*{5pt}
    \caption{Useful taxonomy to illustrate the relationships between concentration and sparsity over a connected Erd\H{o}s-R\'enyi graph. The sequence $\omega_N$ goes to infinity as $N\rightarrow\infty$.}
      \label{tab:taxo}
  \end{center}
\end{table}

--- {\em Uniform concentration regime}. 
The uniform concentration regime is enabled by choosing the following pairwise connection probability:
\beq
p_N=\omega_N \,\frac{\log N}{N}\stackrel{N\rightarrow\infty}{\longrightarrow} p,\quad \omega_N \stackrel{N\rightarrow\infty}{\longrightarrow} \infty,
\label{eq:strongconc}
\eeq
which is tantamount to assuming that~(\ref{eq:pNconn}) holds true with the sequence $c_N$ growing faster than $\log N$. 
Under this regime, the minimal and the maximal degrees of the graph both concentrate asymptotically around the expected degree ($1+(N-1)p_N\sim N p_N$), in the following precise sense: 
\beq
\boxed{
\frac{\bm{d}_{\min}}{N p_N}\stackrel{\textnormal{p}}{\longrightarrow} 1,\quad
\frac{\bm{d}_{\max}}{N p_N}\stackrel{\textnormal{p}}{\longrightarrow} 1,\quad \textnormal{[Uniform concentration]}
}
\label{eq:dmaxmin0}
\eeq
The physical meaning of~(\ref{eq:dmaxmin0}) is that both the minimal and the maximal degrees scale, asymptotically with $N$, as the expected degree. Indeed, Eq.~(\ref{eq:dmaxmin0}) can be restated as: $\bm{d}_{\min}\sim N p_N + f_N$ and $\bm{d}_{\max}\sim N p_N + g_N$, where $f_N$ and $g_N$ are sequences that are asymptotically dominated by $N p_N$.

\begin{figure} [t]
\begin{center}
\includegraphics[scale= 0.3]{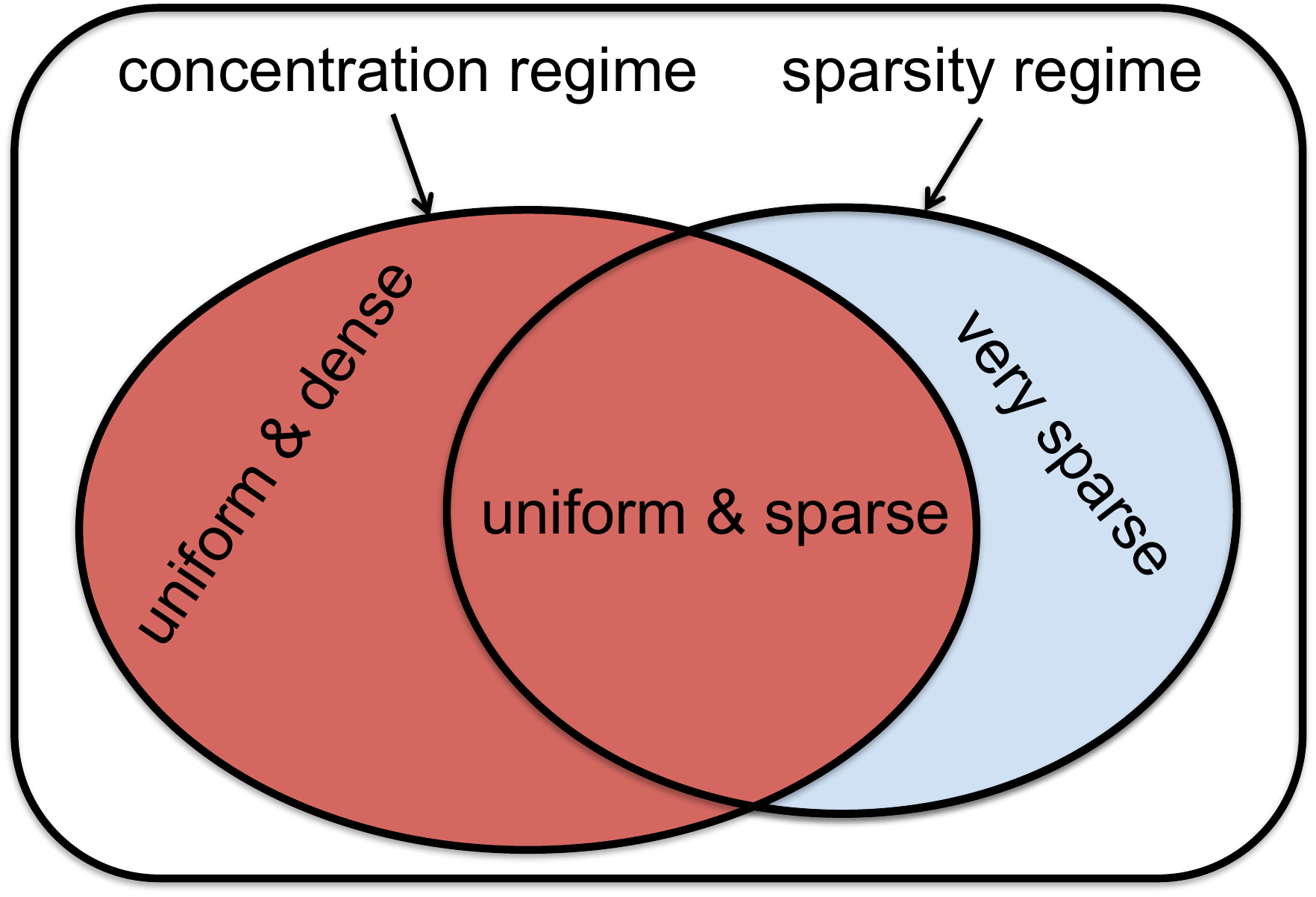}
\caption{
Venn diagram illustrating the relationships between concentration and sparsity over a connected Erd\H{o}s-R\'enyi graph. 
}
\label{fig:Venn}
\end{center}
\end{figure}

Table~\ref{tab:taxo} summarizes the sparsity/concentration taxonomy arising from the previous arguments. 
We are now ready to extract from the above taxonomy the elements that are relevant to the forthcoming treatment: 
\begin{enumerate}
\item Comparing~(\ref{eq:strongconc}) against~(\ref{eq:pNconn}), we see that the regime of concentration does {\em not} include all classes of connected Erd\H{o}s-R\'enyi graphs. In fact, while in~(\ref{eq:pNconn}) $c_N$ is any arbitrary divergent sequence (e.g., we can have $c_N=\log\log N$), according to~(\ref{eq:strongconc}) the sequence $c_N$ should grow with $N$ more than logarithmically. The regime where the graph is connected, whereas~(\ref{eq:strongconc}) is not fulfilled, will be referred to as the {\em very sparse} regime.  
\item According to~(\ref{eq:strongconc}), the regime of concentration can be either sparse or dense. In particular, the regime is dense when $p>0$, and is sparse when $p=0$.
\end{enumerate}
The aforementioned categorizations are illustrated in Fig.~\ref{fig:Venn} by means of a Venn diagram.

After having introduced the necessary details for the generation of the whole network, we are ready to introduce the formal model for partial observations.  

\begin{definition}[Partial observation setting]
In this work we adopt the partial observation setting that has been introduced in~\cite{tomo}: the {\em entire} network graph is generated according to the Erd\H{o}s-R\'enyi model, and the subnetwork of observable measurements, $\mathcal{S}$, has a cardinality $S$ scaling as:
\beq
\frac{S}{N}\stackrel{N\rightarrow\infty}{\longrightarrow} \xi\in (0,1),
\label{eq:cardscaling}
\eeq 
which means that $\xi$ is the (asymptotic) fraction of monitored nodes. 
Since $\xi$ is strictly less than one, condition~(\ref{eq:cardscaling}) conforms to a {\em partial observation} setting.
\end{definition}

\subsection{From Graphs to Combination Matrices}
Once a realization of the random graph occurs, a {\em combination matrix} is obtained by applying a certain {\em combination rule or policy} to this graph. The combination policy defines how the weights $a_{ij}$ are assigned given the particular graph structure. Formally, we have that:
\beq
\bm{A}=\pi(\bm{G}),
\label{eq:combpolicy}
\eeq
where $\pi: \{0,1\}^{N\times N}\rightarrow \mathbb{R}^{N\times N}$ is a deterministic combination policy and the randomness of $\bm{A}$ arises from the randomness of graph $\bm{G}$.

We start by introducing a useful class of combination matrices.
\begin{assumption}[Regular diffusion matrices]
\label{assum:regdiffmat}
We assume here that the combination matrix $\bm{A}$ is symmetric and that: 
\beq
\boxed{
\sum_{\ell=1}^N \bm{a}_{i\ell}=\rho,~~~~~~~
\frac{\kappa}{\bm{d}_{\max}}\,\bm{g}_{i j}\leq
\bm{a}_{i j}\leq \frac{\kappa}{\bm{d}_{\min}}\,\bm{g}_{i j}~~\forall i\neq j
}
\label{eq:fundamentalassum}
\eeq
for some parameters $\rho$ and $\kappa$, with $0<\kappa\leq \rho<1$.~\hfill$\square$
\end{assumption}
For the sake of concreteness, in the following we shall refer to combination matrices obeying the aforementioned assumption as {\em regular diffusion matrices} with parameters $\rho$ and $\kappa$.
We must remark that the most common combination matrices encountered in the literature automatically satisfy Assumption~\ref{assum:regdiffmat}.
For example, some popular choices are the Laplacian and the Metropolis rules reported below, which arise naturally in many applications, for instance, they are one fundamental ingredient of {\em adaptive} networks~\cite{MattaSayedCoopGraphSP2018}. 
The matrix entries corresponding to these combination rules are defined as follows. 
For $i\neq j$:
\beqa
&&\bm{a}_{ij}=\rho \lambda \, \displaystyle{\frac{\bm{g}_{ij}}{\bm{d}_{\max}}},~~~~~~~~~~~~~~~~\textnormal{[Laplacian rule]}
\label{eq:LapMat}
\\
&&\bm{a}_{ij}=\rho\,\displaystyle{\frac{\bm{g}_{ij}}{\max\left\{\bm{d}_i,\bm{d}_j\right\}}},~~~~~~~~\textnormal{[Metropolis rule]}
\label{eq:MetroMat}
\eeqa
whereas the self-weights are determined by the leftmost condition in~(\ref{eq:fundamentalassum}), yielding $\bm{a}_{ii}=\rho - \sum_{\ell\neq i} a_{i\ell}$. For the Laplacian rule, the parameter $\lambda$ fulfills the inequalities $0<\lambda\leq 1$.

\section{Consistent Graph Learning}
\label{sec:ident}
Before continuing, it is useful to summarize the probabilistic model considered in this work. 
First, the network graph $\bm{G}$ is generated according to an Erd\H{o}s-R\'enyi model. 
Then, a combination matrix $\bm{A}$ is generated according to a deterministic combination policy --- see~(\ref{eq:combpolicy}). 
Once a realization $\bm{A}=A$ is given, the stochastic dynamical system in~(\ref{eq:VARmodel}) evolves according to the randomness of the input source $\{\bm{x}_n\}$, which is statistically independent from the random graph $\bm{G}$.  
In the jargon of statistical physics, we could say that $\bm{A}$ is a ``quenched'' variable, i.e., once it is realized, it is frozen as the stochastic dynamical system evolves. 

Let us consider a stream of $n$ consecutive observations taken over the probed subset of nodes $\mathcal{S}$, and collected into the $S\times n$ matrix $\bm{Y}_n$ whose $(\ell,i)$-th entry is, for $\ell\in\mathcal{S}$ and $i=1,2,\ldots,n$:
\beq
[\bm{Y}_n]_{\ell i}=\bm{y}_{\ell}(i).
\label{eq:Ystream}
\eeq
A matrix estimator will be formally defined as some (measurable) function of the data $\widehat{\bm{A}}_{\mathcal{S},n}={\sf f}_n(\bm{Y}_n)$, namely,
\beq
{\sf f}_n: \mathbb{R}^{S\times n}\rightarrow \mathbb{R}^{S\times S}.
\eeq
We focus on the class of asymptotically stable estimators that converge as the number of samples increases. 
In particular, we consider the class of estimators that, for any realization of the combination matrix $\bm{A}$, guarantee the following convergence in probability: 
\beq
\lim_{n\rightarrow\infty}
\P[\|\widehat{\bm{A}}_{\mathcal{S},n}-{\sf h}\|_{\max}>\epsilon|\bm{A}=A]=0.
\label{eq:asystablest}
\eeq
We remark that the limiting estimator in~(\ref{eq:asystablest}), ${\sf h}$, is a deterministic quantity, given $\bm{A}=A$. 
However, this limit will be in general different for the $2^{N(N-1)/2}$ possible realizations of the random Erd\H{o}s-R\'enyi graph, i.e., we should write ${\sf h}={\sf h}(A)$ in~(\ref{eq:asystablest}).
In the following, we will use the following notation
\beq
\widehat{\bm{A}}_{\mathcal{S}}={\sf h}(\bm{A}),
\label{eq:Aslimdef0}
\eeq
where we suppressed the sample subscript $n$ to denote the {\em limiting} matrix estimator $\widehat{\bm{A}}_{\mathcal{S}}$. 
Using~(\ref{eq:Aslimdef0}) in~(\ref{eq:asystablest}) we have:
\beq
\boxed{
\lim_{n\rightarrow\infty}
\P[\|\widehat{\bm{A}}_{\mathcal{S},n}-\widehat{\bm{A}}_{\mathcal{S}}\|_{\max}>\epsilon]=0
}
\label{eq:asystableclass}
\eeq
Owing to dependence on $\bm{A}$, the matrix estimator $\widehat{\bm{A}}_{\mathcal{S}}$ is still random, and its randomness is determined solely by the randomness of the underlying graph. 
Our main goal is to establish that, in the commonly adopted {\em doubly-asymptotic} framework~\cite{ChandrasekaranParriloWillskyAOS2012,AnandkumarValluvanAOS2013} where the network size $N$ becomes large and the number of samples $n$ increases with $N$, it is possible to retrieve consistently (i.e., with probability converging to $1$) the true graph of the probed subset of nodes. In order to achieve this goal, in the next section we start by examining the asymptotic properties of the limiting estimator $\widehat{\bm{A}}_{\mathcal{S}}$ {\em as the network size $N$ goes to infinity}.

\subsection{Analysis of the Limiting Estimator as $N\rightarrow\infty$}
In order to ascertain whether or not it is possible to discriminate interacting (i.e., connected) agents from non-interacting agents, via observation of their output measurements, we now introduce the concept of margins and identifiability gap. 

\begin{figure} [t]
\begin{center}
\includegraphics[scale= 0.3]{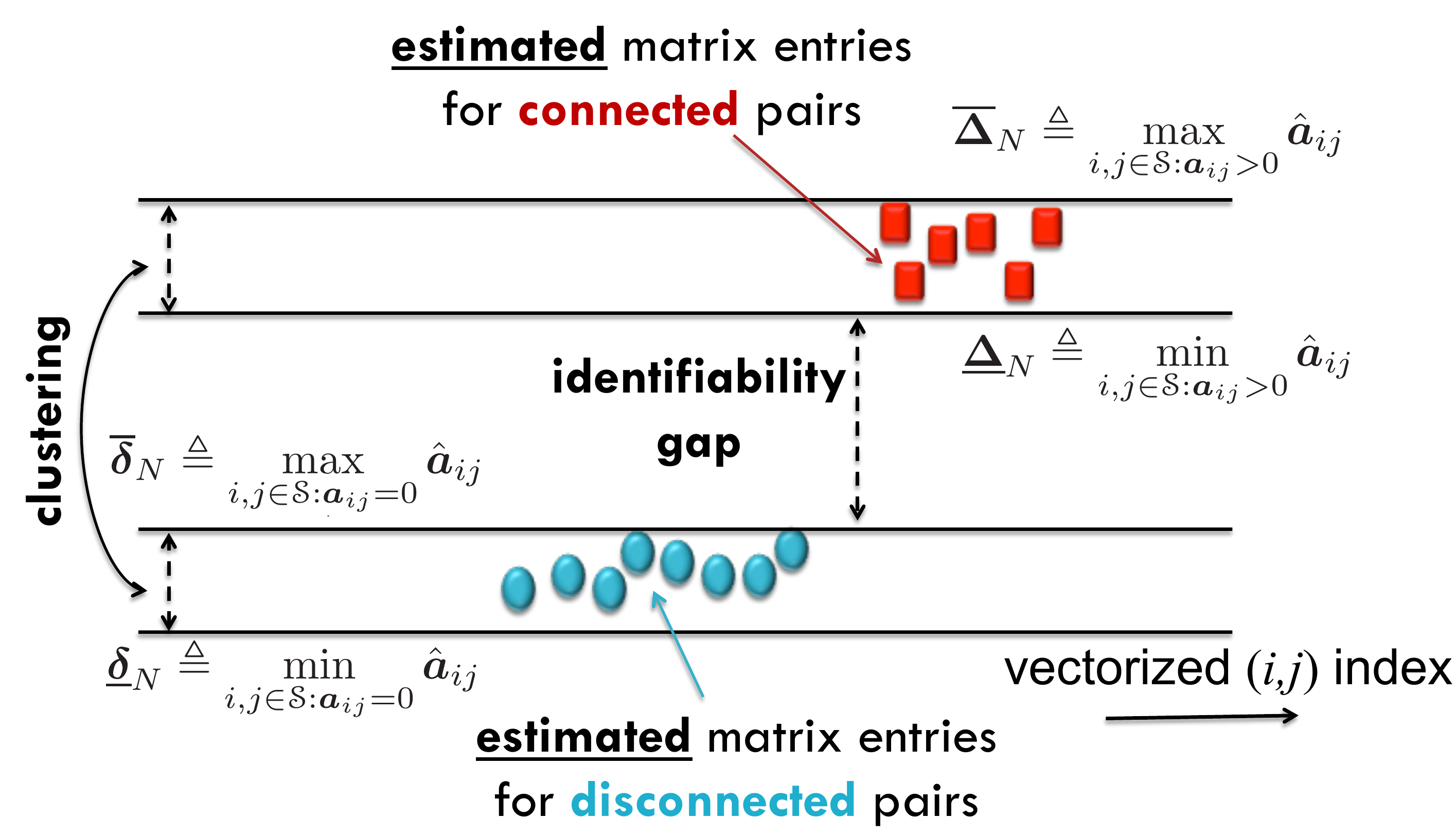}
\caption{Emergence of the identifiability gap.
}
\label{fig:identifiability}
\end{center}
\end{figure}

\begin{definition}[Margins]
Let $\widehat{\bm{A}}_{\mathcal{S}}$ be a certain estimated combination matrix, corresponding to the subset $\mathcal{S}$. 
The lower and upper margins corresponding to the {\em disconnected} pairs are defined as, respectively:\footnote{The definitions in~(\ref{eq:upperlowerdisc}) and~(\ref{eq:upperlowerconn}) are void if the nodes in $\mathcal{S}$ are all connected or all disconnected, respectively. Since, under the Erd\H{o}s-R\'enyi model, these events are irrelevant as $N\rightarrow\infty$, for these singular cases we can formally assign arbitrary values to the margins.}
\beq
\underline{\bm{\delta}}_N
\dfz
\displaystyle{\min_{\substack{i,j\in\mathcal{S}: \bm{a}_{ij}=0 \\ i\neq j}}\widehat{\bm{a}}_{ij}},\quad
\overline{\bm{\delta}}_N
\dfz
\displaystyle{\max_{\substack{i,j\in\mathcal{S}: \bm{a}_{ij}=0 \\ i\neq j}}\widehat{\bm{a}}_{ij}}.
\label{eq:upperlowerdisc}
\eeq
Likewise, the lower and upper margins corresponding to the {\em connected} pairs are defined as, respectively:
\beq
\underline{\bm{\Delta}}_N
\dfz
\displaystyle{\min_{\substack{i,j\in\mathcal{S}: \bm{a}_{ij}>0 \\ i\neq j}}\widehat{\bm{a}}_{ij}},\quad
\overline{\bm{\Delta}}_N
\dfz
\displaystyle{\max_{\substack{i,j\in\mathcal{S}: \bm{a}_{ij}>0 \\ i\neq j}}\widehat{\bm{a}}_{ij}}.
\label{eq:upperlowerconn}
\eeq
~\hfill$\square$
\end{definition}
The aforementioned margins are useful to examine the achievability of structural consistency for an estimator $\widehat{\bm{A}}_{\mathcal{S}}$ --- see Fig.~\ref{fig:identifiability} for an illustration --- and lead to the concept of {\em identifiability gap}.

\begin{definition}[Local structural consistency]
Let $\widehat{\bm{A}}_{\mathcal{S}}$ be an estimated combination matrix. If there exists a sequence $s_N$, a real value $\eta$, and a strictly positive value $\Gamma$, such that, for all $\epsilon>0$:
\beqa
&&\lim_{N\rightarrow\infty}\P[s_N \, \overline{\bm{\delta}}_N < \eta+\epsilon]=1,
\nonumber\\
&&\lim_{N\rightarrow\infty}\P[s_N \,\underline{\bm{\Delta}}_N > \eta+ \Gamma - \epsilon]=1,
\label{eq:igap}
\eeqa
we say that $\widehat{\bm{A}}_{\mathcal{S}}$ achieves local structural consistency, with a bias at most equal to $\eta$, an identifiability gap at least equal to $\Gamma$, and with a scaling sequence $s_N$.~\hfill$\square$ 
\end{definition}

\begin{remark}[Locality]
We use the qualification ``local'' to emphasize that the structure of the subnetwork $\mathcal{S}$ must be inferred from observations gathered {\em locally} in $\mathcal{S}$, even if the nodes of $\mathcal{S}$ undergo the influence of many other nodes belonging to the larger embedding network.~\hfill$\square$
\end{remark}

\begin{remark}[Identifiability gap] 
The condition $s_N \, \overline{\bm{\delta}}_N < \eta+\epsilon$ means that the {\em maximum} entry of $s_N \widehat{\bm{A}}_{\mathcal{S}}$ taken over the {\em disconnected} pairs essentially does not exceed $\eta$. 
Likewise, the condition $s_N \,\underline{\bm{\Delta}}_N > \eta+ \Gamma - \epsilon$ means that the {\em minimum} entry of $s_N \widehat{\bm{A}}_{\mathcal{S}}$ taken over the {\em connected} pairs essentially stays above the value $\eta+\Gamma>\eta$. Combining these two relationships, we conclude that the estimated matrix entries corresponding to connected node pairs stand clearly separated from the entries corresponding to disconnected node pairs. 
The minimum amount of separation is quantified by the gap, $\Gamma$.~\hfill$\square$
\end{remark}

\begin{remark}[Bias]
For the {\em true} combination matrix, the entries corresponding to disconnected pairs are zero. 
In contrast, Eq.~(\ref{eq:igap}) reveals that the scaled entries for disconnected pairs can be close to $\eta$, which results therefore in a {\em bias}. However, and remarkably, {\em this bias does not constitute a problem for consistent classification of connected/non-connected nodes}, i.e., the bias does not affect in any manner identifiability.~\hfill$\square$
\end{remark}

The notion of structural consistency implies the existence of a threshold, comprised between $\eta$ and $\eta+\Gamma$, which correctly separates (in the limit of large $N$) the entries of the matrix estimator, in such a way that the entries corresponding to connected pairs lie above the threshold, whereas the entries corresponding to disconnected pairs lie below the threshold. 

However, an accurate determination of the separation threshold requires some prior knowledge of the monitored system. 
For instance, to set a detection threshold one needs to know the scaling sequence $s_N$. 
In the problems dealt with in this work we will see that $s_N=N p_N$, where $N p_N$ represents the average number of neighbors in the network, and in several practical applications this number is unknown.

As a result, the threshold setting might be a critical issue, and it would be highly desirable to have a {\em universal} (i.e., blind and nonparametric) method to set the threshold. For example, it would be highly desirable to determine a separation threshold using machine learning tools such as a standard clustering algorithm (e.g., a $k$-means clustering with $k=2$).  
It is therefore useful to strengthen the notion of structural consistency so as to incorporate the aforementioned requirement of universality.

\begin{definition}[Universal local structural consistency]
Let $\widehat{\bm{A}}_{\mathcal{S}}$ be an estimated combination matrix. If there exists a sequence $s_N$, a real value $\eta$, and a strictly positive value $\Gamma$, such that, for any $\epsilon>0$:
\beqa
&&\lim_{N\rightarrow\infty}\P[|s_N \,\underline{\bm{\delta}}_N -\eta|<\epsilon]=1, 
\nonumber\\
&&\lim_{N\rightarrow\infty}\P[|s_N \,\overline{\bm{\delta}}_N -\eta|<\epsilon]=1, 
\nonumber\\
&&\lim_{N\rightarrow\infty}
\P[|s_N \,\underline{\bm{\Delta}}_N - (\eta + \Gamma)|<\epsilon]=1,
\nonumber\\
&&\lim_{N\rightarrow\infty}
\P[|s_N \,\overline{\bm{\Delta}}_N - (\eta + \Gamma)|<\epsilon]=1,
\label{eq:igapstrong}
\eeqa
we say that $\widehat{\bm{A}}_{\mathcal{S}}$ achieves uniform local structural consistency, with a bias $\eta$, an identifiability gap $\Gamma$, and with a scaling sequence $s_N$.~\hfill$\square$
\end{definition}

\begin{remark}[Clustering]
According to the notion of {\em universal} structural consistency, the pair of (scaled) lower margins, $s_N\,\underline{\bm{\delta}}_N$ and $s_N\,\overline{\bm{\delta}}_N$, converge to one and the same value, $\eta$, which implies that {\em all the entries of} $s_N \widehat{\bm{A}}_{\mathcal{S}}$ {\em corresponding to disconnected pairs} are sandwiched between these margins --- see Fig.~\ref{fig:identifiability}. 
A similar behavior is observed for the scaled entries over the connected pairs, which converge altogether to $\eta+\Gamma$ since they are sandwiched between $s_N\,\underline{\bm{\Delta}}_N$ and $s_N\,\overline{\bm{\Delta}}_N$. In summary, we conclude that the connected and disconnected agent pairs {\em cluster into well-separated classes} that can be identified, e.g., by means of a universal clustering algorithm.~\hfill$\square$
\end{remark}

\begin{remark}[Scale-irrelevance]
\label{rem:scaleirrel}
The definition of identifiability and consistency do not clarify what would be the right value for the identifiability gap. For example, assume that a certain estimator $\widehat{\bm{A}}_{\mathcal{S}}$ leads to a gap $\Gamma$. 
If we now scale $\widehat{\bm{A}}_{\mathcal{S}}$ by some (positive) constant $c$ to construct another estimator $c\,\widehat{\bm{A}}_{\mathcal{S}}$, then the gap becomes $c\,\Gamma$. This would imply that choosing an arbitrarily large value for $c$ would lead to an arbitrarily large gap, what does this mean?

This observation provides an opportunity to clarify the meaning of the identifiability gap. 
The concept behind identifiability is that there exists a strictly positive $\Gamma$ that allows separating the two sets of connected and disconnected agent pairs. {\em Asymptotically}, the existence of a strictly positive $\Gamma$ enables successful discrimination. 
However, for finite network sizes, the value itself of the gap is not informative about the performance achievable by a classifier that operates on the entries of the estimated matrix. What really determines the classifier performance is the {\em spread of the entries relative to the gap value}. 
A high (resp. low) variability makes classification more (resp. less) difficult. 
For similar reasons, clustering algorithms are clearly unaffected by changes in scale, while they are sensitive to the variability of the entries within each individual cluster.~\hfill$\square$ 
\end{remark}

\subsection{A Consistent Clustering Algorithm}
\label{subsec:clualgo}
The definition of universal local structural consistency implies that the disconnected node pairs cluster around $\eta$, whereas the connected node pairs cluster around the higher value $\eta+\Gamma$. 
Accordingly, {\em for sufficiently large $N$}, there is no doubt that any reasonable clustering algorithm will be able to identify properly these two clusters. For example, a correct {\em asymptotic} guess of the true clusters (i.e., of the true graph) can be obtained by simply choosing an intermediate threshold between the maximum and the minimum matrix entries. 

On the other hand, since in practice we work with {\em finite} network sizes, a number of non-ideal effects could appear, and the non-perfect localization of the clusters could affect the inference performance. 
For these reasons, we start by focusing on a popular and general-purpose clustering algorithm, namely, the $k$-means algorithm (in our case, we know that $k=2$), because it typically offers good performance under many operational conditions. 
However, and somehow paradoxically, while the $k$-means works properly for finite network sizes, it is not obvious that it should be {\em asymptotically} consistent. This is because the $k$-means might suffer when dealing with clusters of very different sizes. Since in our model it is actually permitted that, for large $N$, one cluster can dominate the other one (for instance, when $p_N\rightarrow 0$, the cluster of disconnected nodes becomes predominant), we are not guaranteed that the $k$-means algorithm works properly as $N\rightarrow\infty$.

In order to circumvent this problem, we propose a slight modification of the $k$-means algorithm that is specifically tailored to the universal local structural consistency property. Preliminarily, let $v$ be the $L\times 1$ vector to be clustered, with entries that have been arranged in ascending order.  
The $k$-means algorithm, with $k=2$, attempts to minimize the following cost function:
\beq
\sum_{v_j\in\mathcal{C}_0}(v_j - c_0)^2 + \sum_{v_j\in\mathcal{C}_1}(v_j - c_1)^2,
\label{eq:kmcostfun}
\eeq 
over all possible clusters $\mathcal{C}_0$ and $\mathcal{C}_1$, with $c_0$ and $c_1$ being the cluster centroids, defined as:
\beq
c_0=\frac{1}{|\mathcal{C}_0|}\sum_{v_j\in\mathcal{C}_0} v_j,\quad
c_1=\frac{1}{|\mathcal{C}_1|}\sum_{v_j\in\mathcal{C}_1} v_j.
\eeq 
It is useful to recall that the minimum of the cost function must fulfill the following necessary condition: the midpoint between the two centroids is a threshold that separates the two clusters. 
In our one-dimensional case, with $k=2$, this property implies that it suffices to consider only the cluster configurations $\mathcal{C}_0(j)=\{1,2,\ldots,j\}$ and $\mathcal{C}_1(j)=\{j+1,j+2,\ldots,L\}$, for $j\in\{1,2,\ldots,L-1\}$. (Obviously, if two points, say $v_{i}$ and $v_{i+1}$, are coincident, considering their possible permutations is pointless). Accordingly, we see that any possible partition is identified by an index $j$.

\begin{figure} [t]
\begin{center}
\includegraphics[scale= 0.27]{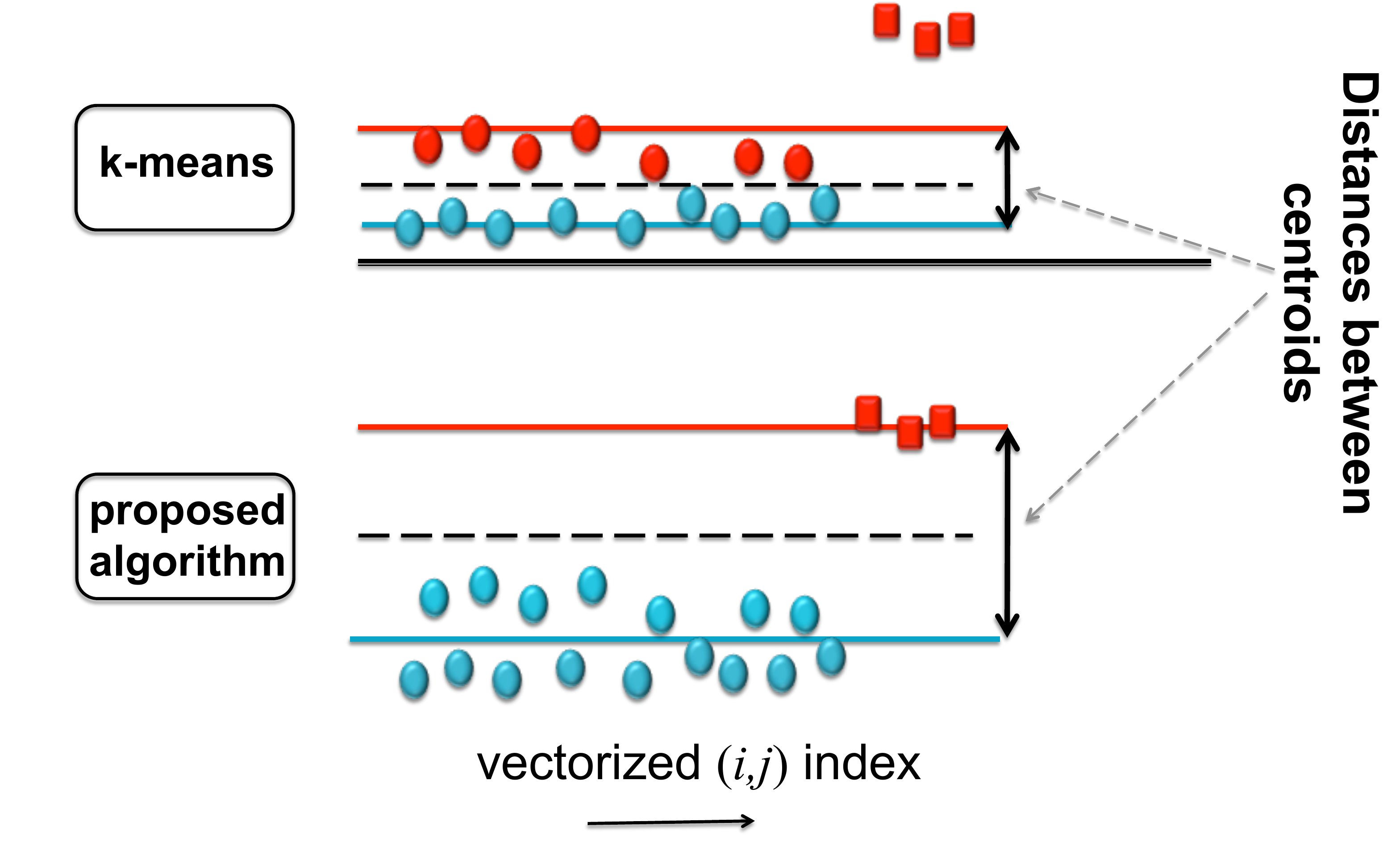}
\caption{Visual comparison between the $k$-means algorithm and the clustering algorithm proposed in this work, for the case of unbalanced clusters. The true clusters are identified by different symbols (circle vs. square). The clusters produced by the algorithms are identified by different colors (blue vs. red). 
}
\label{fig:clustexplain}
\end{center}
\end{figure}

Let us now illustrate a possible (and known) limitation of the $k$-means algorithm, in relation to our graph learning problem. 
For large $N$, the true clusters of connected and disconnected pairs become almost perfectly localized, and, hence, two centroids belonging to the true clusters match the necessary condition (lowermost panel in Fig.~\ref{fig:clustexplain}). 
However, if the sizes of the true clusters are very different, the effect of the smallest ensemble (squares) could be asymptotically negligible, resulting in the configuration illustrated in the uppermost panel of Fig.~\ref{fig:clustexplain}. Here the two centroids estimated by the $k$-means algorithm are both located within the largest ensemble (circles), leading to a wrong partitioning. 
For this reason, we now introduce a modification of the $k$-means algorithm that takes into account also the distance between the centroids. 
Let us now describe the clustering algorithm we propose. 
First, the algorithm enumerates all admissible clusters through an index $j$ spanning the set $\{1,2,\ldots,L-1\}$. 
The set of indices fulfilling the necessary condition of $k$-means are collected in the set $\mathcal{A}=\{j_1,j_2,\ldots\}$. 
At this stage, the classic $k$-means would simply select, among these admissible points, the one ensuring the minimum cost. 
We modify this rule by selecting the index $j^{\star}\in\mathcal{A}$ that maximizes the distance between the clusters' centroids, namely, $j^{\star}=\arg\!\min_{j\in\mathcal{A}} [c_1(j) - c_0(j)]$, with $c_0(j)$ and $c_1(j)$ being the centroids corresponding to the clusters identified by index $j$.
With this modified rule, we want to $i)$ retain the good behavior exhibited by $k$-means in typical situations; and $ii)$ guarantee that the algorithm achieves consistent graph learning. 
Point $i)$ will be confirmed by the numerical experiments in Sec.~\ref{sec:illustexam}. 
Point $ii)$ will be in fact formally established in the next theorem.

Before stating the theorem, let us introduce the input-output relationship pertaining to the clustering algorithm. 
Given an input matrix $M$ with index set $\mathcal{S}$, the clustering algorithm is applied to the off-diagonal entries $m_{ij}$ for $i\in\mathcal{S}$ and $i\neq j$. These entries are classified into classes $\mathcal{C}_0$ and $\mathcal{C}_1$ following the two-step procedure just described. 
Then the clustering algorithm produces as output an estimated adjacency matrix $\widehat{G}$ (with main diagonal set conventionally to $0$) such that, for all $i,j\in\mathcal{S}$ with $i\neq j$:
\beq
\widehat{G}={\sf clu}(M): \widehat{g}_{ij}=\mathbb{I}[m_{i,j}\in\mathcal{C}_1],
\eeq  
where $\mathbb{I}[\mathcal{E}]$ denotes the indicator function, which is equal to $1$ if condition $\mathcal{E}$ is true, and is equal to $0$ otherwise.

\begin{theorem}[Sample consistency of the proposed clustering algorithm]
\label{theor:samplecons}
Let $\widehat{\bm{A}}_{\mathcal{S},n}$ be a stable matrix estimator belonging to class~(\ref{eq:asystableclass}), with the limiting matrix estimator $\widehat{\bm{A}}_{\mathcal{S}}$ achieving universal local structural consistency. 
Let $\bm{G}_{\mathcal{S}}$ be the (random) support graph associated to $\bm{A}_{\mathcal{S}}$ and let
\beq
\widehat{\bm{G}}_{\mathcal{S},n}={\sf clu}(\widehat{\bm{A}}_{\mathcal{S},n})
\eeq
be the subgraph estimated by the proposed clustering algorithm.
Then, a certain scaling law $n(N)$ exists such that:
\beq
\boxed{
\lim_{N\rightarrow\infty}\P\left[\widehat{\bm{G}}_{\mathcal{S},n(N)} = \bm{G}_{\mathcal{S}}\right]=1
}
\eeq
\end{theorem}
\begin{IEEEproof}
See Appendix~\ref{app:samplecons}.
\end{IEEEproof}

\section{Proposed Matrix Estimators}
\label{sec:ME}
We now introduce the three estimators examined in this work. 
Preliminarily, it is useful to introduce the steady-state correlation matrix corresponding to the dynamics in~(\ref{eq:VARmodel}).  
\beq
\bm{R}_0=\lim_{n\rightarrow \infty}\mathbb{E}\left[\bm{y}_n\bm{y}_n^{\top}\right],
\label{eq:limR0def0}
\eeq 
whose existence will be shown soon. 
Exploiting~(\ref{eq:VARmodel}) we see that:
\beq
\mathbb{E}\left[\bm{y}_n\bm{y}_n^{\top}\right]=\sigma^2 \sum_{i=0}^{n-1} \bm{A}^i (\bm{A}^i)^T.
\eeq
Using now the symmetry of $\bm{A}$ along with~(\ref{eq:limR0def0}) we get:
\beq
\bm{R}_0=\sigma^2 \sum_{i=0}^{\infty} \bm{A}^{2i}=\sigma^2 ( I - \bm{A}^2 )^{-1},
\label{eq:R0explicit0}
\eeq
where $I$ is the $N\times N$ identity matrix, and where the convergence of the series is guaranteed by the stability of $\bm{A}$. 
We note explicitly that the randomness of $\bm{R}_0$ (bold notation) comes from randomness of the underlying Erd\H{o}s-R\'enyi graph.

Likewise, we introduce the steady-state one-lag correlation matrix: 
\beq
\bm{R}_1=\lim_{n\rightarrow \infty}\mathbb{E}\left[\bm{y}_n\bm{y}_{n-1}^{\top}\right],
\eeq
which exploiting the dynamics in~(\ref{eq:VARmodel}) can be written as: 
\beq
\bm{R}_1=\bm{A} \bm{R}_0.
\label{eq:verypreliminaryGranger}
\eeq

\subsection{Granger Estimator}
From~(\ref{eq:verypreliminaryGranger}) we obtain the following well-known relationship:
\beq
\bm{A}=\bm{R}_1\bm{R}_0^{-1},
\label{eq:Granger0}
\eeq
a quantity that is also referred to as the best one-step predictor or Granger estimator~\cite{MaterassiSalapakaTAC2012,Geigeretal15}. 
Under the partial observability setting, it is tempting to adapt the structure in~(\ref{eq:Granger0}) by considering only the observable subnet $\mathcal{S}$, and by neglecting the effect of the latent nodes, yielding~\cite{tomo, tomo_journal}:
\beq
\widehat{\bm{A}}^{\textnormal{(Gra)}}_{\mathcal{S}}=[\bm{R}_1]_{\mathcal{S}} ([\bm{R}_0]_{\mathcal{S}})^{-1}.
\eeq
In~\cite{tomo} it is shown that such matrix estimator admits the following representation: 
\beq
\widehat{\bm{A}}^{\textnormal{(Gra)}}_{\mathcal{S}}=\bm{A}_{\mathcal{S}}+\bm{E}^{\textnormal{(Gra)}},
\label{eq:newwidehatAsGra}
\eeq
where (we recall that $\mathcal{S}'$ denotes the subset of unobserved nodes):
\beq
\bm{E}^{\textnormal{(Gra)}}=\bm{A}_{\mathcal{S}\mathcal{S}'} \bm{H} [\bm{A}^2]_{\mathcal{S}'\mathcal{S}}, 
\label{eq:errGraMatform}
\eeq
with:
\beq
\bm{H}=(I_{\mathcal{S}'} - \bm{C})^{-1},\qquad  \bm{C}\dfz[\bm{A}^2]_{\mathcal{S}'}.
\label{eq:Hmatdef}
\eeq
For later use, it is also useful to rewrite~(\ref{eq:errGraMatform}) on an entrywise basis, for all $i,j\in\mathcal{S}$:
\beq
\boxed{
\bm{e}^{\textnormal{(Gra)}}_{ij}=\sum_{\ell,m\in\mathcal{S}'} \bm{a}_{i\ell} \bm{h}_{\ell m} \bm{a}^{(2)}_{m j}
}
\label{eq:basicerrGra}
\eeq

\subsection{One-Lag Correlation Estimator}
In this section we use the one-lag correlation matrix as estimator for the combination matrix.\footnote{When $\sigma^2$ is known, the relationship between $\bm{R}_1=\sigma^2 \bm{A} ( I - \bm{A}^2 )^{-1}$ and $\bm{A}$ is invertible under the full observability regime. 
This can be easily shown by resorting to the spectral decomposition of $\bm{R}_1/\sigma^2$ and $\bm{A}$ (which share the same eigenvectors), and by finding the eigenvalues of $\bm{A}$ from those of $\bm{R}_1/\sigma^2$; this inversion operation can be 
successfully realized because $\bm{A}$ is symmetric and has spectral radius less than one. However, we remark that $\sigma^2$ is assumed unknown and, more importantly, that we work under a partial observability regime, and, hence, the aforementioned inversion procedure does not apply.} 
The reason behind such choice is the following series expansion of the one-lag correlation matrix:
\beqa
\bm{R}_1=\bm{A} \bm{R}_0&=&\sigma^2 \bm{A} ( I - \bm{A}^2 )^{-1}\nonumber\\
&=&
\sigma^2\left(\bm{A} + \bm{A}^3 +\bm{A}^5 +\ldots\right).
\label{eq:R1firstexpansion}
\eeqa
When applied only to the submatrix corresponding to $\mathcal{S}$, Eq.~(\ref{eq:R1firstexpansion}) yields:
\beq
\widehat{\bm{A}}^{\textnormal{(1-lag)}}_{\mathcal{S}}=
[\bm{R}_1]_{\mathcal{S}}=
\sigma^2\left(\bm{A}_{\mathcal{S}} + [\bm{A}^3]_{\mathcal{S}} +[\bm{A}^5]_{\mathcal{S}} +\ldots \right).
\label{eq:onelagestdef}
\eeq
It is convenient to rewrite~(\ref{eq:onelagestdef}) as:
\beq
\widehat{\bm{A}}^{\textnormal{(1-lag)}}_{\mathcal{S}}=\sigma^2\left(\bm{A}_{\mathcal{S}}+\bm{E}^{\textnormal{(1-lag)}}\right),
\label{eq:newwidehatAsOneLag}
\eeq
where:
\beq
\bm{E}^{\textnormal{(1-lag)}}=\sum_{h=1}^{\infty}[\bm{A}^{2 h + 1}]_{\mathcal{S}},
\label{eq:erronelag1}
\eeq
or, for all $i,j\in\mathcal{S}$:
\beq
\boxed{
\bm{e}^{\textnormal{(1-lag)}}_{ij}=\sum_{h=1}^{\infty} \bm{a}_{ij}^{(2 h + 1)}
}
\label{eq:onelagesterrdef}
\eeq

\subsection{Residual Estimator}
Let us introduce the residual vector that computes the (scaled) difference between consecutive time samples:
\beq
\bm{r}_n\dfz \frac{\bm{y}_n - \bm{y}_{n-1}}{\sqrt{2}}.
\eeq
We have clearly $\E[\bm{r}_n \bm{r}_n^{\top}]=\bm{R}_0 - \bm{R}_1=\sigma^2 (I + \bm{A})^{-1}$. 
Accordingly, it makes sense to introduce the following estimator:
\beqa
\widehat{\bm{A}}^{\textnormal{(res)}}_{\mathcal{S}}&=&[\bm{R}_1]_{\mathcal{S}} - [\bm{R}_0]_{\mathcal{S}}
=-\sigma^2\left[(I + \bm{A})^{-1} \right]_{\mathcal{S}}\nonumber\\
&=&
\sigma^2\left(\bm{A}_{\mathcal{S}} - I_{\mathcal{S}} - [\bm{A}^2]_{\mathcal{S}} + [\bm{A}^3]_{\mathcal{S}} + \ldots \right).
\label{eq:newResidualSeries0}
\eeqa
The structure of~(\ref{eq:newResidualSeries0}) motivates the introduction of the matrix:
\beq
\bm{E}^{\textnormal{(res)}}=- I_{\mathcal{S}} + \sum_{h=1}^{\infty}
([\bm{A}^{2h +1}]_{\mathcal{S}} - [\bm{A}^{2h}]_{\mathcal{S}}),
\label{eq:errres1}
\eeq
yielding:
\beq
\widehat{\bm{A}}^{\textnormal{(res)}}_{\mathcal{S}}=\sigma^2\left(\bm{A}_{\mathcal{S}}+\bm{E}^{\textnormal{(res)}}\right).
\label{eq:newwidehatAsRes}
\eeq
Equation~(\ref{eq:errres1}) implies that, for all $i,j\in\mathcal{S}$, with $i\neq j$:
\beq
\boxed{
\bm{e}^{\textnormal{(res)}}_{ij}=
\sum_{h=1}^{\infty}
\left(\bm{a}_{ij}^{(2h +1)} - \bm{a}_{ij}^{(2h)}\right)
}
\label{eq:residerrestdef}
\eeq

\begin{table*}[t]
  \begin{center}
    \begin{tabular}{l|c|c} 
      ~~~~~~~~~~Estimator &  Error bias $\eta$ & Identifiability gap $\Gamma$\\
      \hline
      \hline
      &&
      \\
      Granger: $[\bm{R}_1]_{\mathcal{S}} ([\bm{R}_0]_{\mathcal{S}})^{-1}$ & 
      $\displaystyle{\kappa^2 p \, \frac{(2\rho - \kappa)\,(1 - \xi)}{1 -  (\rho^2 - 2\rho\kappa\xi + \kappa^2 \xi)}}$ & 
      $\kappa$\\
      &&
      \\
      \hline
      &&
      \\
      one-lag: $[\bm{R}_1]_{\mathcal{S}}$ & 
      $\displaystyle{\sigma^2 \kappa^2 p\, \displaystyle{\frac{\rho + \rho \, (\rho-\kappa)^2 + 2\,(\rho-\kappa)}
      {(1-\rho^2)(1 - (\rho - \kappa)^2)^2}}}$ & 
      $\displaystyle{\frac{1 + (\rho-\kappa)^2}{(1 - (\rho-\kappa)^2)^2}}\times \sigma^2\kappa$\\
      &&
      \\
      \hline
      &&
      \\
      residual: $[\bm{R}_1]_{\mathcal{S}} - [\bm{R}_0]_{\mathcal{S}}$ & 
      $-  \displaystyle{\frac{\sigma^2 \kappa^2 p}{(1+\rho)(1 + \rho - \kappa)^2}}$ & 
      $\displaystyle{\frac{\sigma^2\kappa}{(1 + \rho - \kappa)^2}}$\\
      &&
      \\
      \hline
    \end{tabular}
    \vspace*{5pt}
    \caption{Theorem~\ref{theor1}: Biases and gaps of the estimators listed in the leftmost column. For all cases, the scaling sequence is $s_N=N p_N$.}
    \label{tab:Theorem1}
  \end{center}
\end{table*}

\begin{remark}[Error series]
One essential ingredient for the proof of Theorem~\ref{theor1} further ahead is characterizing the convergence properties of the error series in~(\ref{eq:basicerrGra}),~(\ref{eq:onelagesterrdef}) and~(\ref{eq:residerrestdef}). 
The different behavior of these series can give rise to different performance of the graph learners. 
For example, comparing~(\ref{eq:onelagesterrdef}) against~(\ref{eq:residerrestdef}), we see that, due to the subtraction of the correlation matrix $\bm{R}_0$, the error of the residual estimator takes the structure of an alternating series. This feature might be useful to reduce the error and, hence, to improve the performance of the graph learner. We will explore this aspect in Sec.~\ref{sec:illustexam}.~\hfill$\square$
\end{remark}

\section{Achievability of Universal Local Structural Consistency}
\label{sec:MainTh}
The following theorem establishes the fundamental consistency properties of the estimators presented in Sec.~\ref{sec:ME}.

\begin{theorem}[Universal local structural consistency]
\label{theor1}
Let $\bm{A}$ be a regular diffusion matrix with parameters $\rho$ and $\kappa$, with the network graph drawn according to an Erd\H{o}s-R\'enyi random graph model $\mathscr{G}(N,p_N)$ where the fraction of observable nodes, $S/N$, converges to some nonzero value $\xi$.

Then, under the uniform concentration regime where:
\beq
p_N=\omega_N\displaystyle{\frac{\log N}{N}}\rightarrow p,~~~~\textnormal{with } \omega_N\rightarrow \infty,
\eeq
the Granger, the one-lag, and the residual estimators achieve universal local structural consistency, with a scaling sequence $s_N=N p_N$, and with the error biases and identifiability gaps listed in Table~\ref{tab:Theorem1}. 
\end{theorem}
\begin{IEEEproof}
The proof of the theorem is provided in Appendix~\ref{app:Theorem1}, and relies heavily on a number of auxiliary lemmas and theorems reported in the appendices. In particular, the core of the proof is the following. 
First, Theorems~$2$ and~$3$ in Appendix~\ref{app:Matrecursions} construct {\em uniform} (w.r.t. $N$) bounds on the entries of the combination matrix and related matrices. These bounds are useful to characterize the error associated to the three estimators listed in Table~\ref{tab:Theorem1}. Then, exploiting the asymptotic concentration property of the maximal and minimal degrees, it is possible to prove the convergence of the matrix series relevant for the computation of the errors, namely, Eqs.~(\ref{eq:basicerrGra}),~(\ref{eq:onelagesterrdef}) and~(\ref{eq:residerrestdef}). 
These convergence properties are finally used to compute the bias and the identifiability gaps listed in Table~\ref{tab:Theorem1}. 
\end{IEEEproof}

\begin{remark}[Role of degree concentration]
For the class of regular diffusion matrices in Assumption~\ref{assum:regdiffmat}, concentration of the degrees induces concentration of the nonzero entries of the combination matrix. 
This clearly creates an identifiability gap in the {\em true} matrix $\bm{A}_{\mathcal{S}}$.
However, what is critical for graph recovery is the existence of an identifiability gap in an {\em estimated} matrix $\widehat{\bm{A}}_{\mathcal{S}}$. 
The existence of graph learners possessing such a gap is in fact proved in Theorem~\ref{theor1}. 
Let us provide some intuition behind this result.

To this aim, we start by examining the useful representations of the matrix estimators in Eqs.~(\ref{eq:newwidehatAsGra}),~(\ref{eq:newwidehatAsOneLag}), and~(\ref{eq:newwidehatAsRes}). 
From these representations, we see that the existence of an identifiability gap in a matrix estimator $\widehat{\bm{A}}_{\mathcal{S}}$ depends on $\bm{A}_{\mathcal{S}}$, but will depend strongly also on some error matrix, say it $\bm{E}$. 
Since each entry in $\bm{E}$ is a function of the entries in $\bm{A}$ (in general, also of the {\em latent} nodes belonging to the unobserved subset $\mathcal{S}'$), a key point is to understand how the $\bm{a}_{ij}$'s combine with each other to produce $\bm{E}$.
For what observed before, the $\bm{a}_{ij}$'s exhibit {\em concentration} (in their nonzero values). On the other hand, they exhibit also {\em randomness} (in the {\em location} of the nonzero entries, due to the random graph model). 
The attributes of concentration and randomness are critical to reveal the nontrivial result shown in Theorem~\ref{theor1}: For suitably structured error matrices (i.e., the matrices corresponding to the three aforementioned estimators), it will be seen that the $\bm{a}_{ij}$'s combine with each other so as to induce a concentration in $\bm{E}$, which in turn determines the emergence of an identifiability gap in $\widehat{\bm{A}}_{\mathcal{S}}$. 
In summary, using appropriate matrix estimators the overall influence of latent agents is quantified by an error matrix, whose entries converge to some deterministic quantity, equally for both connected and disconnected agents. 
In this way, the connections among the probed agents stick out consistently from the error floor.

Moreover, in light of Theorem~\ref{theor:samplecons}, the universal local structural consistency of the {\em limiting} estimators $\widehat{\bm{A}}_{\mathcal{S}}$ implies consistency of the matrix estimators $\widehat{\bm{A}}_{\mathcal{S},n}={\sf f}_n(\bm{Y}_n)$ (i.e., of the {\em real} estimators based on the measured samples) as the sample size $n$ grows with the network size $N$ with some law $n(N)$. The forthcoming section is devoted to establishing which law $n(N)$ is sufficient to achieve consistent learning.
~\hfill$\square$
\end{remark}

\begin{figure*}
\begin{minipage}{.33\linewidth}
{\centering{\bf ~~~~~~Granger}\par\medskip}
\vspace*{5pt}
\centering
{
\includegraphics[scale=.34]{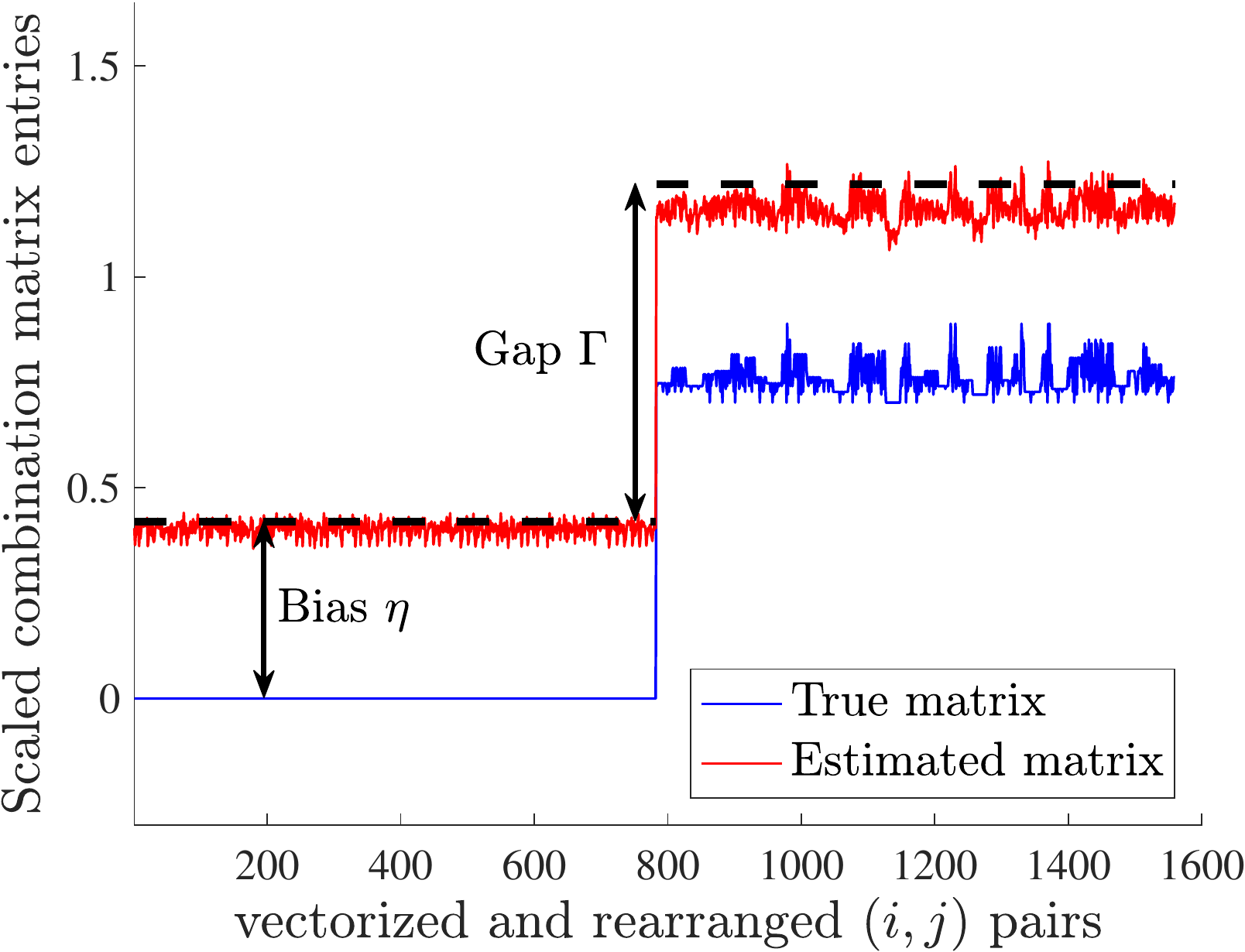}}
\end{minipage}
\begin{minipage}{.33\linewidth}
{\centering{\bf ~~~~~~one-lag}\par\medskip}
\vspace*{5pt}
\centering
{
\includegraphics[scale=.34]{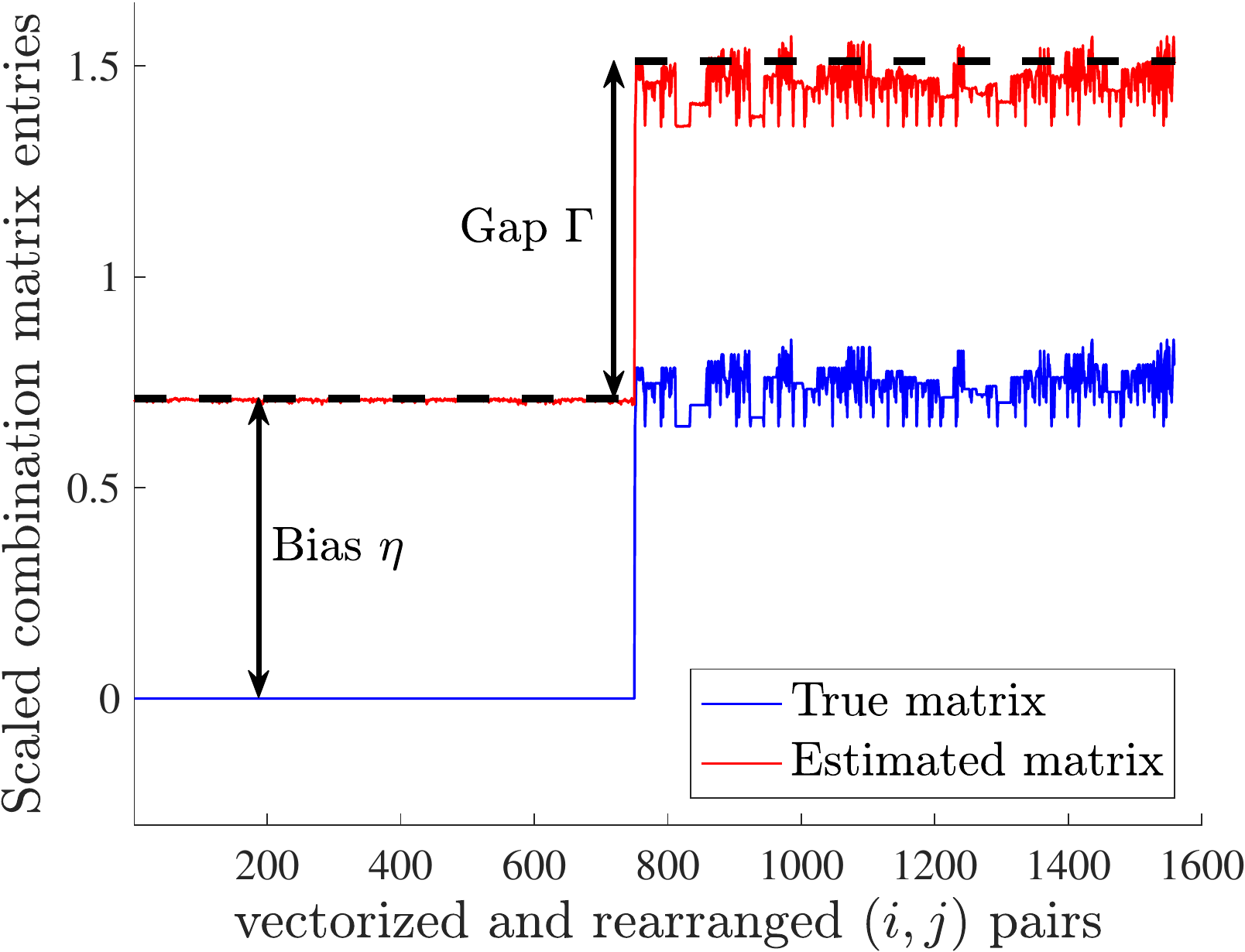}}
\end{minipage}
\begin{minipage}{.33\linewidth}
{\centering{\bf ~~~~~~residual}\par\medskip}
\vspace*{5pt}
\centering
{
\includegraphics[scale=.34]{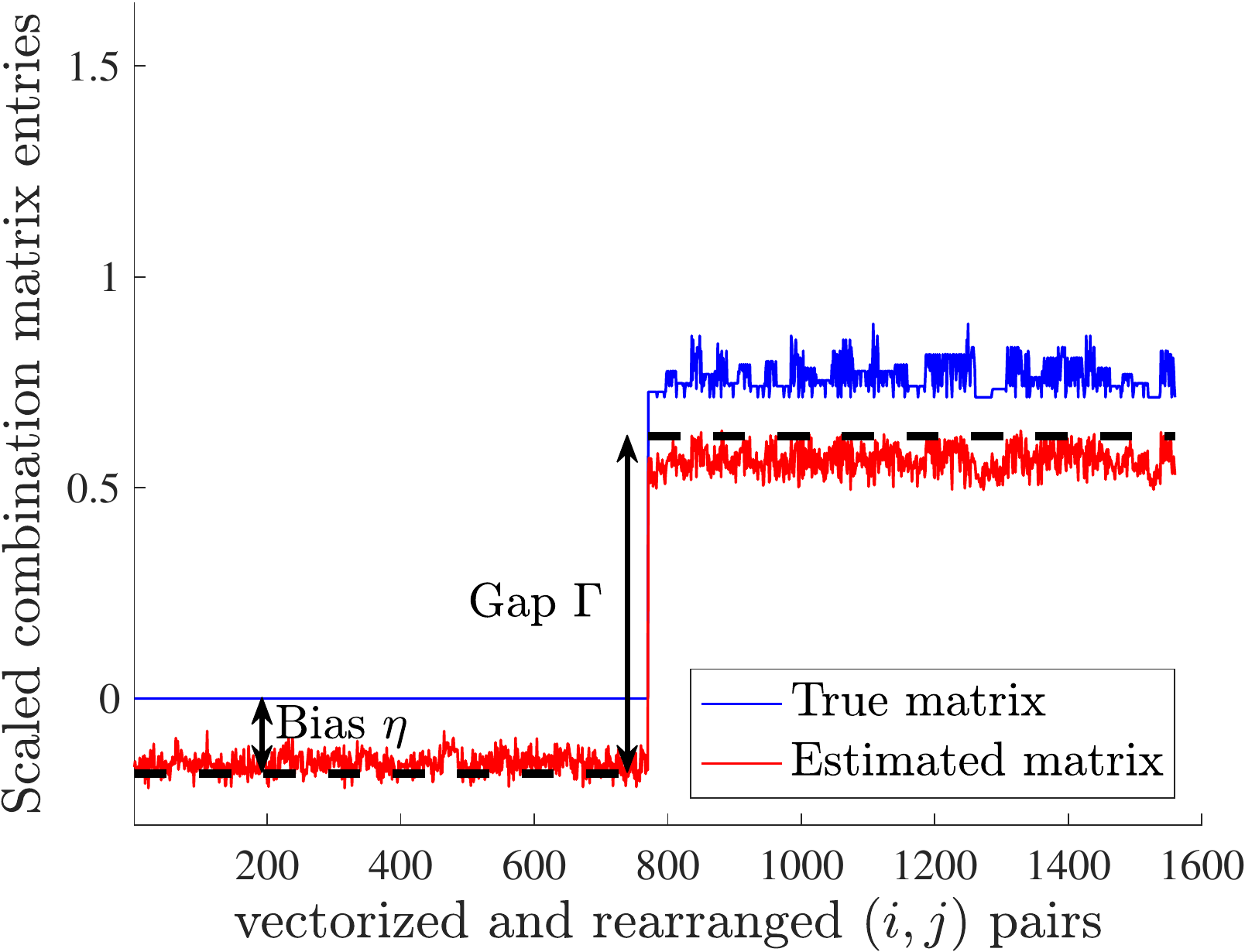}}
\end{minipage}
\caption{
Graphical illustration of Theorem~\ref{theor1}. The three panels refer to the matrix estimators considered in the theorem. 
In each panel, the entries of the {\em true} matrix $\bm{A}_{\mathcal{S}}$ are vectorized following column-major ordering, and the (vectorized) $(i,j)$ pairs are rearranged in such a way that the zero entries appear before the nonzero entries. The same ordering used for the true matrix is applied to the entries of the {\em estimated} matrix, $\widehat{\bm{A}}_{\mathcal{S}}$.
Broken lines display the theoretical values: the lowermost line refers to the error gap, $\eta$, and the uppermost line refers to the quantity $\eta+\Gamma$, where $\Gamma$ is the identifiability gap.}
\label{fig:ExemplifTheor1}
\end{figure*}

The main message conveyed by Theorem~\ref{theor1} is illustrated in Fig.~\ref{fig:ExemplifTheor1}, where we depict: $i)$ the entries of the true combination matrix (disconnected pairs are displayed in blue, connected pairs are displayed in red), vectorized and ordered as shown in the figure, and magnified by $N p_N$; and $ii)$ the entries of the estimated combination matrix, magnified by $N p_N$, vectorized and ordered with the same ordering used for the the true combination matrix. The three panels refer to the three estimators considered in this work, as detailed in the panel titles.
The essential features illustrated in Sec.~\ref{sec:ident} are clearly visible in Fig.~\ref{fig:ExemplifTheor1}. 
First, we can appreciate the emergence of the gap $\Gamma$ and of the bias $\eta$, which, remarkably, match well the theoretical predictions summarized in Table~\ref{tab:Theorem1}, as indicated by the broken lines (theoretical limiting values). 
It is also seen how the bias does not affect separability between the groups of connected and disconnected agents.
Second, we see how the entries are clustered around the values predicted by the theorem. 
Moreover, we remark that the performance of the three learning algorithms cannot be anticipated from the claim of Theorem~\ref{theor1}. Referring to Fig.~\ref{fig:ExemplifTheor1}, the observed values of biases and gaps do not contain information useful to compare the estimators against each other. What plays a role in their performance is the spread of the matrix entries around their limiting values. In the particular example shown in Fig.~\ref{fig:ExemplifTheor1}, the one-lag estimators seems to offer a better concentration (i.e., a lower spread) than the residual estimator, which seems in turn to exhibit a lower spread as compared to the Granger estimator.

Another useful conclusion stemming from Theorem~\ref{theor1} pertains to the dependence of the main quantities in Table~\ref{tab:Theorem1} on the system parameters. We start by examining the identifiability gap. 
First, we have chosen for all three estimator the same scaling sequence, namely, $s_N=N p_N$. With this choice, we see that the gap of the Granger estimator is equal to the bounding constant $\kappa$ that characterizes the regular diffusion matrix~(\ref{eq:fundamentalassum}). 

Let us now focus on the gaps pertaining to the one-lag and to the residual estimators. 
Both of them exhibit two main differences with respect to the Granger estimator. 
First, they depend also on the variance of the input process, $\sigma^2$. This behavior should be expected, since in the Granger estimator the one-lag correlation matrix multiplies the inverse of the correlation matrix, and, hence, the effect of $\sigma^2$ disappears. In contrast, the one-lag and the residual estimators to not cancel this effect. 
Second, when $\kappa\neq\rho$ (a condition that occurs, for instance, for the Laplacian rule), the term $\sigma^2 \kappa$ multiplies a factor that is a function of $\rho-\kappa$. This factor is greater than one for the one-lag estimator, whereas it is smaller than one for the residual estimator. 
However, due to scale-irrelevance, we see that, in the absence of any information about the spread of the matrix entries, the dependence alone upon $\sigma^2$, as well as a magnified/reduced gap do not imply any conclusion about the performance of the pertinent estimators.

Let us switch to the analysis of the bias. Table~\ref{tab:Theorem1} shows how the biases pertaining to the different estimators depend upon the main problem parameters. All estimators depend upon the product $N p_N$ (through the scaling sequence $s_N$), upon the constant $\kappa$, and upon the combination matrix spectral radius $\rho$. 
Notably, only the bias of the Granger estimator depends upon the fraction of monitored nodes $\xi$. This finding makes sense, since the Granger estimator is based upon inversion of a partial matrix (which clearly varies with the number of latent variables), whereas the one-lag and the residual estimator are natively determined by pairwise correlations. In comparison, only the bias of the Granger estimator does not depend upon $\sigma^2$, and this behavior is easily grasped in light of the explanation in the previous paragraph.

\section{Sample Complexity}
Let us consider the following estimators.
\beqa
&&\widehat{\bm{A}}_{\mathcal{S},n}^{\textnormal{(Gra)}}~=[\widehat{\bm{R}}_{1,n}]_{\mathcal{S}}([\widehat{\bm{R}}_{0,n}]_{\mathcal{S}})^{-1},
\label{eq:sampleGra}
\\
&&\widehat{\bm{A}}_{\mathcal{S},n}^{\textnormal{(1-lag)}}=[\widehat{\bm{R}}_{1,n}]_{\mathcal{S}},~~
\label{eq:sampleonelag}
\\
&&\widehat{\bm{A}}_{\mathcal{S},n}^{\textnormal{(res)}}~~=[\widehat{\bm{R}}_{1,n}]_{\mathcal{S}}-[\widehat{\bm{R}}_{0,n}]_{\mathcal{S}}.
\label{eq:sampleres}
\eeqa
Estimators~(\ref{eq:sampleGra}),(\ref{eq:sampleonelag}) and~(\ref{eq:sampleres}) are the {\em sample} versions of~(\ref{eq:Granger0}), (\ref{eq:onelagestdef}) and~(\ref{eq:newResidualSeries0}), respectively. 
These sample estimates are obtained by replacing the true correlation matrices with their empirical counterparts. 
The estimator in~(\ref{eq:sampleGra}) is assumed to be unspecified when $[\widehat{\bm{R}}_{0,n}]_{\mathcal{S}}$ is singular (which happens when $n\leq N$). In order to get consistent estimators, we must ensure that the probability of observing a singular correlation matrix is asymptotically vanishing. 
In order to bypass the instability related to matrix inversion (which could affect also the sample complexity) we now introduce a regularized Granger estimator. 

For $i\in\mathcal{S}$, the $i$-th row of the regularized Granger estimator $\widehat{\bm{A}}^{\textnormal{(reGra)}}_{\mathcal{S},n}$ is a solution to the constrained optimization problem (here $x$ is a row vector):
\beq
\min_{x\in\mathbb{R}^S} 
\left\|
x \, [\widehat{\bm{R}}_{0,n}]_{\mathcal{S}} - [\widehat{\bm{R}}_{1,n}]_i
\right\|_{\infty}~~\textnormal{s.t. }\|x\|_1\leq 1.
\label{eq:regularsamGra}
\eeq
The $\ell_1$-constraint on the admissible row vectors $x$ in~(\ref{eq:regularsamGra}) arises from the fact that the rows of our target estimator, the {\em limiting} Granger estimator, have $\ell_1$-norm bounded by $1$ --- see~(\ref{eq:normboundedby2}) in Appendix~\ref{app:samplecomple}.

It is also useful to observe that, in view of~(\ref{eq:sampleGra}), when the sample correlation matrix is invertible, the non-regularized Granger estimator is the only matrix that yields a zero residual in~(\ref{eq:regularsamGra}). As a result, whenever  $\|\widehat{\bm{A}}^{\textnormal{(Gra)}}_{\mathcal{S},n}\|_{\infty}\leq 1$ the plain and regularized Granger estimators coincide. 

Now, assuming that the empirical correlation matrices converge to the true correlation matrices as $n\rightarrow\infty$, it is easily seen that the four estimators~(\ref{eq:sampleGra})--(\ref{eq:sampleres}) converge to their limiting counterparts and, hence, guarantee condition~(\ref{eq:asystableclass}). 
Since in Theorem~\ref{theor1} we established that the {\em limiting} Granger, one-lag and residual estimators achieve universal local structural consistency, in view of Theorem~\ref{theor:samplecons} this property implies that all the estimators~(\ref{eq:sampleGra})--(\ref{eq:sampleres}) are able to learn well the underlying subgraph of probed nodes, provided that the number of samples $n$ grows as the network size $N$ diverges. The goal of a sample complexity analysis is to establish how this number of samples should scale with $N$ in order to grant correct graph learning. The forthcoming theorem provides an answer for the class of models~(\ref{eq:VARmodel}) with Gaussian input source.  
\begin{theorem}[Sample complexity of the proposed estimators]
\label{theor:samplecomplexity}
Assume that model~(\ref{eq:VARmodel}) holds with Gaussian source data $\{\bm{x}_n\}$, and with initial state $\bm{y}_0$ distributed according to the stationary distribution of the process (so that the overall process $\bm{y}_n$ is stationary).  
Let $\bm{A}$ be a regular diffusion matrix with parameters $\rho$ and $\kappa$, with the network graph drawn according to an Erd\H{o}s-R\'enyi random graph model $\mathscr{G}(N,p_N)$ where the fraction of observable nodes, $S/N$, converges to some nonzero value $\xi$. 
Let $\widehat{\bm{A}}_{\mathcal{S},n}$ be either~(\ref{eq:regularsamGra}), (\ref{eq:sampleonelag}) or~(\ref{eq:sampleres}). 

Then, under the uniform concentration regime where:
\beq
p_N=\omega_N\displaystyle{\frac{\log N}{N}}\stackrel{N\rightarrow\infty}{\longrightarrow} p,~~~~\textnormal{with } \omega_N\rightarrow \infty,
\label{eq:unifconctheosample}
\eeq
we have that:
\beq
\lim_{N\rightarrow\infty} \P[{\sf clu}(\widehat{\bm{A}}_{\mathcal{S},n(N)}) = \bm{G}_{\mathcal{S}}]=1,
\eeq
provided that the number of samples is on the order of:
\beq
\boxed{
n(N)=\Omega\left(\, (N p_N)^2 \log S \, \right)
}
\label{eq:samplecomplaw}
\eeq
Moreover, in the dense regime where $p>0$, the same result holds for the (non-regularized) sample Granger estimator~(\ref{eq:sampleGra}).
\end{theorem}
\begin{IEEEproof}
See Appendix~\ref{app:samplecomple}.
\end{IEEEproof}
We see from~(\ref{eq:samplecomplaw}) that under the dense regime the growth is essentially quadratic in $N$, since $p_N$ converges to some positive constant $p$. 
In order to see what happens under the sparse regime, it is useful to apply~(\ref{eq:unifconctheosample}) and rewrite the sample size in~(\ref{eq:samplecomplaw}) as:
\beq
\omega^2_N (\log N)^3,
\eeq 
revealing that the specific sample complexity under the sparse regime depends on the specific speed of growth of the sequence $\omega_N$, which in fact regulates the sparsity of the problem.

The scaling law found in~(\ref{eq:samplecomplaw}) is significant since it matches well with the scaling laws that have been found in the literature in relation to other graphical models. 
Over graphical models with independent samples, in~\cite{AnandkumarValluvanAOS2013,AnandkumarTanHuangWillskyJMLR2012} it is shown that the sample complexity grows logarithmically with $N$, but proportionally to the inverse square of the minimum absolute entry of the precision matrix. In our setting, the inverse square of the minimum matrix entry scales as $(N p_N)^2$. 
Still with reference to graphical models with independent samples, in~\cite{ChandrasekaranParriloWillskyAOS2012} it is shown that when the degree grows as a polynomial in $\log N$, the sample complexity is on the order of $N$ times a polynomial in $\log N$. 
If we compare this result against our result for the case where $N p_N$ grows poly-logarithmically, an extra-factor $N$ appears. However, this additional factor can be ascribed to the fact that the algorithm in~\cite{ChandrasekaranParriloWillskyAOS2012} manages to estimate consistently both the subgraph of probed nodes and a low-rank matrix that summarizes the influence of latent nodes.
For the case of linear dynamical systems like~(\ref{eq:VARmodel}), an algorithm is proposed in~\cite{JalaliSanghaviICML2012} whose sample complexity scales as $\log N$ when the node degree is kept fixed. 
However, when the degree grows with $N$ (as in our setting), then the sample complexity in~\cite{JalaliSanghaviICML2012} contains an additional $(N p_N)^3$ factor (even if the authors of~\cite{JalaliSanghaviICML2012} indicate that they suspect the exponent could be reduced to $2$). 
In addition, to grant graph recovery a lower bound on the minimum combination-matrix entry is assumed in~\cite{JalaliSanghaviICML2012}, but this bound is not allowed to scale with $N$.

The sample complexity found in Theorem~\ref{theor:samplecomplexity} is primarily determined by the error in estimating the empirical correlations, which is examined in Lemma~\ref{lem:concorr}. 
The limitations of the empirical correlations in high-dimensional settings are well known. 
For example, the correlation matrix is singular when $n\leq N$. Even for larger $n$ the number of parameters to be estimated (i.e., the matrix entries) grows quadratically, thus requiring a high number of samples to get good performance. 
On the other hand, when focusing on the $\|\cdot\|_{\max}$ error norm, it is known that a maximum error up to an $\epsilon$ can be obtained when the number of samples grows {\em logarithmically} with the dimension~\cite{BentoIbrahimiMontanari,HanLuLiuJMLR,LohWainwrightAOS2012,RaoKipnisJavidiEldarGoldsmithCDC2016}. 
This effect is summarized by the $\log S$ factor appearing in~(\ref{eq:samplecomplaw}). 
Unfortunately, this is not the main factor determining the sample complexity.
The main factor is instead the $N p_N$ term. 
This term arises because, in order to guarantee that the matrix entries are well classifiable, we need to guarantee a maximum error up to $\epsilon/(N p_N)$, which is related to the scaling law of the smallest combination-matrix entry. 

Finally, we observe that for the one-lag and the residual estimators the error on empirical correlations automatically relates to the error on the combination matrix. This is not so automatic for the Granger estimator in~(\ref{eq:sampleGra}), where matrix inversion introduces some source of instability that can affect sample complexity. 
This is why in~(\ref{eq:regularsamGra}) we introduced a regularization procedure that helps achieve the sample complexity in~(\ref{eq:samplecomplaw}). 
We notice that in the dense regime, where we would end up in any case with a quadratic growth, the same complexity is achieved by the (non-regularized) Granger estimator.
 
\section{Illustrative Examples}
\label{sec:illustexam}
In this section we report the results of computer simulations that we have conducted to illustrate the theoretical results developed in the previous sections. 
It is useful to highlight the common structure adopted in the forthcoming examples.
\begin{enumerate}
\item
We consider the Granger, the one-lag, and the residual estimator, which are computed both under the assumption of unlimited sample sizes ({\em theoretical} correlation matrices, $R_0$ and $R_1$), as well as under the assumption of finite sample size ({\em empirical} correlation matrices). 
\item
We compare the performance of the estimators in terms of the probability of correct graph recovery, evaluated through Monte Carlo simulations.
\item
The classification/thresholding stage are implemented through the clustering algorithm described in Sec.~\ref{subsec:clualgo}, i.e., we do not rely on any a-priori knowledge of the system parameter to set the classification threshold. 
\item
Finally, the input source samples, $\bm{x}_i(n)$, are i.i.d. samples from a standard Gaussian distribution, and we set $\sigma=1$.
\end{enumerate}

\begin{figure} [t]
\begin{center}
\includegraphics[scale= 0.45]{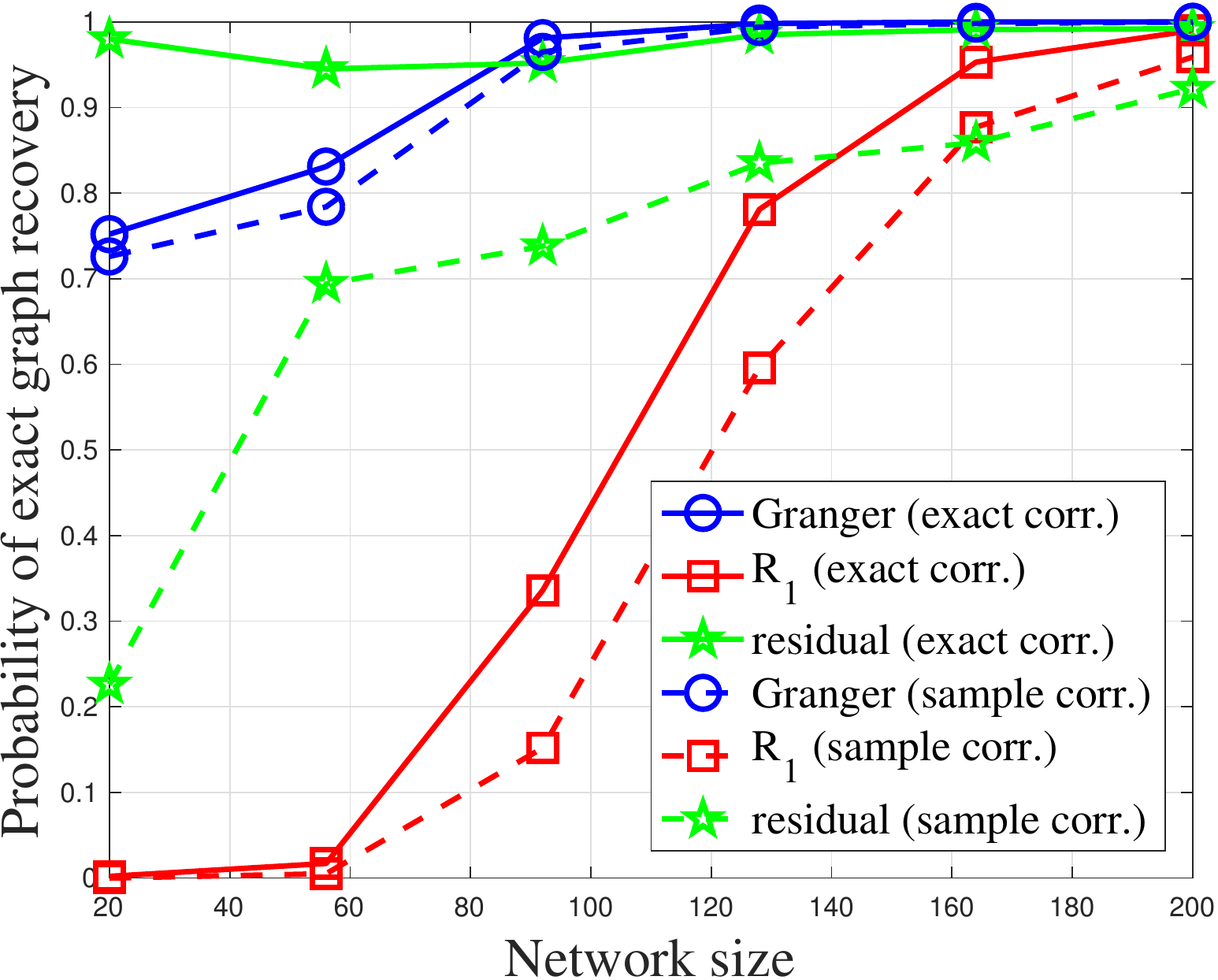}
\caption{
{\em Uniform-and-dense} regime with fraction of monitored nodes $\xi=0.6$ and a Metropolis combination rule~(\ref{eq:MetroMat}).
The figure displays the probability of correct graph recovery (estimated over $10^3$ Monte Carlo runs) for the three estimators detailed in the legend, for both the cases where the correlation matrices are perfectly known, or where they are estimated from time samples. In particular, as regards the Granger estimator based on sample correlations, we used~(\ref{eq:sampleGra}).
The constant connection probability is $p_N=p=0.1$, and the parameters of the Metropolis combination matrix are $\rho=\kappa=0.99$. The number of time samples used for the empirical correlation estimators scales as $(N p_N)^2 \log S$, with the last point being equal to $5\times 10^5$. 
}
\label{fig:strongdense1}
\end{center}
\end{figure}

\begin{figure} [t]
\begin{center}
\includegraphics[scale= 0.45]{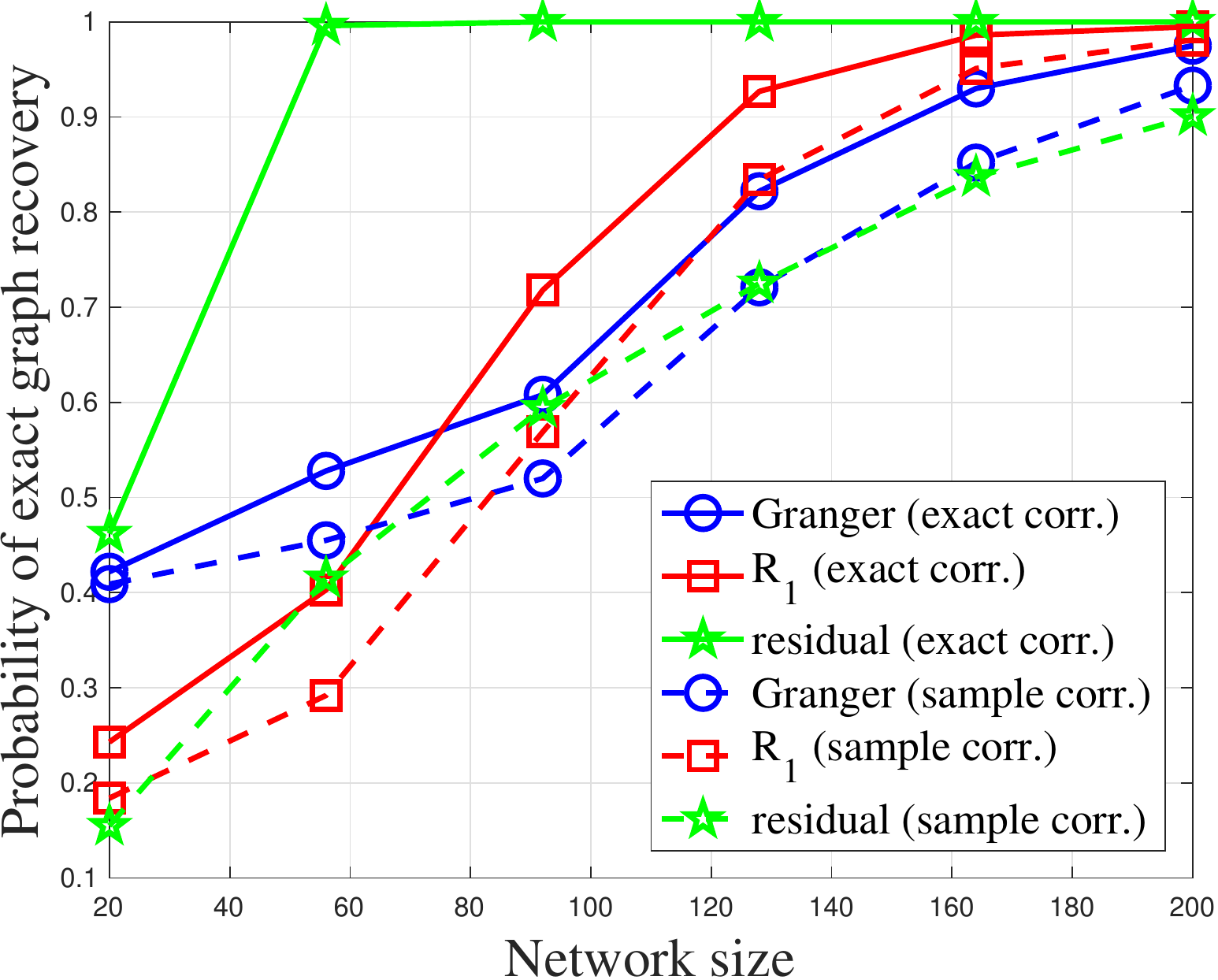}
\caption{
{\em Uniform-and-dense} regime with fraction of monitored nodes $\xi=0.2$ and a Laplacian combination rule~(\ref{eq:LapMat}).
The figure displays the probability of correct graph recovery (estimated over $10^3$ Monte Carlo runs) for the three estimators detailed in the legend, for both the cases where the correlation matrices are perfectly known, or where they are estimated from time samples. In particular, as regards the Granger estimator based on sample correlations, we used~(\ref{eq:sampleGra}).
The constant connection probability is $p_N=p=0.1$, and the parameters of the Laplacian combination matrix are $\rho=0.99$ and $\lambda=0.9$ (yielding $\kappa=\rho\lambda=0.891$). The number of time samples used for the empirical correlation estimators scales as $(N p_N)^2 \log S$, with the last point being equal to $5\times 10^5$.
}
\label{fig:strongdense2}
\end{center}
\end{figure}

In Fig.~\ref{fig:strongdense1} we offer a sample of the behavior observed in the uniform-and-dense regime (constant connection probability). As a general comment, we see that all the estimators, both with unlimited or limited sample sizes, match well the theoretical predictions. 
Let us now start by examining the unlimited-samples estimators.
In this particular example, we observe that, after an initial transient with small network sizes, the Granger estimator outperforms the other two estimators. A legitimate question at this stage is whether or not this is a general behavior.

\begin{figure} [t]
\begin{center}
\includegraphics[scale= 0.45]{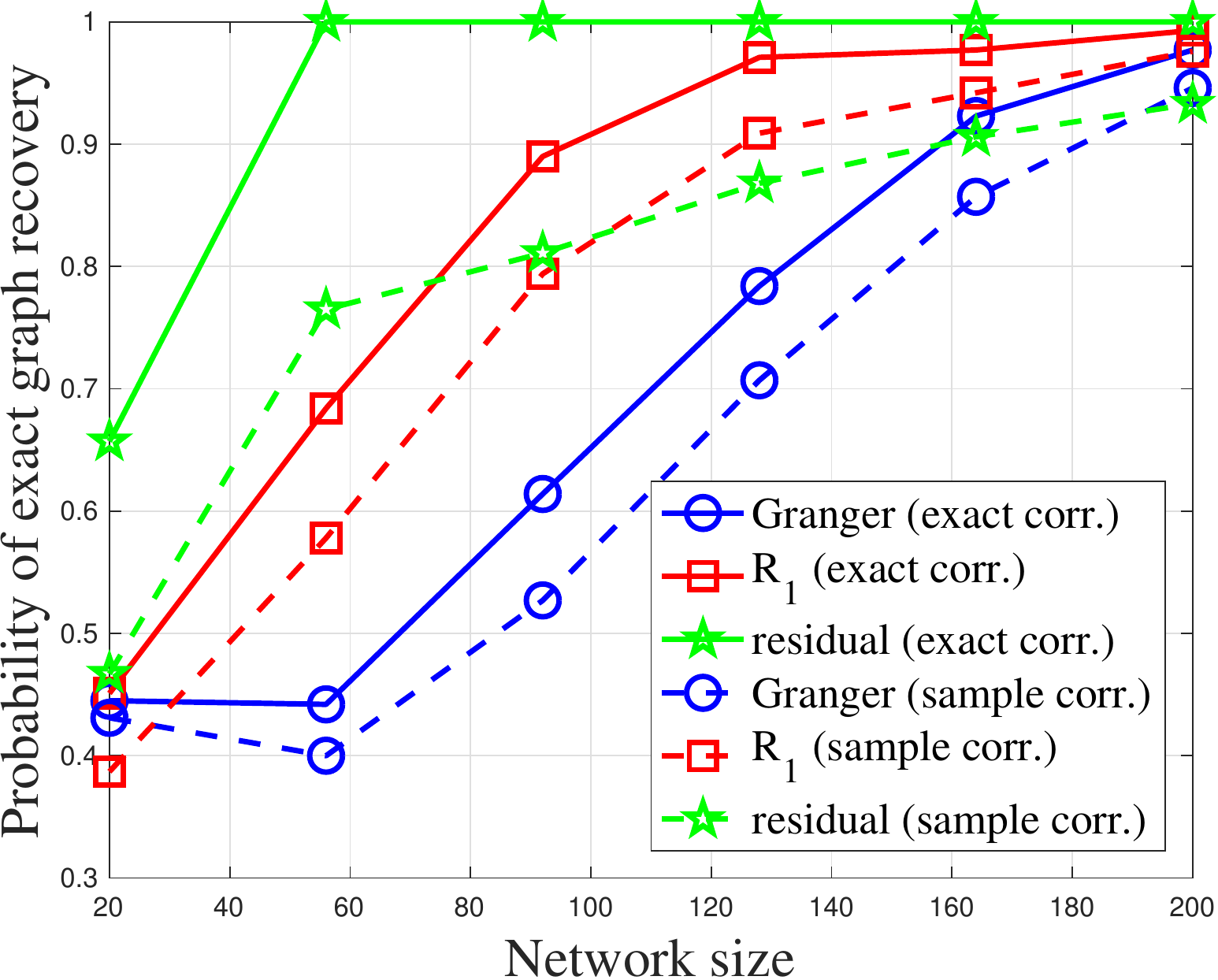}
\caption{
{\em Uniform-and-sparse} regime with fraction of monitored nodes $\xi=0.2$ and a Laplacian combination rule~(\ref{eq:LapMat}).
The figure displays the probability of correct graph recovery (estimated over $10^3$ Monte Carlo runs) for the three estimators detailed in the legend, for both the cases where the correlation matrices are perfectly known, or where they are estimated from time samples. In particular, as regards the Granger estimator based on sample correlations, we used~(\ref{eq:regularsamGra}).
The connection probability is $p_N=0.25 \,(\log N)/\sqrt{N}$, and the parameters of the Laplacian combination matrix are $\rho=0.99$ and $\lambda=0.9$ (yielding $\kappa=\rho\lambda=0.891$). The number of time samples used for the empirical correlation estimators scales as $(N p_N)^2 \log S$, with the last point being equal to $5\times 10^5$.
}
\label{fig:strongsparse}
\end{center}
\end{figure}

In order to answer, in Fig.~\ref{fig:strongdense2} we consider a different setting, and we see that the residual estimator outperforms the one-lag estimator, which in turn outperforms the Granger estimator. 
One reason behind the discrepancies between Figs.~\ref{fig:strongdense1} and~\ref{fig:strongdense2} has to be ascribed to the fact that, in the latter case, the fraction of observed nodes is reduced substantially. 
Accordingly, the Granger estimator seems more sensitive to the level of observability, whereas it tends to be more performing when one approaches the regime of full observability. 
This sounds reasonable since under full observability the Granger estimator coincides exactly with the {\em true} combination matrix while the one-lag and the residual estimators do not. 
In summary, one notable conclusion is that estimators that perform very well under full observability need not be as good under partial observability. 
In order to complete the analysis, in Fig.~\ref{fig:strongsparse} we show one example pertaining to the uniform-and-sparse regime, where we see that conclusions similar to those drawn for the uniform-and-dense regime apply.

We continue by examining the performance of the {\em actual} estimators, namely, the estimators based on sample correlation matrices. 
All these estimators have been implemented with a number of samples growing with the scaling law predicted by Theorem~\ref{theor:samplecomplexity}. 
Two main trends emerge. 
First, with the predicted scaling law, all the estimators increase their performance as the network size grows. 
Second, a difference emerges between the residual estimator and the other two estimators. 
In particular, for the same number of samples, the residual estimator (dashed green line) is more distant from its unlimited-samples performance (solid green line) than what is observed for the Granger and the one-lag estimators. 
As a consequence, even when the limiting residual estimator is the best one, the {\em sample} residual estimator (for the number of samples considered in the simulations) is outperformed by the other two estimators, implying that the residual estimator is more sensitive to the error associated to the convergence of the sample correlation matrices. On the other hand, increasing the number of samples one can in principle approximate with arbitrary fidelity the limiting performance (solid green line) of the residual estimator.

\begin{table*}[t]
  \begin{center}
    \begin{tabular}{l|c|c} 
      ~~~~~~~~~~Estimator & Uniform concentration regime & Very sparse regime\\
      \hline
      \hline
      &&
      \\
      Granger: $[\bm{R}_1]_{\mathcal{S}} ([\bm{R}_0]_{\mathcal{S}})^{-1}$ & 
      Universal local structural consistency &
      Mild consistency\footnotemark~\cite{tomo}\\
      &&
      \\
      \hline
      &&
      \\
      one-lag: $[\bm{R}_1]_{\mathcal{S}}$ & 
      Universal local structural consistency &
      ?\\
      &&
      \\
      \hline
      &&
      \\
      residual: $[\bm{R}_1]_{\mathcal{S}} - [\bm{R}_0]_{\mathcal{S}}$ & 
      Universal local structural consistency &
      ?\\
      &&
      \\
      \hline
    \end{tabular}
    \vspace*{5pt}
    \caption{Summary of the known and yet unknown consistency properties of the estimators considered in this work.
    }
    \label{tab:summary}
  \end{center}
\end{table*}
\footnotetext{Actually, for the Granger estimator under the very sparse regime, the available consistency result is a {\em mild consistency} result provided in terms of average fraction of errors.}

\section{Concluding Remarks and Open Issues}
This work examined the problem of graph learning when data can be collected from a limited subset of nodes. 
The goal is to learn the topology of the subgraph of monitored nodes. 
We considered one useful class of generative graph models, namely, the random graph Erd\H{o}s-R\'enyi model. 
Three learning algorithms were developed, which were proved to be structurally consistent in the thermodynamic limit where the overall network size grows without bound. 
We explored various regimes of connectivity, including the often overlooked regime of {\em dense} connectivity. 
One revealing conclusion stemming from our analysis is that the {\em statistical concentration} of node degrees plays a central role for consistent graph learning, sometimes even more important than the {\em sparsity} of connections.

Several works in the literature of graph learning assume sparsity in the graph of connections. The role of sparsity in graph learning can be (at least) twofold. On one hand, sparsity can be leveraged to reduce the complexity associated to the estimators (we have seen this effect also in our analysis of the sample complexity). On the other hand, in the presence of latent variables, sparsity can be leveraged since it reduces the (unseen) effect of the latent unobserved nodes on the probed nodes~\cite{tomo_journal, ChandrasekaranParriloWillskyAOS2012}. One interesting conclusion stemming from our analysis is that, under the setting considered in this work (Erd\H{o}s-R\'enyi graphs and regular diffusion matrices) an important role is played by another structural property of the graph, namely, the statistical concentration of node degrees.

A succinct summary of the results is given in Table~\ref{tab:summary}. 
This work focused on the uniform concentration regimes, whereas some results for the very sparse regime are already available from previous works. 
Moreover, we used question marks to highlight some open issues. 
We see that results about consistency under the very sparse regime (which was not dealt with here) are not available for the one-lag and for the residual estimators. 
It is also interesting to note that the results for the very sparse regime obtained, e.g., in~\cite{tomo,tomo_journal} ground on and require tools different from those adopted here. Moreover, we observe that we are not aware of results available for the disconnected regime.

There are further open issues that may deserve attention. 
One issue concerns {\em directed} graphs, which are relevant, e.g., in the context of causation. 
We believe that the graph-edge-domain approach developed in this work (more specifically, the recursive inequalities on the matrix-power entries obtained in Appendix~\ref{app:Matrecursions} can be generalized to get insights about the directed graph setting. 

This work focused on the Erd\H{o}s-R\'enyi model, and used certain regularity assumptions on the diffusion matrix. 
A useful extension would be to examine structural consistency for other graph models, and/or under different regularity  assumptions. For example, one might have a certain heterogeneity in the network (e.g., different connectivity across nodes), so as to observe multiple gaps in the estimated matrices, and an interesting question would be whether or not consistency can be achieved under these conditions.

Another open problem regards the search for {\em optimal} learning algorithms, where optimality can be formulated by taking into account the learning performance as well as the different sources of complexity.
In this connection, establishing the fundamental limits of the system in terms of computational and sample complexity would constitute a significant advance.

\begin{appendices}

\section{Useful Properties of Maximal and Minimal Degrees}

We denote by $\bm{\mathcal{B}}_i(N,q)$, with $i=1,2,\ldots,K$, a sequence of $K$ binomial random variables (not necessarily independent) with success probability $q$ over $N$ independent Bernoulli trials. Moreover, we denote by $\bm{\mathcal{B}}_{\max}(N,q,K)$ and $\bm{\mathcal{B}}_{\min}(N,q,K)$ the maximum and the minimum over this sequence, respectively.
The following two relationships are standard inequalities arising from the application of the Chernoff bounding technique, and will be the fundamental building blocks to characterize the asymptotic behavior of several random quantities arising in our problem.
The inequalities are as follows.\footnote{For any $t>0$, we can write: 
\beq
\P[\bm{\mathcal{B}}_i(N,q)\geq x]=\P[e^{\bm{\mathcal{B}}_i(N,q)\,t}\geq e^{x t}]\leq
e^{- x t}\,\E[e^{\bm{\mathcal{B}}_i(N,q)\,t}],
\label{eq:MarkovChernoff}
\eeq
where the latter inequality is an application of Markov's inequality. Since a binomial variable of parameters $N$ and $q$ is the sum of $N$ independent Bernoulli variables with success probability equal to $q$, we can further write:
\beq
\E[e^{\bm{\mathcal{B}}_i(N,q)\,t}]=(q e^t + 1 - q)^{N}=(1 + q (e^t - 1)^{N}\leq e^{N q (e^t - 1)},
\label{eq:MGFbinom}
\eeq
where the latter inequality follows by observing that, for $z>0$, one has $(1+z)^N\leq e^{N z}$. Combining~(\ref{eq:MGFbinom}) with~(\ref{eq:MarkovChernoff}) yields~(\ref{eq:basicmax}). Equation~(\ref{eq:basicmin}) is worked out with a similar technique.
} 
For any $t>0$:
\beqa
\P[\bm{\mathcal{B}}_{\max}(N,q,K)\geq x]
&\leq&
K 
e^{- x t + N q(e^t - 1)},
\label{eq:basicmax}
\\
\P[\bm{\mathcal{B}}_{\min}(N,q,K)\leq x]
&\leq&
K 
e^{x t - N q(1 - e^{-t})}.
\label{eq:basicmin}
\eeqa
We now apply these fundamental bounds to some specific random variables that are of interest in our setting. 

We start by characterizing the behavior of the variables:
\beq
\bm{\mathcal{B}}_{\max}(N,p_N,N), \qquad \bm{\mathcal{B}}_{\min}(N,p_N,N),
\eeq
which, as we will see, are useful to characterize the behavior of the maximal and minimal degree of the graphs that we use in this work. The forthcoming lemma contains fundamental (classic) results about the asymptotic behavior of $\bm{\mathcal{B}}_{\max}(N,p_N,N)$ and $\bm{\mathcal{B}}_{\min}(N,p_N,N)$ under the different regimes for the probability $p_N$.
\begin{lemma}[Asymptotic scaling of $\bm{\mathcal{B}}_{\max}(N,p_N,N)$ and $\bm{\mathcal{B}}_{\min}(N,p_N,N)$] 
\label{lem:Bmaxmin}
Let the probability $p_N$ scale with $N$ according to~(\ref{eq:strongconc}). 
Then:
\beq
\boxed{
\frac{\bm{\mathcal{B}}_{\max}(N,p_N,N)}{N p_N}\stackrel{\textnormal{p}}{\longrightarrow} 1,~~
\frac{\bm{\mathcal{B}}_{\min}(N,p_N,N)}{N p_N}\stackrel{\textnormal{p}}{\longrightarrow} 1
}
\label{eq:Bmaxminstrong}
\eeq
\end{lemma}
\begin{IEEEproof} 
The following inequality, holding for all $\epsilon>0$, is easily obtained from~(\ref{eq:basicmax}) by setting $K=N$, $q=p_N$, $x=(1+\epsilon)Np_N$, $t=\log(1+\epsilon)$, and $g_{\epsilon}\dfz 1+[1+\epsilon][\log(1+\epsilon) - 1]$:
\beq
\P[\bm{\mathcal{B}}_{\max}(N,p_N,N)\geq (1+\epsilon) N p_N]
\leq
N e^{- Np_N g_{\epsilon}}.
\label{eq:dmaxsc}
\eeq 
Using now~(\ref{eq:strongconc}) in~(\ref{eq:dmaxsc}) we get: 
\beq
\P[\bm{\mathcal{B}}_{\max}(N,p_N,N)\geq (1+\epsilon) N p_N]
\leq
N^{1 - \omega_N g_{\epsilon}}\stackrel{N\rightarrow\infty}\longrightarrow 0,
\label{eq:dmaxsc}
\eeq 
which follows because $g_{\epsilon}>0$ for all $\epsilon>0$, as $g_0=0$ and $dg_{\epsilon}/d\epsilon >0$ for all $\epsilon>0$.

Likewise, the following inequality, holding for all $0<\epsilon<1$, is easily obtained from~(\ref{eq:basicmin}) by setting $K=N$, $q=p_N$, $x=(1-\epsilon)N p_N$, $t=-\log(1-\epsilon)$, and $h_{\epsilon}\dfz 1-[1-\epsilon][1 - \log(1-\epsilon)]$: 
\beq
\P[\bm{\mathcal{B}}_{\min}(N,p_N,N)\leq (1-\epsilon) N p_N]
\leq
N^{1 - \omega_N h_{\epsilon}}\stackrel{N\rightarrow\infty}\longrightarrow 0,
\label{eq:dmindense}
\eeq 
which follows because $h_{\epsilon}>0$ for all $0<\epsilon<1$, as $h_0=0$ and $dh_{\epsilon}/d\epsilon >0$ for all $0<\epsilon<1$. 
By joining~(\ref{eq:dmaxsc}) with~(\ref{eq:dmindense}), and observing that $\bm{\mathcal{B}}_{\max}(N,p_N,N)\geq\bm{\mathcal{B}}_{\min}(N,p_N,N)$, we conclude that~(\ref{eq:Bmaxminstrong}) holds true.
\end{IEEEproof}

As simple corollaries to Lemma~\ref{lem:Bmaxmin}, we can now obtain the characterization of the maximal and minimal degrees.

\begin{corollary}[Behavior of $\bm{d}_{\max}$ and $\bm{d}_{\min}$]
\label{cor:dmaxmin}
If the connection probability of the Erd\H{o}s-R\'enyi model obeys~(\ref{eq:strongconc}), then we have:
\beq
\boxed{
\frac{\bm{d}_{\max}}{N p_N}\stackrel{\textnormal{p}}{\longrightarrow} 1,\qquad
\frac{\bm{d}_{\min}}{N p_N}\stackrel{\textnormal{p}}{\longrightarrow} 1,\qquad \textnormal{[Uniform concentration]}
}
\label{eq:dmaxmin}
\eeq
\end{corollary}
\begin{IEEEproof}
The degree of a single node is equal to $1$ plus (because in our setting the degree counts also the node itself) a binomial random variable with parameters $N-1$ and $p_N$. Therefore, we have the following representation:
\beqa
\bm{d}_{\max}&=&1 + \bm{\mathcal{B}}_{\max}(N-1,p_N,N),\label{eq:dmaxappdef}\\
\bm{d}_{\min}&=&1 + \bm{\mathcal{B}}_{\min}(N-1,p_N,N).\label{eq:dminappdef}
\eeqa
In order to obtain useful bounds involving $\bm{d}_{\max}$ and $\bm{d}_{\min}$, let us introduce a modified sequence of binomial variables, obtained by adding one more Bernoulli trial to each binomial variable $\bm{\mathcal{B}}_i(N,p_N)$, with $i=1,2,\ldots,N$. The corresponding maximum and minimum taken over the modified sequence will be denoted by $\widetilde{\bm{\mathcal{B}}}_{\max}(N,p_N,N)$ and $\widetilde{\bm{\mathcal{B}}}_{\min}(N,p_N,N)$, respectively. 
Since a Bernoulli variable can be either zero or one, from~(\ref{eq:dmaxappdef}) and~(\ref{eq:dminappdef}) we get readily the following bounds:
\beqa
\bm{d}_{\max}&\leq& 1 + \widetilde{\bm{\mathcal{B}}}_{\max}(N,p_N,N),\\
\bm{d}_{\min}&\geq& \widetilde{\bm{\mathcal{B}}}_{\min}(N,p_N,N),
\eeqa
and, hence, the claims of the corollary follow readily from Lemma~\ref{lem:Bmaxmin}, with the factor $1$ playing no role as $N\rightarrow\infty$. 
\end{IEEEproof}

\subsection{Another Useful Concentration Result} 
\begin{lemma}[Maximum and minimum of $N^2$ binomial variables with success probability $p_N^2$]
\label{lem:fundmaxN2lemma}
Assume that the success probability obeys~(\ref{eq:strongconc}). Then we have that:
\beq
\boxed{
\displaystyle{\frac{\bm{\mathcal{B}}_{\max}(N,p_N^2,N^2)}{N p_N}\stackrel{\textnormal{p}}{\longrightarrow} p},
\qquad
\displaystyle{\frac{\bm{\mathcal{B}}_{\min}(N,p_N^2,N^2)}{N p_N}\stackrel{\textnormal{p}}{\longrightarrow} p}
}
\label{eq:fundlemmaNsquare}
\eeq
\end{lemma}
\begin{IEEEproof}
If $p_N\rightarrow p>0$, we can set $p'_N=p_N^2$, and obviously $p'_N$ converges to $p^2>0$, implying that the binomial variables of parameters $N$ and $p'_N$ are generated under the uniform concentration regime. The result in~(\ref{eq:fundlemmaNsquare}) then readily follows from~(\ref{eq:basicmax}) and~(\ref{eq:basicmin}). 

If $p_N\rightarrow p=0$, it suffices to prove the claim for the maximum. 
Applying~(\ref{eq:basicmax}) we can write:
\beqa
\P[\bm{\mathcal{B}}_{\max}(N,p_N^2,N^2)\geq \epsilon N p_N]
&\leq&
N^2
e^{- N p_N [\epsilon t - p_N(e^t - 1)]}\nonumber\\
&=&N^2
N^{- \omega_N [\epsilon t - p_N(e^t - 1)]}.\nonumber\\
\label{eq:maxN2pN2}
\eeqa
where, in the last step, we used the equality $N p_N=\omega_N\log N$ that follows from~(\ref{eq:strongconc}). 
Moreover, since we are considering the case where $p_N\rightarrow 0$ as $N\rightarrow\infty$, for any $\epsilon'>0$ and for sufficiently large $N$ we will have $p_N<\epsilon'$, so that, asymptotically, it is legitimate to replace~(\ref{eq:maxN2pN2}) with:
\beq
\P[\bm{\mathcal{B}}_{\max}(N,p_N^2,N^2)\geq \epsilon N p_N]
\leq
N^2
N^{- \omega_N [\epsilon t - \epsilon'(e^t - 1)]}.
\label{eq:maxN2pN2tris}
\eeq
Now, choosing $\epsilon'$ small enough so that $\epsilon'(e^t - 1)<\epsilon t$,
we finally get:
\beq
\P[\bm{\mathcal{B}}_{\max}(N,p_N^2,N^2)\geq \epsilon N p_N]
\stackrel{N\rightarrow\infty}{\longrightarrow} 0,
\label{eq:maxN2pN2fin}
\eeq
which completes the proof of the lemma.
\end{IEEEproof}

\begin{table*}[t]
  \begin{center}
    \begin{tabular}{l|l|l|c} 
&~~~~~~~~~~~~{\bf Random variable} & ~~~~~~~~~~~~{\bf Random variable} & {\bf Limit (in probability)}\\
      \hline
      \hline
1)&
$\mathfrak{M}_a\dfz\displaystyle{\max_{\substack{i,j\in\mathcal{N} \\ i\neq j}}\bm{a}_{i j}}$  &
$\mathfrak{m}_a\dfz\displaystyle{\min_{\substack{i,j\in\mathcal{N} \\ i\neq j}}\bm{a}_{i j}}$  &
 $0$\\	\hline
2)&
$\mathfrak{M}_{a,\textnormal{self}}\dfz\displaystyle{\max_{i\in\mathcal{N}} \bm{a}_{i i}}$	&
$\mathfrak{m}_{a,\textnormal{self}}\dfz\displaystyle{\min_{i\in\mathcal{N}} \bm{a}_{i i}}$	& 
$\rho-\kappa$\\ \hline
3)&
$\mathfrak{M}_{a_2,\textnormal{self}}\dfz\displaystyle{\max_{i\in\mathcal{N}} \bm{a}^{(2)}_{i i}}$ &
$\mathfrak{m}_{a_2,\textnormal{self}}\dfz\displaystyle{\min_{i\in\mathcal{N}} \bm{a}^{(2)}_{i i}}$ & 
$(\rho-\kappa)^2$\\ \hline
4)&
$\mathfrak{M}_{a,\textnormal{sum}}\dfz\displaystyle{\max_{\substack{i,j\in\mathcal{N} \\ i\neq j}}\sum_{\substack{\ell\in\mathcal{N} \\ \ell \neq j}} \bm{a}_{i \ell}}$	&
$\mathfrak{m}_{a,\textnormal{sum}}\dfz\displaystyle{\min_{\substack{i,j\in\mathcal{N} \\ i\neq j}}\sum_{\substack{\ell\in\mathcal{N} \\ \ell \neq j}} \bm{a}_{i \ell}}$	& 
$\rho$\\ \hline
5)&
$\mathfrak{M}_{c,\textnormal{self}} \dfz\displaystyle{\max_{\ell \in \mathcal{S}'}  \bm{c}_{\ell \ell}} $ & 
$\mathfrak{m}_{c,\textnormal{self}} \dfz\displaystyle{\min_{\ell \in \mathcal{S}'}  \bm{c}_{\ell \ell}} $ & 
$(\rho-\kappa)^2$\\  \hline
6)&
$\mathfrak{M}^{(\mathcal{S}')}_{a,\textnormal{sum}}\dfz\displaystyle{\max_{\substack{\ell,m \in \mathcal{S}' \\ \ell\neq m}} 
\sum_{\substack{h\in\mathcal{S}' \\ h \neq \ell,m}} \bm{a}_{h m}}$ &
$\mathfrak{m}^{(\mathcal{S}')}_{a,\textnormal{sum}}\dfz\displaystyle{\min_{\substack{\ell,m \in \mathcal{S}' \\ \ell\neq m}} 
\sum_{\substack{h\in\mathcal{S}' \\ h \neq \ell,m}} \bm{a}_{h m}}$ & 
$\kappa (1-\xi)$\\  \hline
7)&
$\widetilde{\mathfrak{M}}^{(\mathcal{S}')}_{a,\textnormal{sum}}\dfz
\displaystyle{\max_{i\in\mathcal{S}}\sum_{\substack{\ell\in\mathcal{S}'}}\bm{a}_{i\ell}}$ &
$\widetilde{\mathfrak{m}}^{(\mathcal{S}')}_{a,\textnormal{sum}}\dfz
\displaystyle{\min_{i\in\mathcal{S}}\sum_{\substack{\ell\in\mathcal{S}'}}\bm{a}_{i\ell}}$ & 
$\kappa (1-\xi)$\\  \hline
8)&
$\widetilde{\widetilde{\mathfrak{M}}}^{(\mathcal{S}')}_{a,\textnormal{sum}}\dfz
\displaystyle{\max_{i,j\in\mathcal{S}}
\sum_{\ell\in\mathcal{S}'}\bm{a}_{i\ell}
\sum_{\substack{h\in\mathcal{N} \\ h\neq j}}\bm{a}_{h j}
\sum_{\substack{m\in\mathcal{S}' \\ m\neq \ell,h}}
\bm{a}_{m h}
}$
& 
$\widetilde{\widetilde{\mathfrak{m}}}^{(\mathcal{S}')}_{a,\textnormal{sum}}\dfz
\displaystyle{\min_{i,j\in\mathcal{S}}
\sum_{\ell\in\mathcal{S}'}\bm{a}_{i\ell}
\sum_{\substack{h\in\mathcal{N} \\ h\neq j}}\bm{a}_{h j}
\sum_{\substack{m\in\mathcal{S}' \\ m\neq \ell,h}}
\bm{a}_{m h}
}$&
$\kappa^3 (1-\xi)^2$\\  \hline
9)&
$\widetilde{\mathfrak{M}}^{(\mathcal{S}')}\dfz \displaystyle{\max_{i\in\mathcal{S}}\sum_{\substack{\ell,m\in\mathcal{S}' \\ \ell\neq m}}\bm{a}_{i\ell}\bm{a}_{\ell m}}$ &
$\widetilde{\mathfrak{m}}^{(\mathcal{S}')}\dfz \displaystyle{\min_{i\in\mathcal{S}}\sum_{\substack{\ell,m\in\mathcal{S}' \\ \ell\neq m}}\bm{a}_{i\ell}\bm{a}_{\ell m}}$ &
$\kappa^2 (1-\xi)^2$\\  \hline
10)&
$\widetilde{\widetilde{\mathfrak{M}}}^{(\mathcal{S}')}\dfz\displaystyle{\max_{\substack{i,j\in\mathcal{S} \\ i\neq j}}\sum_{\substack{\ell, m\in\mathcal{S}' \\ \ell\neq m}}\bm{a}_{i\ell}\bm{a}_{m j}}$ &
$\widetilde{\widetilde{\mathfrak{m}}}^{(\mathcal{S}')}\dfz\displaystyle{\min_{\substack{i,j\in\mathcal{S} \\ i\neq j}}\sum_{\substack{\ell, m\in\mathcal{S}' \\ \ell\neq m}}\bm{a}_{i\ell}\bm{a}_{m j}}$ &
$\kappa^2 (1-\xi)^2$\\  \hline
11)&
$\mathfrak{M}_{c,\textnormal{sum}} \dfz\displaystyle{\max_{\substack{\ell,m \in \mathcal{S}' \\ \ell\neq m}}
\sum_{\substack{h\in\mathcal{S}'\\h\neq m}} \bm{c}_{\ell h}}$ &
$\mathfrak{m}_{c,\textnormal{sum}} \dfz\displaystyle{\min_{\substack{\ell,m \in \mathcal{S}' \\ \ell\neq m}}
\sum_{\substack{h\in\mathcal{S}'\\h\neq m}} \bm{c}_{\ell h}}$ &
$\rho^2 - 2\rho\kappa\xi+\kappa^2 \xi$\\  \hline
12)&
$\mathfrak{M}\dfz\displaystyle{\max_{\substack{i,j\in\mathcal{N} \\ i \neq j}}  \sum_{\substack{\ell\in\mathcal{N} \\ \ell \neq i,j}}
\bm{a}_{i \ell}\bm{a}_{\ell j}}$ &
$\mathfrak{m}\dfz\displaystyle{\min_{\substack{i,j\in\mathcal{N} \\ i \neq j}}  \sum_{\substack{\ell\in\mathcal{N} \\ \ell \neq i,j}}
\bm{a}_{i \ell}\bm{a}_{\ell j}}$ &
$\left.
\begin{array}{ll}
N p_N \,\mathfrak{M} \\
\\
N p_N \,\mathfrak{m}
\end{array}
\right\}
\stackrel{\textnormal p}{\longrightarrow} \kappa^2 p$\\ \hline
13)&
$\mathfrak{M}^{(\mathcal{S}')}\dfz\displaystyle{\max_{\substack{i,j\in\mathcal{N} \\ i \neq j}}  \sum_{\substack{\ell\in\mathcal{S}' \\ \ell \neq i,j}}\bm{a}_{i \ell}\bm{a}_{\ell j}}$ &
$\mathfrak{m}^{(\mathcal{S}')}\dfz\displaystyle{\min_{\substack{i,j\in\mathcal{N} \\ i \neq j}}  \sum_{\substack{\ell\in\mathcal{S}' \\ \ell \neq i,j}}\bm{a}_{i \ell}\bm{a}_{\ell j}}$ &
$\left.
\begin{array}{ll}
N p_N\, \mathfrak{M}^{(\mathcal{S}')} \\
\\
N p_N\, \mathfrak{m}^{(\mathcal{S}')}
\end{array}
\right\}
\stackrel{\textnormal p}{\longrightarrow} \kappa^2 p (1-\xi)$\\
\hline
14)&
$\mathfrak{M}^{(\mathcal{S}')}_{a_3,\textnormal{sum}}\dfz\displaystyle{\max_{\substack{i,j\in\mathcal{S} \\ i\neq j}}\sum_{\substack{\ell, m\in\mathcal{S}' \\ \ell\neq m}}\bm{a}_{i\ell}\bm{a}_{\ell m}\bm{a}_{m j}}$ &
$\mathfrak{m}^{(\mathcal{S}')}_{a_3,\textnormal{sum}}\dfz\displaystyle{\min_{\substack{i,j\in\mathcal{S} \\ i\neq j}}\sum_{\substack{\ell, m\in\mathcal{S}' \\ \ell\neq m}}\bm{a}_{i\ell}\bm{a}_{\ell m}\bm{a}_{m j}}$ &
$\left.
\begin{array}{ll}
Np_N\,\mathfrak{M}^{(\mathcal{S}')}_{a_3,\textnormal{sum}}\\
\\
Np_N\,\mathfrak{m}^{(\mathcal{S}')}_{a_3,\textnormal{sum}}
\end{array}
\right\}
\stackrel{\textnormal{p}}{\longrightarrow}\kappa^3 p (1-\xi)^2$\\
\hline
    \end{tabular}
\vspace*{5pt}
    \caption{Random variables and convergences relevant for the proofs of the theorems.}
        \label{tab:Mdefs}
  \end{center}
\end{table*}

\section{Useful Recursion on Matrix Powers}
\label{app:Matrecursions}

The next two theorems establish some useful properties of the powers of matrix $\bm{A}$, and of the powers of matrix $\bm{C}\dfz[\bm{A}^2]_{\mathcal{S}'}$ that has been introduced in~(\ref{eq:Hmatdef}). 
In particular, the theorems will provide upper and lower bounds for the entries of the matrix powers in terms of two appropriately constructed stochastic processes. 
These upper and lower bounds will be useful to examine the concentration behavior of the estimators proposed in this work, namely, Theorem~\ref{TheoremUsefulRecursionGranger} will be useful for the Granger estimator, whereas Theorem~\ref{TheoremPowersA} will be useful for the one-lag and for the residual estimators.

In the following, the symbol $\bm{a}^{(k)}_{i j}$ denotes the $(i,j)$-th entry of the $k$-th matrix power $\bm{A}^k$, and $\mathcal{N}$ denotes the set of all nodes. 
Table~\ref{tab:Mdefs} lists some useful random variables that are necessary to state and prove the theorems. 

\begin{theorem}[Useful bounds on the sum of powers of $\bm{A}$]
\label{TheoremPowersA}
The combination matrix $\bm{A}$ fulfills the following bounds for all $i,j\in\mathcal{N}$. 
We have:
\beq
\underline{\bm{\Sigma}}^{\textnormal{(even)}}_\alpha
\leq
\sum_{h=1}^\infty \bm{a}^{(2h)}_{i i}
\leq
\overline{\bm{\Sigma}}^{\textnormal{(even)}}_\alpha,~
\underline{\bm{\Sigma}}^{\textnormal{(odd)}}_\alpha
\leq
\sum_{h=1}^\infty \bm{a}^{(2h+1)}_{i i}
\leq
\overline{\bm{\Sigma}}^{\textnormal{(odd)}}_\alpha,
\label{eq:selfsum}
\eeq
and for $i\neq j$:
\beq
\bm{a}_{ij}\,\underline{\bm{\Sigma}}^{\textnormal{(even)}}_\beta
+
\mathfrak{m}\,\underline{\bm{\Sigma}}^{\textnormal{(even)}}_\gamma
\leq
\sum_{h=1}^\infty \bm{a}^{(2h)}_{i j}
\leq
\bm{a}_{ij}\,\overline{\bm{\Sigma}}^{\textnormal{(even)}}_\beta
+
\mathfrak{M}\,\overline{\bm{\Sigma}}^{\textnormal{(even)}}_\gamma,
\label{eq:offdiagsumeven}
\eeq
and
\beq
\bm{a}_{ij}\,\underline{\bm{\Sigma}}^{\textnormal{(odd)}}_\beta
+
\mathfrak{m}\,\underline{\bm{\Sigma}}^{\textnormal{(odd)}}_\gamma
\leq
\sum_{h=1}^\infty \bm{a}^{(2h+1)}_{i j}
\leq
\bm{a}_{ij}\,\overline{\bm{\Sigma}}^{\textnormal{(odd)}}_\beta
+
\mathfrak{M}\,\overline{\bm{\Sigma}}^{\textnormal{(odd)}}_\gamma,
\label{eq:offdiagsumodd}
\eeq
where the quantities $\mathfrak{m}, \mathfrak{M}$ and all the ``$\bm{\Sigma}$'' are random variables that do not depend on $i$ or $j$. 

If $\bm{A}$ is a regular diffusion matrix of parameters $\rho$ and $\kappa$ according to Assumption~\ref{assum:regdiffmat}, then under the uniform concentration regime the following convergences in probability hold as $N\rightarrow\infty$ (in the following formulas $\zeta=\rho-\kappa$).
\beq
\begin{array}{lll}
\overline{\bm{\Sigma}}^{\textnormal{(even)}}_{\alpha} \textnormal{ and }\underline{\bm{\Sigma}}^{\textnormal{(even)}}_{\alpha}
\stackrel{\textnormal{p}}{\longrightarrow}
\displaystyle{\frac{\zeta^2}{1 - \zeta^2}},
\\
\overline{\bm{\Sigma}}^{\textnormal{(even)}}_{\beta} \textnormal{ and } \underline{\bm{\Sigma}}^{\textnormal{(even)}}_{\beta}
\stackrel{\textnormal{p}}{\longrightarrow}
\displaystyle{\frac{2\,\zeta}{(1 - \zeta^2)^2}},
\\
\overline{\bm{\Sigma}}^{\textnormal{(even)}}_{\gamma} \textnormal{ and } \underline{\bm{\Sigma}}^{\textnormal{(even)}}_{\gamma}
\stackrel{\textnormal{p}}{\longrightarrow}
\displaystyle{\frac{1 + \zeta^2 + 2\rho\,\zeta}{(1-\rho^2)(1 - \zeta^2)^2}}
\end{array}
\label{eq:SigmaEvenConv}
\eeq
\beq
\begin{array}{lll}
\overline{\bm{\Sigma}}^{\textnormal{(odd)}}_{\alpha} \textnormal{ and }\underline{\bm{\Sigma}}^{\textnormal{(odd)}}_{\alpha}
\stackrel{\textnormal{p}}{\longrightarrow}
\displaystyle{\frac{\zeta^3}{1 - \zeta^2}},
\\
\overline{\bm{\Sigma}}^{\textnormal{(odd)}}_{\beta} \textnormal{ and }\underline{\bm{\Sigma}}^{\textnormal{(odd)}}_{\beta}
\stackrel{\textnormal{p}}{\longrightarrow}
\displaystyle{\frac{3\,\zeta^2 - \zeta^4}{(1 - \zeta^2)^2}},
\\
\overline{\bm{\Sigma}}^{\textnormal{(odd)}}_{\gamma} \textnormal{ and }\underline{\bm{\Sigma}}^{\textnormal{(odd)}}_{\gamma}
\stackrel{\textnormal{p}}{\longrightarrow}
\displaystyle{\frac{\rho + \rho \, \zeta^2 + 2\,\zeta}{(1-\rho^2)(1 - \zeta^2)^2}}
\end{array}
\label{eq:SigmaOddConv}
\eeq
\end{theorem}

\begin{IEEEproof}
The proof of the theorem will be articulated in the following steps.
\begin{enumerate}
\item[i)]
The entries of the combination matrix $\bm{A}$ are bounded as follows:
\beq
\underline{\bm{\alpha}}_k
\leq 
\bm{a}^{(k)}_{i i}
\leq \overline{\bm{\alpha}}_k,
\label{eq:aiiclaim}
\eeq
and, for $i\neq j$:
\beq
\underline{\bm{\beta}}_k\,\bm{a}_{i j} + 
\underline{\bm{\gamma}}_k\,\mathfrak{m}
\leq
\bm{a}^{(k)}_{i j}
\leq
\overline{\bm{\beta}}_k\,\bm{a}_{i j} + 
\overline{\bm{\gamma}}_k\,\mathfrak{M},
\label{eq:aijclaim}
\eeq
where, for $k\geq 2$, the (random) sequences $\overline{\bm{\alpha}}_k$, $\underline{\bm{\alpha}}_k$, $\overline{\bm{\beta}}_k$, $\underline{\bm{\beta}}_k$, $\overline{\bm{\gamma}}_k$, and $\underline{\bm{\gamma}}_k$, are determined by the following recursions:
\beq
\begin{array}{lll}
\overline{\bm{\alpha}}_{k+1}=\mathfrak{M}_{a,\textnormal{self}}\,\overline{\bm{\alpha}}_k  
+
\mathfrak{M}_a\,\rho^{k},
\\
\\
\overline{\bm{\beta}}_{k+1}=\overline{\bm{\alpha}}_k + \mathfrak{M}_{a,\textnormal{self}}\,\overline{\bm{\beta}}_k, 
\\
\\
\overline{\bm{\gamma}}_{k+1}=\overline{\bm{\beta}}_k 
+
\mathfrak{M}_{a,\textnormal{sum}}\,\overline{\bm{\gamma}}_k
\end{array}
\label{eq:alphabetagamma}
\eeq
with the initialization choices:
\beq
\overline{\bm{\alpha}}_2=\mathfrak{M}_{a_2,\textnormal{self}},\quad
\overline{\bm{\beta}}_2= 2\,\mathfrak{M}_{a,\textnormal{self}},
\quad 
\overline{\bm{\gamma}}_2=1,
\label{eq:UBinit2}
\eeq
and
\beq
\begin{array}{lll}
\underline{\bm{\alpha}}_{k+1}=\mathfrak{m}_{a,\textnormal{self}}\, \underline{\bm{\alpha}}_k, 
\\
\\
\underline{\bm{\beta}}_{k+1}=\underline{\bm{\alpha}}_k + \mathfrak{m}_{a,\textnormal{self}}\,\underline{\bm{\beta}}_k, 
\\
\\
\underline{\bm{\gamma}}_{k+1}=\underline{\bm{\beta}}_k 
+
\mathfrak{m}_{a,\textnormal{sum}}\,\underline{\bm{\gamma}}_k
\end{array}
\label{eq:alphabetagammamin}
\eeq
with the initialization choices:
\beq
\underline{\bm{\alpha}}_2=\mathfrak{m}_{a_2,\textnormal{self}},\quad
\underline{\bm{\beta}}_2= 2\,\mathfrak{m}_{a,\textnormal{self}},
\quad 
\underline{\bm{\gamma}}_2=1.
\label{eq:LBinit2}
\eeq
\item[ii)]
Let us introduce the following series:
\beqa
\overline{\bm{\Sigma}}^{\textnormal{(even)}}_\alpha&\dfz&
\sum_{h=1}^\infty \overline{\bm{\alpha}}_{2h},
~~~~
\underline{\bm{\Sigma}}^{\textnormal{(even)}}_\alpha\dfz
\sum_{h=1}^\infty \underline{\bm{\alpha}}_{2h},
\nonumber\\
\overline{\bm{\Sigma}}^{\textnormal{(even)}}_\beta&\dfz&
\sum_{h=1}^\infty \overline{\bm{\beta}}_{2h},
~~~~
\underline{\bm{\Sigma}}^{\textnormal{(even)}}_\beta\dfz
\sum_{h=1}^\infty \underline{\bm{\beta}}_{2h},
\nonumber\\
\overline{\bm{\Sigma}}^{\textnormal{(even)}}_\gamma&\dfz&
\sum_{h=1}^\infty \overline{\bm{\gamma}}_{2h},
~~~~
\underline{\bm{\Sigma}}^{\textnormal{(even)}}_\gamma\dfz
\sum_{h=1}^\infty \underline{\bm{\gamma}}_{2h},
\label{eq:SigmaevenDef}
\eeqa
and
\beqa
\overline{\bm{\Sigma}}^{\textnormal{(odd)}}_\alpha&\dfz&
\sum_{h=1}^\infty \overline{\bm{\alpha}}_{2h+1},
~~~~
\underline{\bm{\Sigma}}^{\textnormal{(odd)}}_\alpha\dfz
\sum_{h=1}^\infty \underline{\bm{\alpha}}_{2h+1},
\nonumber\\
\overline{\bm{\Sigma}}^{\textnormal{(odd)}}_\beta&\dfz&
\sum_{h=1}^\infty \overline{\bm{\beta}}_{2h+1},
~~~~
\underline{\bm{\Sigma}}^{\textnormal{(odd)}}_\beta\dfz
\sum_{h=1}^\infty \underline{\bm{\beta}}_{2h+1},
\nonumber\\
\overline{\bm{\Sigma}}^{\textnormal{(odd)}}_\gamma&\dfz&
\sum_{h=1}^\infty \overline{\bm{\gamma}}_{2h+1},
~~~~
\underline{\bm{\Sigma}}^{\textnormal{(odd)}}_\gamma\dfz
\sum_{h=1}^\infty \underline{\bm{\gamma}}_{2h+1},\nonumber\\
\label{eq:SigmaoddDef}
\eeqa
with:
\beqa
\overline{\bm{\Sigma}}_\alpha&\dfz&
\overline{\bm{\Sigma}}^{\textnormal{(even)}}_\alpha
+
\overline{\bm{\Sigma}}^{\textnormal{(odd)}}_\alpha,
~~~~
\underline{\bm{\Sigma}}_\alpha\dfz
\underline{\bm{\Sigma}}^{\textnormal{(even)}}_\alpha
+
\underline{\bm{\Sigma}}^{\textnormal{(odd)}}_\alpha,
\nonumber\\
\overline{\bm{\Sigma}}_\beta
&\dfz&
\overline{\bm{\Sigma}}^{\textnormal{(even)}}_\beta
+
\overline{\bm{\Sigma}}^{\textnormal{(odd)}}_\beta,
~~~~
\underline{\bm{\Sigma}}_\beta
\dfz
\underline{\bm{\Sigma}}^{\textnormal{(even)}}_\beta
+
\underline{\bm{\Sigma}}^{\textnormal{(odd)}}_\beta,
\nonumber\\
\overline{\bm{\Sigma}}_\gamma
&\dfz&
\overline{\bm{\Sigma}}^{\textnormal{(even)}}_\gamma
+
\overline{\bm{\Sigma}}^{\textnormal{(odd)}}_\gamma,
~~~~
\underline{\bm{\Sigma}}_\gamma
\dfz
\underline{\bm{\Sigma}}^{\textnormal{(even)}}_\gamma
+
\underline{\bm{\Sigma}}^{\textnormal{(odd)}}_\gamma.\nonumber\\
\label{eq:SigmaOverallDef2}
\eeqa
We will first show that the following convergence takes place:
\beq
\begin{array}{lll}
\overline{\bm{\Sigma}}_{\alpha} \textnormal{ and } \underline{\bm{\Sigma}}_{\alpha}
&
\stackrel{\textnormal{p}}{\longrightarrow}
\displaystyle{\frac{\zeta^2}{1 - \zeta}},
\\
\\
\overline{\bm{\Sigma}}_{\beta} \textnormal{ and } \underline{\bm{\Sigma}}_{\beta}
&
\stackrel{\textnormal{p}}{\longrightarrow}
\displaystyle{
\frac{1 - (1 - \zeta)^2}
{(1 - \zeta)^2}
}, 
\\
\\
\overline{\bm{\Sigma}}_{\gamma} \textnormal{ and } \underline{\bm{\Sigma}}_{\gamma}
&
\stackrel{\textnormal{p}}{\longrightarrow}
\displaystyle{
\frac{1}
{(1 - \rho)(1 - \zeta)^2}
}
\end{array}
\label{eq:SigmaTotConv}
\eeq
and then prove~(\ref{eq:SigmaEvenConv}) and~(\ref{eq:SigmaOddConv}).

\end{enumerate}
We start by examining the relationships pertaining to the main diagonal terms, namely, Eq.~(\ref{eq:aiiclaim}). 
For $k=2$ the claim is trivially true with the values of $\overline{\bm{\alpha}}_2$ and $\underline{\bm{\alpha}}_2$ in~(\ref{eq:UBinit2}) and~(\ref{eq:LBinit2}), because of the definitions $3)$ and $4)$ in Table~\ref{tab:Mdefs}. We shall therefore reason by induction to prove that~(\ref{eq:aiiclaim}) holds for an arbitrary $k$. In particular, we assume the claim verified for $k$, and manage to prove that it is verified for $k+1$. 
To this aim, we start by writing the diagonal terms of the $(k+1)$-th matrix power as: 
\beq
\bm{a}^{(k+1)}_{i i}=
\sum_{\ell\in\mathcal{N}}
\bm{a}_{i \ell}\bm{a}^{(k)}_{\ell i}
=
\bm{a}_{i i}\bm{a}^{(k)}_{i i}
+
\sum_{\substack{\ell\in\mathcal{N} \\ \ell\neq i}}
\bm{a}_{i \ell}\bm{a}^{(k)}_{\ell i},
\eeq
and, hence:
\beq
\bm{a}_{i i}\bm{a}^{(k)}_{i i}
\leq
\bm{a}^{(k+1)}_{i i}
\leq
\bm{a}_{i i}\bm{a}^{(k)}_{i i} + \mathfrak{M}_a \, \rho^{k},
\label{eq:selfaii}
\eeq
where we have used the fact that $\sum_{\ell\in\mathcal{N}}\bm{a}_{\ell i}^{(k)}=\rho^k$ along with definition $1)$ in Table~\ref{tab:Mdefs}.
Since we have assumed that~(\ref{eq:aiiclaim}) is true for $k$, from~(\ref{eq:selfaii}) we also have:
\beq
\mathfrak{m}_{a,\textnormal{self}}\,\underline{\bm{\alpha}}_k
\leq
\bm{a}_{ii}^{(k+1)}
\leq
\mathfrak{M}_{a,\textnormal{self}}\,\overline{\bm{\alpha}}_k
+
\mathfrak{M}_{a}\,\rho^{k},
\eeq
from which we conclude that~(\ref{eq:aiiclaim}) holds true, with the sequences $\overline{\bm{\alpha}}_k$ and $\underline{\bm{\alpha}}_k$ obeying the recursions in~(\ref{eq:alphabetagamma}) and~(\ref{eq:alphabetagammamin}), respectively.

We switch to the proof of~(\ref{eq:aijclaim}). 
For all $i,j\in\mathcal{N}$, with $i\neq j$, we have:
\beqa
\bm{a}_{ij}^{(2)}&=&\sum_{\ell\in\mathcal{N}}\bm{a}_{i\ell}\bm{a}_{\ell j}=
(\bm{a}_{ii} + \bm{a}_{jj})\bm{a}_{ij} + \sum_{\substack{\ell\in\mathcal{N} \\ \ell\neq i,j}}\bm{a}_{i\ell}\bm{a}_{\ell j}\nonumber\\
&\leq& 2\,\mathfrak{M}_{a,\textnormal{self}}\,\bm{a}_{ij} + \mathfrak{M}
\leq
2\,\mathfrak{M}_{a,\textnormal{self}}\,\mathfrak{M}_a + \mathfrak{M},
\label{eq:fundamentalaij2}
\eeqa
where we have applied the definitions $1)$, $2)$ and $12)$ in Table~\ref{tab:Mdefs}. 
First, we observe that~(\ref{eq:fundamentalaij2}) implies that~(\ref{eq:aijclaim}) holds true in the case $k=2$, with the choices detailed in~(\ref{eq:UBinit}) and~(\ref{eq:LBinit}). 
Moreover, we have:
\beq
\bm{a}^{(k+1)}_{i j}=
\sum_{\ell\in\mathcal{N}}
\bm{a}_{i \ell}\bm{a}^{(k)}_{\ell j}
=
\bm{a}_{i j}\bm{a}^{(k)}_{j j} + 
\sum_{\substack{\ell\in\mathcal{N} \\ \ell\neq j}}
\bm{a}_{i \ell}\bm{a}^{(k)}_{\ell j}.
\eeq
Therefore, since~(\ref{eq:aiiclaim}) holds true, and assuming that~(\ref{eq:aijclaim}) holds for an arbitrary $k\geq 2$, we have:
\beqa
\bm{a}^{(k+1)}_{i j}
&\leq&
\overline{\bm{\alpha}}_k \bm{a}_{i j}
+
\overline{\bm{\beta}}_k
\sum_{\substack{\ell\in \mathcal{N} \\ \ell\neq j}}
\bm{a}_{i \ell}\bm{a}_{\ell j} +
\overline{\bm{\gamma}}_k
\,\mathfrak{M}
\sum_{\substack{\ell\in \mathcal{N} \\ \ell\neq j}}
\bm{a}_{\ell j}\nonumber\\
&\leq&
(\overline{\bm{\alpha}}_k + \overline{\bm{\beta}}_k\,\mathfrak{M}_{a,\textnormal{self}}) \bm{a}_{i j}
+
( \overline{\bm{\beta}}_k + \overline{\bm{\gamma}}_k\,\mathfrak{M}_{a,\textnormal{sum}} )\,\mathfrak{M},\nonumber\\
\label{eq:arecur1}
\eeqa
which shows that the rightmost inequality in~(\ref{eq:aijclaim}) holds with the sequences $\overline{\bm{\beta}}_k$ and $\overline{\bm{\gamma}}_k$ obeying~(\ref{eq:alphabetagamma}), and with the initialization choice in~(\ref{eq:UBinit2}).

Next we focus on proving part ii). 
First, we notice that all the involved series are convergent. We will explain this conclusion with reference to the upper bounding sequences, with the case of the lower bounding sequences being dealt with similarly. 
The system of recursions in~(\ref{eq:alphabetagamma}) can be solved by calculating first $\overline{\bm{\alpha}}_k$, then $\overline{\bm{\beta}}_k$ (after substituting $\overline{\bm{\alpha}}_k$) and finally $\overline{\bm{\gamma}}_k$ (after substituting $\overline{\bm{\beta}}_k$). Applying Lemma~\ref{lem:lemmaseries}, it can be verified that all the obtained solutions are linear combinations of geometric sequences with ratio strictly smaller than one, from which convergence of the series $\overline{\bm{\Sigma}}_\alpha$, $\overline{\bm{\Sigma}}_\beta$ and $\overline{\bm{\Sigma}}_\gamma$ follows.

First of all, it is convenient to rewrite the series in~(\ref{eq:SigmaOverallDef2}) in the following more explicit form:
\beqa
\overline{\bm{\Sigma}}_\alpha&\dfz&
\sum_{k=2}^\infty \overline{\bm{\alpha}}_k,
~~~~
\underline{\bm{\Sigma}}_\alpha\dfz
\sum_{k=2}^\infty \underline{\bm{\alpha}}_k,
\nonumber\\
\overline{\bm{\Sigma}}_\beta&\dfz&
\sum_{k=2}^\infty \overline{\bm{\beta}}_k,
~~~~
\underline{\bm{\Sigma}}_\beta\dfz
\sum_{k=2}^\infty \underline{\bm{\beta}}_k,
\nonumber\\
\overline{\bm{\Sigma}}_\gamma&\dfz&
\sum_{k=2}^\infty \overline{\bm{\gamma}}_k,
~~~~
\underline{\bm{\Sigma}}_\gamma\dfz
\sum_{k=2}^\infty \underline{\bm{\gamma}}_k.
\label{eq:SigmasExplicit}
\eeqa
Let us consider the first line in~(\ref{eq:alphabetagamma}). By summing over index $k$, and using the definition of $\overline{\bm{\Sigma}}_{\alpha}$ in~(\ref{eq:SigmasExplicit}), we can write: 
\beq
\overline{\bm{\Sigma}}_{\alpha}
=
\overline{\bm{\alpha}}_2 + 
\mathfrak{M}_{a,\textnormal{self}}\,\overline{\bm{\Sigma}}_{\alpha} + \frac{\mathfrak{M}_a\, \rho^2}{1-\rho}
\Rightarrow
\overline{\bm{\Sigma}}_{\alpha}=\frac{\mathfrak{M}_{a_2,\textnormal{self}} + \bm{\epsilon}}{1 - \mathfrak{M}_{a,\textnormal{self}}}.
\label{eq:Sigmaa0}
\eeq
where we have set $\bm{\epsilon}=\frac{\mathfrak{M}_a\,\rho^2}{1-\rho}$.
Likewise, operating on the second line in~(\ref{eq:alphabetagamma}), and using the definition of $\overline{\bm{\Sigma}}_{\beta}$ in~(\ref{eq:SigmasExplicit}), we get:
\beq
\overline{\bm{\Sigma}}_{\beta}
=
\overline{\bm{\beta}}_2 
+ 
\overline{\bm{\Sigma}}_{\alpha} 
+ 
\mathfrak{M}_{a,\textnormal{self}}\,\overline{\bm{\Sigma}}_{\beta}
\Rightarrow
\overline{\bm{\Sigma}}_{\beta}=\frac{2\,\mathfrak{M}_{a,\textnormal{self}} + \overline{\bm{\Sigma}}_{\alpha}}
{1 - \mathfrak{M}_{a,\textnormal{self}}}.
\label{eq:Sigmab0}
\eeq
Finally, repeating the above procedure over the third line of~(\ref{eq:alphabetagamma}), we obtain:
\beq
\overline{\bm{\Sigma}}_{\gamma}=\overline{\bm{\gamma}}_2 + 
\overline{\bm{\Sigma}}_{\beta}
+
\mathfrak{M}_{a,\textnormal{sum}}\,\overline{\bm{\Sigma}}_{\gamma}
\Rightarrow
\overline{\bm{\Sigma}}_{\gamma}=\frac{1 + \overline{\bm{\Sigma}}_{\beta}}{1 - \mathfrak{M}_{a,\textnormal{sum}}},
\eeq
which, using~(\ref{eq:Sigmaa0}) and~(\ref{eq:Sigmab0}), after straightforward algebra yields:
\beq
\overline{\bm{\Sigma}}_{\gamma}=\frac{1 - \mathfrak{M}_{a,\textnormal{self}}^2 + \mathfrak{M}_{a_2,\textnormal{self}} + \bm{\epsilon}}{(1 - \mathfrak{M}_{a,\textnormal{sum}})(1 - \mathfrak{M}_{a,\textnormal{self}})^2}.
\label{eq:Sigmag0}
\eeq
Now, in view of the convergences in probability proved in Lemma~\ref{lem:listofconv} --- specifically, in view of~(\ref{eq:Maself}),~(\ref{eq:Ma2self}),~(\ref{eq:masum}) --- it is legitimate to replace the pertinent variables in~(\ref{eq:Sigmaa0}),~(\ref{eq:Sigmab0}), and~(\ref{eq:Sigmag0}), with their limits (that are reported in Table~\ref{tab:Mdefs}). 
After some lengthy, though straightforward, algebra, this replacement leads to~(\ref{eq:SigmaTotConv}).

We now move on to prove~(\ref{eq:SigmaOddConv}).  
Using the definitions of $\overline{\bm{\Sigma}}_{\alpha}^{\textnormal{(even)}}$ and $\overline{\bm{\Sigma}}_{\alpha}^{\textnormal{(odd)}}$ in~(\ref{eq:SigmaevenDef}) and~(\ref{eq:SigmaoddDef}), and summing over $k$ the first line in~(\ref{eq:alphabetagamma}), we see that:
\beqa
\overline{\bm{\Sigma}}_{\alpha}^{\textnormal{(odd)}}
&=&\sum_{\substack{k=3 \\k\textnormal{ odd}}}^{\infty}\overline{\bm{\alpha}}_{k}=
\sum_{\substack{k=2 \\ k\textnormal{ even}}}^{\infty}\overline{\bm{\alpha}}_{k+1}\nonumber\\
&=&
\mathfrak{M}_{a,\textnormal{self}}\,\sum_{\substack{k=2 \\ k\textnormal{ even}}}^{\infty}\overline{\bm{\alpha}}_{k}
+
\mathfrak{M}_a\,\sum_{\substack{k=2 \\ k\textnormal{ even}}}^{\infty}\rho^k\nonumber\\
&=&
\mathfrak{M}_{a,\textnormal{self}}\,\overline{\bm{\Sigma}}_{\alpha}^{\textnormal{(even)}}
+
\mathfrak{M}_a\,\sum_{k=1}^{\infty}\rho^{2 k}\nonumber\\
&=&
\mathfrak{M}_{a,\textnormal{self}}\,\overline{\bm{\Sigma}}_{\alpha}
-
\mathfrak{M}_{a,\textnormal{self}}\,\overline{\bm{\Sigma}}_{\alpha}^{\textnormal{(odd)}}
+
\mathfrak{M}_a\,\frac{\rho^2}{1-\rho^2},\nonumber\\ 
\eeqa
yielding:
\beq
\overline{\bm{\Sigma}}_{\alpha}^{\textnormal{(odd)}}=\displaystyle{
\frac{
\mathfrak{M}_{a,\textnormal{self}}\,\overline{\bm{\Sigma}}_{\alpha}
+
\bm{\epsilon}'}
{1 + \mathfrak{M}_{a,\textnormal{self}}}
}
\label{eq:SigmaaoddFin}
\eeq
where we defined $\bm{\epsilon}'=\mathfrak{M}_a\,\frac{\rho^2}{1-\rho^2}$.
Using now~(\ref{eq:SigmaOverallDef2}) and~(\ref{eq:SigmaaoddFin}), we further obtain:
\beq
\overline{\bm{\Sigma}}_{\alpha}^{\textnormal{(even)}}=\displaystyle{
\frac{
\overline{\bm{\Sigma}}_{\alpha}
-
\bm{\epsilon}'}
{1 + \mathfrak{M}_{a,\textnormal{self}}}
}
\label{eq:SigmaaevenFin}
\eeq
Proceeding in a similar way, from the second relationship in~(\ref{eq:alphabetagamma}), we obtain:
\beqa
\overline{\bm{\Sigma}}_{\beta}^{\textnormal{(odd)}}
&=&
\overline{\bm{\Sigma}}_{\alpha}^{\textnormal{(even)}}
+
\mathfrak{M}_{a,\textnormal{self}}\,\overline{\bm{\Sigma}}_{\beta}^{\textnormal{(even)}}
\nonumber\\
&=&
\overline{\bm{\Sigma}}_{\alpha}^{\textnormal{(even)}}
+
\mathfrak{M}_{a,\textnormal{self}}\,\overline{\bm{\Sigma}}_{\beta}
-
\mathfrak{M}_{a,\textnormal{self}}\,\overline{\bm{\Sigma}}_{\beta}^{\textnormal{(odd)}},
\eeqa
yielding:
\beq
\overline{\bm{\Sigma}}_{\beta}^{\textnormal{(odd)}}=\displaystyle{
\frac{\mathfrak{M}_{a,\textnormal{self}}\,\overline{\bm{\Sigma}}_{\beta}
+
\overline{\bm{\Sigma}}_{\alpha}^{\textnormal{(even)}}
}
{1 + \mathfrak{M}_{a,\textnormal{self}}}
}
\label{eq:SigmaboddFin}
\eeq
and
\beq
\overline{\bm{\Sigma}}_{\beta}^{\textnormal{(even)}}=
\displaystyle{
\frac{\overline{\bm{\Sigma}}_{\beta} - \overline{\bm{\Sigma}}_{\alpha}^{\textnormal{(even)}}
}
{1 + \mathfrak{M}_{a,\textnormal{self}}}
}
\label{eq:SigmabevenFin}
\eeq
Finally, from the third relationship in~(\ref{eq:alphabetagamma}), we can write:
\beqa
\overline{\bm{\Sigma}}_{\gamma}^{\textnormal{(odd)}}
&=&
\overline{\bm{\Sigma}}_{\beta}^{\textnormal{(even)}}
+
\mathfrak{M}_{a,\textnormal{sum}}\,\overline{\bm{\Sigma}}_{\gamma}^{\textnormal{(even)}}
\nonumber\\
&=&
\overline{\bm{\Sigma}}_{\beta}^{\textnormal{(even)}}
+
\mathfrak{M}_{a,\textnormal{sum}}\,\overline{\bm{\Sigma}}_{\gamma}
-
\mathfrak{M}_{a,\textnormal{sum}}\,\overline{\bm{\Sigma}}_{\beta}^{\textnormal{(odd)}},
\eeqa
yielding:
\beq
\overline{\bm{\Sigma}}_{\gamma}^{\textnormal{(odd)}}=\displaystyle{
\frac{\mathfrak{M}_{a,\textnormal{sum}}\,\overline{\bm{\Sigma}}_{\gamma}
+
\overline{\bm{\Sigma}}_{\beta}^{\textnormal{(even)}}
}
{1 + \mathfrak{M}_{a,\textnormal{sum}}}
}
\label{eq:SigmagoddFin}
\eeq
and
\beq
\overline{\bm{\Sigma}}_{\gamma}^{\textnormal{(even)}}=
\displaystyle{
\frac{\overline{\bm{\Sigma}}_{\gamma} - \overline{\bm{\Sigma}}_{\beta}^{\textnormal{(even)}}
}
{1 + \mathfrak{M}_{a,\textnormal{sum}}}
}
\label{eq:SigmagevenFin}
\eeq
The limiting results in~(\ref{eq:SigmaEvenConv}) and~(\ref{eq:SigmaOddConv}) are now obtained by replacing, in~(\ref{eq:SigmaaoddFin}),~(\ref{eq:SigmaaevenFin}),~(\ref{eq:SigmaboddFin}),~(\ref{eq:SigmabevenFin}),~(\ref{eq:SigmagoddFin}), and~(\ref{eq:SigmagevenFin}), the pertinent random variables with their limiting counterparts shown in Table~\ref{tab:Mdefs}. 
\end{IEEEproof}

\begin{theorem}[Useful bounds on the matrix $\bm{H}$]
\label{TheoremUsefulRecursionGranger}
The matrix $\bm{H}$ in~(\ref{eq:Hmatdef}) fulfills the following bounds for all $\ell,m\in\mathcal{S}'$. 
We have:
\beq
1 + \underline{\bm{\Phi}}_{\alpha}\leq\bm{h}_{\ell \ell}\leq 1 + \overline{\bm{\Phi}}_{\alpha},
\label{eq:hself}
\eeq
and for $\ell\neq m$:
\beq
\underline{\bm{\Phi}}_{\beta}\,\bm{a}_{\ell m} +
\underline{\bm{\Phi}}_{\gamma}\leq\bm{h}_{\ell m}
\leq 
\overline{\bm{\Phi}}_{\beta}\, \bm{a}_{\ell m} 
+ 
\overline{\bm{\Phi}}_{\gamma}.
\label{eq:hcross}
\eeq
If $\bm{A}$ is a regular diffusion matrix of parameters $\rho$ and $\kappa$ according to Assumption~\ref{assum:regdiffmat}, then under the uniform concentration regime the following convergences in probability hold as $N\rightarrow\infty$ (in the following formulas $\zeta=\rho-\kappa$).
\beq
\begin{array}{lll}
\overline{\bm{\Phi}}_{\alpha} \textnormal{ and } \underline{\bm{\Phi}}_{\alpha}
\stackrel{\textnormal{p}}{\longrightarrow}
\displaystyle{
\frac{\zeta^2}{1 - \zeta^2}
},
\\
\\
\overline{\bm{\Phi}}_{\beta}\textnormal{ and } \underline{\bm{\Phi}}_{\beta}
\stackrel{\textnormal{p}}{\longrightarrow}\displaystyle{
\frac{2\,\zeta}
{(1 - \zeta^2)^2}
},
\\
\\
N p_N \overline{\bm{\Phi}}_{\gamma} \textnormal{ and } N p_N\underline{\bm{\Phi}}_{\gamma}
\stackrel{\textnormal{p}}{\longrightarrow} \phi
\end{array}
\label{eq:SigmaGrangerConv}
\eeq
where:
\beq
\phi\dfz
\kappa^2 p\,\frac{1 - \zeta^2
+ 
2\,\zeta [ 2\,\zeta\,(1 - \xi) + \kappa (1 - \xi) ]
}
{[1 -  (\rho^2 - 2\rho\kappa\xi + \kappa^2 \xi) ][1 - \zeta^2]^2}.
\eeq
\end{theorem}

\begin{IEEEproof}
The proof of the theorem is articulated in the following steps.
\begin{enumerate}
\item[i)]
The entries of the matrix $\bm{C}^k$ are bounded as follows:
\beq
\underline{\bm{\alpha}}'_k
\leq 
\bm{c}^{(k)}_{\ell \ell}
\leq \overline{\bm{\alpha}}'_k,
\label{eq:cllclaim}
\eeq
and, for $\ell\neq m$:
\beq
\underline{\bm{\beta}}'_k \, \bm{a}_{\ell m} + \underline{\bm{\gamma}}'_k
\leq 
\bm{c}^{(k)}_{\ell m}
\leq
\overline{\bm{\beta}}'_k \, \bm{a}_{\ell m} + \overline{\bm{\gamma}}'_k,
\label{eq:clmclaim}
\eeq
where, for $k\geq 1$, the (random) sequences $\overline{\bm{\alpha}}'_k$, $\underline{\bm{\alpha}}'_k$, $\overline{\bm{\beta}}'_k$, $\underline{\bm{\beta}}'_k$, $\overline{\bm{\gamma}}'_k$, and $\underline{\bm{\gamma}}'_k$, are determined by the following recursions:
\beq
\begin{array}{lll}
\overline{\bm{\alpha}}'_{k+1}&\!\!\!\!\!=\mathfrak{M}_{c,\textnormal{self}}\,\overline{\bm{\alpha}}'_k 
+
(2\,\mathfrak{M}_{a,\textnormal{self}}\,\mathfrak{M}_a + \mathfrak{M} )\,\rho^{2 k},
\\
\\
\overline{\bm{\beta}}'_{k+1}&\!\!\!\!\!=2\,\mathfrak{M}_{a,\textnormal{self}}\,\overline{\bm{\alpha}}'_k + \mathfrak{M}_{c,\textnormal{self}}\,\overline{\bm{\beta}}'_k, 
\\
\\
\overline{\bm{\gamma}}'_{k+1}&\!\!\!\!\!=\mathfrak{M} \, \overline{\bm{\alpha}}'_k
+(2\,\mathfrak{M}_{a,\textnormal{self}}\,\mathfrak{M}^{(\mathcal{S}')} + \mathfrak{M} \, \mathfrak{M}^{(\mathcal{S}')}_{a,\textnormal{sum}})\overline{\bm{\beta}}'_k
\\
\\
&\!\!\!\!\!+
\mathfrak{M}_{c,\textnormal{sum}}\,\overline{\bm{\gamma}}'_k
\end{array}
\label{eq:GrangerUpperRecursion}
\eeq
with the initialization choices:
\beq
\overline{\bm{\alpha}}'_1=\mathfrak{M}_{c,\textnormal{self}},\quad
\overline{\bm{\beta}}'_1= 2\,\mathfrak{M}_{a,\textnormal{self}},
\quad 
\overline{\bm{\gamma}}'_1=\mathfrak{M},
\label{eq:UBinit}
\eeq
and
\beq
\begin{array}{lll}
\underline{\bm{\alpha}}'_{k+1}&\!\!\!\!\!=\mathfrak{m}_{c,\textnormal{self}}\,\underline{\bm{\alpha}}'_k,
\\
\\
\underline{\bm{\beta}}'_{k+1}&\!\!\!\!\!=2\,\mathfrak{m}_{a,\textnormal{self}}\,\underline{\bm{\alpha}}'_k + \mathfrak{m}_{c,\textnormal{self}}\,\underline{\bm{\beta}}'_k, 
\\
\\
\underline{\bm{\gamma}}'_{k+1}&\!\!\!\!\!=\mathfrak{m}\,\underline{\bm{\alpha}}'_k 
+ (2\,\mathfrak{m}_{a,\textnormal{self}}\,\mathfrak{m}^{(\mathcal{S}')} + \mathfrak{m}\,\mathfrak{m}^{(\mathcal{S}')}_{a,\textnormal{sum}})\underline{\bm{\beta}}'_k
\\
\\
&\!\!\!\!\!+
\mathfrak{m}_{c,\textnormal{sum}}\,\underline{\bm{\gamma}}'_k
\end{array}
\label{eq:GrangerLowerRecursion}
\eeq
with the initialization choices:
\beq
\underline{\bm{\alpha}}'_1=\mathfrak{m}_{c,\textnormal{self}},\quad
\underline{\bm{\beta}}'_1= 2\,\mathfrak{m}_{a,\textnormal{self}},
\quad 
\underline{\bm{\gamma}}'_1=\mathfrak{m}.
\label{eq:LBinit}
\eeq
\item[ii)] 
Let us introduce the following series:
\beqa
\overline{\bm{\Phi}}_\alpha&\dfz&
\sum_{k=1}^\infty \overline{\bm{\alpha}}'_k,
~~~~
\underline{\bm{\Phi}}_\alpha\dfz
\sum_{k=1}^\infty \underline{\bm{\alpha}}'_k,
\nonumber\\
\overline{\bm{\Phi}}_\beta&\dfz&
\sum_{k=1}^\infty \overline{\bm{\beta}}'_k,
~~~~
\underline{\bm{\Phi}}_\beta\dfz
\sum_{k=1}^\infty \underline{\bm{\beta}}'_k,
\nonumber\\
\overline{\bm{\Phi}}_\gamma&\dfz&
\sum_{k=1}^\infty \overline{\bm{\gamma}}'_k,
~~~~
\underline{\bm{\Phi}}_\gamma\dfz
\sum_{k=1}^\infty \underline{\bm{\gamma}}'_k.
\label{eq:Sigmasoveralldef}
\eeqa
\end{enumerate}

We start with the inequalities pertaining to the main diagonal of matrices $\bm{C}^k$, namely, with~(\ref{eq:cllclaim}). 
For $k=1$ we use the definitions of $\overline{\bm{\alpha}}'_1$ and $\underline{\bm{\alpha}}'_1$ in~(\ref{eq:UBinit}) and~(\ref{eq:LBinit}), respectively, to see that~(\ref{eq:cllclaim}) is trivially satisfied in view of definition $5)$ in Table~\ref{tab:Mdefs}. 
We shall now prove that~(\ref{eq:cllclaim}) holds for an arbitrary $k$ by induction. Assume thus that~(\ref{eq:cllclaim}) is true for $k$, we need to show that it is true for $k+1$.
From the definition of matrix $\bm{C}$ in~(\ref{eq:Hmatdef}), we can write its terms on the main diagonal as:
\beq
\bm{c}^{(k+1)}_{\ell \ell}=
\sum_{h\in\mathcal{S}'}
\bm{c}_{\ell h}\bm{c}^{(k)}_{h \ell}
=
\bm{c}_{\ell \ell}\bm{c}^{(k)}_{\ell \ell}
+
\sum_{\substack{h\in\mathcal{S}' \\ h\neq \ell}}
\bm{c}_{\ell h}\bm{c}^{(k)}_{h \ell}.
\eeq
In view of~(\ref{eq:fundamentalaij2}) we can write:
\beq
\sum_{\substack{h\in\mathcal{S}' \\ h\neq\ell}}\bm{c}_{\ell h}\bm{c}^{(k)}_{h \ell}=
\sum_{\substack{h\in\mathcal{S}' \\ h\neq\ell}}\bm{a}^{(2)}_{\ell h}\bm{c}^{(k)}_{h \ell} 
\leq
(2\,\mathfrak{M}_{a,\textnormal{self}}\,\mathfrak{M}_a + \mathfrak{M})\sum_{\substack{h\in\mathcal{S}' \\ h\neq\ell}} \bm{c}^{(k)}_{h \ell}.
\eeq
Moreover, we observe that:
\beq
\sum_{\substack{h\in\mathcal{S}' \\ h\neq\ell}}\bm{c}_{h \ell}\leq
\sum_{h\in\mathcal{N}}\bm{c}_{h \ell}=
\sum_{h\in\mathcal{N}}\bm{a}^{(2)}_{h \ell}=\rho^{2k},
\eeq
where the last equality follows from the first relationship in~(\ref{eq:fundamentalassum}) and from the symmetry of matrix $\bm{A}$. 
Therefore, we conclude that:
\beq
\mathfrak{m}_{c,\textnormal{self}}\,\bm{c}^{(k)}_{\ell \ell}
\leq
\bm{c}^{(k+1)}_{\ell \ell}
\leq
\mathfrak{M}_{c,\textnormal{self}}\,\bm{c}^{(k)}_{\ell \ell} + (2\,\mathfrak{M}_{a,\textnormal{self}}\,\mathfrak{M}_a + \mathfrak{M})\, \rho^{2 k},
\label{eq:cellell}
\eeq
Since we have assumed that~(\ref{eq:cllclaim}) holds true for $k$, we can further apply~(\ref{eq:cllclaim})  into~(\ref{eq:cellell}), yielding:
\beq
\mathfrak{m}_{c,\textnormal{self}}\,\underline{\bm{\alpha}}'_k
\leq
\bm{c}_{\ell\ell}^{(k+1)}\leq
\mathfrak{M}_{c,\textnormal{self}}\,\overline{\bm{\alpha}}'_k
+
(2\,\mathfrak{M}_{a,\textnormal{self}}\,\mathfrak{M}_a + \mathfrak{M})\,\rho^{2 k},
\eeq
which reveals that~(\ref{eq:cllclaim}) holds true, with the sequences $\overline{\bm{\alpha}}'_k$ and $\underline{\bm{\alpha}}'_k$ obeying the recursions in~(\ref{eq:GrangerUpperRecursion}) and~(\ref{eq:GrangerLowerRecursion}), respectively.

Let us move on to examine the case $\ell\neq m$. For all $\ell,m\in\mathcal{S}'$, with $\ell\neq m$, we can use~(\ref{eq:fundamentalaij2}) to conclude that:
\beq
\bm{c}_{\ell m}
\leq
2\,\mathfrak{M}_{a,\textnormal{self}}\,\bm{a}_{\ell m} + \mathfrak{M}.
\label{eq:cellm}
\eeq
Equation~(\ref{eq:cellm}) shows that the upper bound in~(\ref{eq:clmclaim}) holds for $k=1$, with the choices in~(\ref{eq:UBinit}). 
Now, rewriting the relationship $\bm{C}^{k+1}=\bm{C}\bm{C}^k$ on an entrywise basis, we have:   
\beqa
\bm{c}^{(k+1)}_{\ell m}&=&
\sum_{h\in\mathcal{S}'}
\bm{c}_{\ell h}\bm{c}^{(k)}_{h m}
\nonumber\\
&=&
\bm{c}_{\ell m}\bm{c}^{(k)}_{m m} + 
\sum_{\substack{h\in\mathcal{S}' \\ h\neq m}}
\bm{c}_{\ell h}\bm{c}^{(k)}_{h m}.
\label{eq:clmkplus1}
\eeqa
Accordingly, if we assume that~(\ref{eq:clmclaim}) holds true for an arbitrary $k$, from~(\ref{eq:clmkplus1}) we get:  
\beqa
\bm{c}^{(k+1)}_{\ell m}
&\leq&
\overline{\bm{\alpha}}'_k (2\,\mathfrak{M}_{a,\textnormal{self}}\,\bm{a}_{\ell m} + \mathfrak{M})
\nonumber\\
&+&
\overline{\bm{\beta}}'_k
\sum_{\substack{h\in\mathcal{S}' \\ h\neq m}}
\bm{c}_{\ell h}\bm{a}_{h m} + 
\overline{\bm{\gamma}}'_k
\sum_{\substack{h\in\mathcal{S}' \\ h\neq m}}
\bm{c}_{\ell h}.
\label{eq:crecur1}
\eeqa
Focusing on the first summation, using definition $5)$ in Table~\ref{tab:Mdefs} to bound the term $\bm{c}_{\ell}$, and~(\ref{eq:cellm}) to bound the term $\bm{c}_{\ell h}$, we get:
\beqa
\sum_{\substack{h\in\mathcal{S}' \\ h \neq m}}
\bm{c}_{\ell h}\bm{a}_{h m}&=&
\bm{c}_{\ell \ell}\bm{a}_{\ell m}
+
\sum_{\substack{h\in\mathcal{S}' \\ h \neq \ell, m}}
\bm{c}_{\ell h}\bm{a}_{h m}\nonumber\\
&\leq&
\mathfrak{M}_{c,\textnormal{self}}\,\bm{a}_{\ell m}
+
2\,\mathfrak{M}_{a,\textnormal{self}}\,
\sum_{\substack{h\in\mathcal{S}' \\ h \neq \ell, m}}
\bm{a}_{\ell h}\bm{a}_{h m}
\nonumber\\
&+& \mathfrak{M}\sum_{\substack{h\in\mathcal{S}' \\ h \neq \ell, m}}
\bm{a}_{h m}
\nonumber\\
&\leq&
\mathfrak{M}_{c,\textnormal{self}}\,\bm{a}_{\ell m}
+
2\,\mathfrak{M}_{a,\textnormal{self}}\,\mathfrak{M}^{(\mathcal{S}')}
+ \mathfrak{M}\,\mathfrak{M}^{(\mathcal{S}')}_{a,\textnormal{sum}},\nonumber\\
\label{eq:zoom}
\eeqa
where, in the latter inequality, we have further applied the definitions $7)$ and $13)$ in Table~\ref{tab:Mdefs}.
Using now~(\ref{eq:zoom}) into~(\ref{eq:crecur1}), we get:
\beqa
\lefteqn{\bm{c}_{\ell m}^{(k+1)}\leq
( 
\underbrace{
2\,\mathfrak{M}_{a,\textnormal{self}}\,\overline{\bm{\alpha}}'_k + \mathfrak{M}_{c,\textnormal{self}}\,\overline{\bm{\beta}}'_k 
}_{\overline{\bm{\beta}}'_{k+1}}
)
\,\bm{a}_{\ell m}}\nonumber\\
&&\!\!\!\!\!\!\!\!+
[
\underbrace{
\mathfrak{M}\,\overline{\bm{\alpha}}'_k 
+ (2\,\mathfrak{M}_{a,\textnormal{self}}\,\mathfrak{M}^{(\mathcal{S}')} + \mathfrak{M}\,\mathfrak{M}^{(\mathcal{S}')}_{a,\textnormal{sum}})\overline{\bm{\beta}}'_k
+
\mathfrak{M}_{c,\textnormal{sum}}\,\overline{\bm{\gamma}}'_k 
}_{\overline{\bm{\gamma}}'_{k+1}}
],
\nonumber\\
\eeqa
which shows that~(\ref{eq:clmclaim}) holds also for $k+1$, with $\overline{\bm{\beta}}'_k$, and $\overline{\bm{\gamma}}'_k$ obeying the recursion in~(\ref{eq:GrangerUpperRecursion}). The proof of~(\ref{eq:GrangerLowerRecursion}) is similar. 
This concludes part i) of the theorem. 
Part ii) comes readily as a corollary of part i), after noticing that the matrix $\bm{H}$ defined in~(\ref{eq:Hmatdef}) can be expressed as $\bm{H}=(I_{\mathcal{S}'} - \bm{C})^{-1}=I_{\mathcal{S}'} + \bm{C} + \bm{C}^2+\ldots$, and summing the inequalities in~(\ref{eq:cllclaim}) and~(\ref{eq:clmclaim}) over index $k$. Convergence of all the involved series is guaranteed by Lemma~\ref{lem:lemmaseries}. 

Next we focus on proving part iii). 
Preliminarily, we observe that all the involved series are convergent. We will explain this conclusion with reference to the upper bounding sequences, with the case of the lower bounding sequences being dealt with similarly. 
The system of recursions in~(\ref{eq:GrangerUpperRecursion}) can be solved by calculating first $\overline{\bm{\alpha}}'_k$, then $\overline{\bm{\beta}}'_k$ (after substituting $\overline{\bm{\alpha}}'_k$) and finally $\overline{\bm{\gamma}}'_k$ (after substituting $\overline{\bm{\beta}}'_k$). Applying Lemma~\ref{lem:lemmaseries}, it can be verified that all the obtained solutions are linear combinations of geometric sequences with ratio strictly smaller than one, from which convergence of the series $\overline{\bm{\Phi}}_\alpha$, $\overline{\bm{\Phi}}_\beta$ and $\overline{\bm{\Phi}}_\gamma$ follows.

Let us consider the first line in~(\ref{eq:GrangerUpperRecursion}). By summing over index $k$, and recalling the definition of $\overline{\bm{\Phi}}_{\alpha}$ in~(\ref{eq:Sigmasoveralldef}), we can write: 
\beq
\overline{\bm{\Phi}}_{\alpha}
=
\overline{\bm{\alpha}}'_1 + 
\mathfrak{M}_{c,\textnormal{self}}\,\overline{\bm{\Phi}}_{\alpha} + \bm{\epsilon}
\Rightarrow
\overline{\bm{\Phi}}_{\alpha}=\frac{\mathfrak{M}_{c,\textnormal{self}} + \bm{\epsilon}}{1 - \mathfrak{M}_{c,\textnormal{self}}},
\label{eq:Sigmaa}
\eeq
where we have set:
\beq
\bm{\epsilon}=(2\,\mathfrak{M}_{a,\textnormal{self}}\,\mathfrak{M}_a + \mathfrak{M})\,\frac{\rho^2}{1-\rho^2}.
\eeq
Likewise, from the second line in~(\ref{eq:GrangerUpperRecursion}), summing over index $k$, and recalling the definition of $\overline{\bm{\Phi}}_{\beta}$ in~(\ref{eq:Sigmasoveralldef}), we can write:
\beq
\overline{\bm{\Phi}}_{\beta}
=
\overline{\bm{\beta}}'_1 
+ 
2\,\mathfrak{M}_{a,\textnormal{self}}\,\overline{\bm{\Phi}}_{\alpha} 
+ 
\mathfrak{M}_{c,\textnormal{self}}\,\overline{\bm{\Phi}}_{\beta},
\eeq
yielding:
\beq
\overline{\bm{\Phi}}_{\beta}=\frac{2\,\mathfrak{M}_{a,\textnormal{self}}\,(1+ \overline{\bm{\Phi}}_{\alpha})}
{1 - \mathfrak{M}_{c,\textnormal{self}}}.
\label{eq:Sigmab}
\eeq
Using now~(\ref{eq:Sigmaa}) into~(\ref{eq:Sigmab}), we get:
\beq
\overline{\bm{\Phi}}_{\beta}
=
\frac{2\,\mathfrak{M}_{a,\textnormal{self}}}
{(1 - \mathfrak{M}_{c,\textnormal{self}})^2}(1+ \bm{\epsilon}).
\eeq
Finally, repeating the same procedure on the third line in~(\ref{eq:GrangerUpperRecursion}), we obtain:
\beqa
\overline{\bm{\Phi}}_{\gamma}&=&\overline{\bm{\gamma}}'_1 + 
\mathfrak{M}\,\overline{\bm{\Phi}}_{\alpha}
+ (2\,\mathfrak{M}_{a,\textnormal{self}}\,\mathfrak{M}^{(\mathcal{S}')} + \mathfrak{M}\,\mathfrak{M}^{(\mathcal{S}')}_{a,\textnormal{sum}}) \overline{\bm{\Phi}}_{\beta}\nonumber\\
&+&
\mathfrak{M}_{c,\textnormal{sum}}\,\overline{\bm{\Phi}}_{\gamma},
\eeqa
which yields the solution:
\beqa
\overline{\bm{\Phi}}_{\gamma}&=&
\frac{\mathfrak{M}(1 + \overline{\bm{\Phi}}_{\alpha})
+ (2\,\mathfrak{M}_{a,\textnormal{self}}\,\mathfrak{M}^{(\mathcal{S}')} 
+ \mathfrak{M}\,\mathfrak{M}^{(\mathcal{S}')}_{a,\textnormal{sum}}) \overline{\bm{\Phi}}_{\beta}}
{1 - \mathfrak{M}_{c,\textnormal{sum}}}.\nonumber\\
\label{eq:Sigmag}
\eeqa
Now, in view of the convergences in probability proved in Lemma~\ref{lem:listofconv} --- specifically, in view of~(\ref{eq:Maself}),~(\ref{eq:MaSprimesum}),~(\ref{eq:Mcself}),~(\ref{eq:Mcsum}) and~(\ref{eq:Msimply}) --- it is legitimate to replace the pertinent variables with their limits (that are reported in Table~\ref{tab:Mdefs}) in~(\ref{eq:Sigmaa}),~(\ref{eq:Sigmab}), and~(\ref{eq:Sigmag}), which, after some lengthy, though straightforward, algebra, leads to~(\ref{eq:SigmaGrangerConv}).
\end{IEEEproof}

\section{Proof of Theorem~\ref{theor1}}
\label{app:Theorem1}
We start by proving an auxiliary lemma. 

\begin{lemma}[Sufficient conditions for universal local structural consistency] 
\label{lem:suffcond}
Let the network graph be drawn according to an Erd\H{o}s-R\'enyi random graph model, and let $\bm{A}$ be a regular diffusion matrix with parameters $\rho$ and $\kappa$. 
Let $\mathcal{S}$ be the set of observable nodes and consider then an estimator:
\beq
\widehat{\bm{A}}_{\mathcal{S}}=\bm{A}_{\mathcal{S}}+\bm{E}.
\label{eq:errdecompappendix}
\eeq 
Assume that, for all $i,j\in\mathcal{S}$, with $i\neq j$:
\beq
\underline{\bm{w}}_N\, \bm{a}_{ij} + \underline{\bm{z}}_N
\leq
\bm{e}_{ij}
\leq
\overline{\bm{w}}_N\, \bm{a}_{ij} + \overline{\bm{z}}_N,
\label{eq:errgenbound0}
\eeq
where the quantities $\underline{\bm{w}}_N$, $\overline{\bm{w}}_N$, $\underline{\bm{z}}_N$ and $\overline{\bm{z}}_N$ do depend on the network size, $N$, but they do not depend on $(i,j)$, and fulfill the following convergences:
\beqa
&&\underline{\bm{w}}_N\stackrel{\textnormal{p}}{\longrightarrow}w,~~~~~~~~~~
\overline{\bm{w}}_N\stackrel{\textnormal{p}}{\longrightarrow}w,
\nonumber\\
&&N p_N \underline{\bm{z}}_N\stackrel{\textnormal{p}}{\longrightarrow}z,~~~~
~N p_N\overline{\bm{z}}_N\stackrel{\textnormal{p}}{\longrightarrow}z.
\label{eq:variousconvXY}
\eeqa
Then, under the uniform concentration regime, the estimator $\widehat{\bm{A}}_{\mathcal{S}}$ achieves universal local structural consistency, with a scaling sequence $s_N=N p_N$, bias $\eta=z$, and identifiability gap $\Gamma=\kappa\,(1 + w)$.
\end{lemma}
\begin{IEEEproof}
Using~(\ref{eq:upperlowerdisc}) and~(\ref{eq:errdecompappendix}) we can write:
\beq
\underline{\bm{\delta}}_N=\min_{\substack{i,j\in\mathcal{S}:\bm{a}_{ij}=0 \\ i\neq j}} \bm{e}_{ij},\qquad
\overline{\bm{\delta}}_N=\max_{\substack{i,j\in\mathcal{S}:\bm{a}_{ij}=0 \\ i\neq j}} \bm{e}_{ij},
\eeq
and, hence, from~(\ref{eq:errgenbound0}) we get:
\beq
N p_N \, \underline{\bm{z}}_N
\leq 
N p_N \,\underline{\bm{\delta}}_N
\leq 
N p_N \,\overline{\bm{\delta}}_N
\leq 
N p_N \, \overline{\bm{z}}_N.
\eeq
Using~(\ref{eq:variousconvXY}), we conclude that:
\beq
N p_N \,\overline{\bm{\delta}}_N\stackrel{\textnormal{p}}{\longrightarrow} z,\qquad
N p_N \,\underline{\bm{\delta}}_N\stackrel{\textnormal{p}}{\longrightarrow} z.
\eeq 
Let us now examine the connected pairs. 
From~(\ref{eq:errgenbound0}) we know that:
\beq
(1+\underline{\bm{w}}_N)\,\bm{a}_{ij} + \underline{\bm{z}}_N
\leq
\bm{a}_{ij} + \bm{e}_{ij}
\leq (1+\overline{\bm{w}}_N)\,\bm{a}_{ij} + \overline{\bm{z}}_N,
\eeq
which, used along with~(\ref{eq:upperlowerconn}) gives:
\beq
N p_N \underline{\bm{\Delta}}_N
\geq
(1+\underline{\bm{w}}_N)\,[N p_N \min_{\substack{i,j\in\mathcal{S}:\bm{a}_{ij}>0\\ i\neq j}} \bm{a}_{ij}] + N p_N \underline{\bm{z}}_N,
\eeq
and
\beq
N p_N \overline{\bm{\Delta}}_N
\leq
(1+\overline{\bm{w}}_N)\,[N p_N \max_{\substack{i,j\in\mathcal{S}:\bm{a}_{ij}>0\\ i\neq j}} \bm{a}_{ij}] + N p_N \overline{\bm{z}}_N.
\eeq
Now we observe that, in view of Assumption~\ref{assum:regdiffmat} we can write, for all pairs $(i,j)$ where $\bm{a}_{ij}>0$:
\beq
\kappa\,\frac{N p_N}{\bm{d}_{\max}}
\leq
N p_N \, \bm{a}_{ij}
\leq\kappa\,\frac{N p_N}{\bm{d}_{\min}}.
\label{eq:aijminbasic}
\eeq
Under the uniform concentration regime --- see~(\ref{eq:dmaxmin}) --- we conclude from~(\ref{eq:aijminbasic}) that:
\beq
N p_N \,\underline{\bm{\Delta}}_N\stackrel{\textnormal{p}}{\longrightarrow} \kappa\,(1 + w) + z,~
N p_N \,\overline{\bm{\Delta}}_N\stackrel{\textnormal{p}}{\longrightarrow} \kappa\,(1 + w) + z.
\eeq
It remains to apply the definition of bias and identifiability gap in~(\ref{eq:igapstrong}) to get the claim of the lemma.
\end{IEEEproof}

In order to prove Theorem~\ref{theor1}, it is necessary examine separately the Granger estimator and the other two estimators, namely, the one-lag and the residual estimators.

\vspace*{5pt}
\begin{IEEEproof}[Proof for the Granger estimator]
The proof for the Granger estimator boils down to combining Theorem~\ref{TheoremUsefulRecursionGranger} with Lemma~\ref{lem:suffcond}.

From~(\ref{eq:basicerrGra}) we can write, for $i,j\in\mathcal{S}$, with $i\neq j$:
\beqa
\bm{e}^{\textnormal{(Gra)}}_{ij}&=&
\sum_{\ell,m\in\mathcal{S}'}\bm{a}_{i\ell}\bm{h}_{\ell m}\bm{a}^{(2)}_{m j}
\nonumber\\
&=&
\sum_{\ell\in\mathcal{S}'}\bm{a}_{i\ell}\bm{h}_{\ell \ell}\bm{a}^{(2)}_{\ell j}
+
\sum_{\substack{\ell,m\in\mathcal{S}' \\ \ell\neq m}}\bm{a}_{i\ell}\bm{h}_{\ell m}\bm{a}^{(2)}_{m j}
\nonumber\\
&\leq&
(1 + \overline{\bm{\Phi}}_{\alpha})
\sum_{\ell\in\mathcal{S}'}\bm{a}_{i\ell}\bm{a}^{(2)}_{\ell j}
\label{eq:boundselfh}
\\
&+&
\overline{\bm{\Phi}}_{\beta} 
\sum_{\substack{\ell,m\in\mathcal{S}' \\ \ell\neq m}}\bm{a}_{i\ell}\bm{a}_{\ell m}\bm{a}^{(2)}_{m j} 
+ 
\overline{\bm{\Phi}}_{\gamma}
\sum_{\substack{\ell,m\in\mathcal{S}' \\ \ell\neq m}}\bm{a}_{i\ell}\bm{a}^{(2)}_{m j},
\label{eq:boundmutualh}
\nonumber
\\
\label{eq:ineqchainGra}
\eeqa
where the inequality is obtained by bounding the entries of the matrix $\bm{H}$, specifically, we have that~(\ref{eq:boundselfh}) follows from~(\ref{eq:hself}), whereas~(\ref{eq:boundmutualh}) follows from~(\ref{eq:hcross}).
Using now~(\ref{eq:fundamentalaij2}) in~(\ref{eq:ineqchainGra}), we get:
\beqa
\bm{e}^{\textnormal{(Gra)}}_{ij}&\leq&
(1 + \overline{\bm{\Phi}}_{\alpha})
\sum_{\ell\in\mathcal{S}'}\bm{a}_{i\ell}
(2\,\mathfrak{M}_{a,\textnormal{self}}\,\bm{a}_{\ell j} + \mathfrak{M})
\nonumber\\
&+&
\overline{\bm{\Phi}}_{\beta} 
\sum_{\substack{\ell,m\in\mathcal{S}' \\ \ell\neq m}}\bm{a}_{i\ell}\bm{a}_{\ell m}
(2\,\mathfrak{M}_{a,\textnormal{self}}\,\bm{a}_{m j} + \mathfrak{M})
\nonumber\\
&+& 
2\overline{\bm{\Phi}}_{\gamma}
\mathfrak{M}_{a,\textnormal{self}} \sum_{\substack{\ell,m\in\mathcal{S}' \\ \ell\neq m}}\bm{a}_{i\ell}\bm{a}_{m j}
\nonumber\\
&+&
\overline{\bm{\Phi}}_{\gamma}
\sum_{\substack{\ell,m\in\mathcal{S}' \\ \ell\neq m}}\bm{a}_{i\ell}
\sum_{\substack{\ell'\in\mathcal{N} \\ \ell'\neq m,j}}\bm{a}_{m\ell'}\bm{a}_{\ell' j}.
\eeqa
which can be recast in the following convenient form:
\beqa
\bm{e}^{\textnormal{(Gra)}}_{ij}
&\leq&
(1+\overline{\bm{\Phi}}_{\alpha})
\left[
2\,\mathfrak{M}_{a,\textnormal{self}}\,
\mathfrak{M}^{(\mathcal{S}')}\,
+
\mathfrak{M}\,
\widetilde{\mathfrak{M}}^{(\mathcal{S}')}_{a,\textnormal{sum}}
\right]
\nonumber\\
&+&
\overline{\bm{\Phi}}_{\beta}
\left[
2\,\mathfrak{M}_{a,\textnormal{self}}\, 
\mathfrak{M}^{(\mathcal{S}')}_{a_3,\textnormal{sum}}
+
\mathfrak{M}\,
\widetilde{\mathfrak{M}}^{(\mathcal{S}')}
\right]
\nonumber\\ 
&+& 
\overline{\bm{\Phi}}_{\gamma}
\left[
2\,\mathfrak{M}_{a,\textnormal{self}}\,
\widetilde{\widetilde{\mathfrak{M}}}^{(\mathcal{S}')}
+
\widetilde{\widetilde{\mathfrak{M}}}^{(\mathcal{S}')}_{a,\textnormal{sum}}
\right]
\dfz\overline{\bm{z}}_N,\nonumber\\
\label{eq:eijineq}
\eeqa
where we have used definitions $7)$, $8)$, $9)$, $10)$ and $14)$ listed in Table~\ref{tab:Mdefs}.
The arguments leading to this result can be repeated by replacing upper bounds with lower bounds, and maxima with minima (e.g., $\mathfrak{M}$ replaced by $\mathfrak{m}$, or $\overline{\bm{\Phi}}_\alpha$ replaced by $\underline{\bm{\Phi}}_\alpha$), yielding:
\beqa
\bm{e}^{\textnormal{(Gra)}}_{ij}
&\geq&
(1+\underline{\bm{\Phi}}_{\alpha})
\left[
2\,\mathfrak{m}_{a,\textnormal{self}}\,
\mathfrak{m}^{(\mathcal{S}')}\,
+
\mathfrak{m}\,
\widetilde{\mathfrak{m}}^{(\mathcal{S}')}_{a,\textnormal{sum}}
\right]
\nonumber\\
&+&
\underline{\bm{\Phi}}_{\beta}
\left[
2\,\mathfrak{m}_{a,\textnormal{self}}\, 
\mathfrak{m}^{(\mathcal{S}')}_{a_3,\textnormal{sum}}
+
\mathfrak{m}\,
\widetilde{\mathfrak{m}}^{(\mathcal{S}')}
\right]
\nonumber\\ 
&+& 
\underline{\bm{\Phi}}_{\gamma}
\left[
2\,\mathfrak{m}_{a,\textnormal{self}}\,
\widetilde{\widetilde{\mathfrak{m}}}^{(\mathcal{S}')}
+
\widetilde{\widetilde{\mathfrak{m}}}^{(\mathcal{S}')}_{a,\textnormal{sum}}
\right]\dfz\underline{\bm{z}}_N.\nonumber\\
\label{eq:eijineqbis}
\eeqa
Now, under the uniform concentration regime we can use the pertinent convergences in probability listed in Table~\ref{tab:Mdefs}, in conjunction with the convergences in~(\ref{eq:SigmaGrangerConv}), which, after some tedious but straightforward algebra, lead to:
\beq
N p_N\,\overline{\bm{z}}_N\stackrel{\textnormal{p}}{\longrightarrow}\eta,\qquad
N p_N\,\underline{\bm{z}}_N\stackrel{\textnormal{p}}{\longrightarrow}\eta
\label{eq:upperbounderrorGra}
\eeq
where $\eta$ is the bias corresponding to the Granger estimator in Table~\ref{tab:Theorem1}. 
Accordingly, we can conclude that the error for the Granger estimator fulfills the hypotheses of Lemma~\ref{lem:suffcond}, with the choice $\overline{\bm{w}}_N=\underline{\bm{w}}_N=0$, and with the quantities $\overline{\bm{z}}_N$ and $\underline{\bm{z}}_N$ defined in~(\ref{eq:eijineq}) and~(\ref{eq:eijineqbis}), respectively. This concludes the proof for the claim pertaining to the behavior of the Granger estimator under the uniform concentration regime.
\end{IEEEproof}

\begin{IEEEproof}[Proof for the one-lag and for the residual estimators]
The proof of the claim for the one-lag and for the residual estimators boils down to combining Theorem~\ref{TheoremPowersA} with Lemma~\ref{lem:suffcond}. 
In light of the definitions of bias and gap, it suffices to prove the claim with $\sigma^2=1$, and then scale the values obtained for the bias and the gap by an arbitrary $\sigma^2$ --- see also Remark~\ref{rem:scaleirrel}. 

Using~(\ref{eq:onelagesterrdef}), the error corresponding to the one-lag estimator can be written as, for all $i,j\in \mathcal{S}$, with $i\neq j$:
\beq
\bm{e}^{\textnormal{(1-lag)}}_{ij}=\sum_{h=1}^\infty \bm{a}^{(2h+1)}_{ij},
\eeq
and, hence, using the bounds in~(\ref{eq:offdiagsumodd}), we can write: 
\beq
\underline{\bm{\Sigma}}^{\textnormal{(odd)}}_\beta
\,\bm{a}_{ij}
+
\underline{\bm{\Sigma}}^{\textnormal{(odd)}}_\gamma\,\mathfrak{m}
\leq
\bm{e}^{\textnormal{(1-lag)}}_{ij}
\leq
\overline{\bm{\Sigma}}^{\textnormal{(odd)}}_\beta
\,\bm{a}_{ij}
+
\overline{\bm{\Sigma}}^{\textnormal{(odd)}}_\gamma\,\mathfrak{M}.
\eeq
Now we see that there are two contributions in the error. 
The first contributions ($\underline{\bm{\Sigma}}^{\textnormal{(odd)}}_\beta$ for the lower bound, and $\overline{\bm{\Sigma}}^{\textnormal{(odd)}}_\beta$ for the upper bound) multiply the entries of the combination matrix, $\bm{a}_{ij}$. Accordingly, they will play a role for connected agents.
The second contributions ($\underline{\bm{\Sigma}}^{\textnormal{(odd)}}_\gamma\,\mathfrak{m}$ for the lower bound, and $\overline{\bm{\Sigma}}^{\textnormal{(odd)}}_\gamma\,\mathfrak{M}$ for the upper bound) play a role for {\em all} agents, whether or note they are connected. 

Now, using the convergence results in~(\ref{eq:SigmaOddConv}) and in~(\ref{eq:Msimply}), simple algebraic calculations lead to:
\beq
N p_N \, \underline{\bm{\Sigma}}^{\textnormal{(odd)}}_\gamma\,\mathfrak{m}
\stackrel{\textnormal{p}}{\longrightarrow}\eta,
\quad
N p_N \, \overline{\bm{\Sigma}}^{\textnormal{(odd)}}_\gamma\,\mathfrak{M}
\stackrel{\textnormal{p}}{\longrightarrow}\eta,
\eeq 
where $\eta$ is equal to the bias of the one-lag estimator as defined in the pertinent row of Table~\ref{tab:Theorem1}. 
For what concerns $\underline{\bm{\Sigma}}^{\textnormal{(odd)}}_\beta$ and $\overline{\bm{\Sigma}}^{\textnormal{(odd)}}_\beta$, from~(\ref{eq:SigmaOddConv}) we see that: 
\beq
\underline{\bm{\Sigma}}^{\textnormal{(odd)}}_\beta
\stackrel{\textnormal{p}}{\longrightarrow}\Gamma/\kappa - 1,
\quad
\overline{\bm{\Sigma}}^{\textnormal{(odd)}}_\beta
\stackrel{\textnormal{p}}{\longrightarrow} \Gamma/\kappa - 1,
\eeq
where $\Gamma$ is equal to the bias of the one-lag estimator as defined in the pertinent row of Table~\ref{tab:Theorem1} (recall that we are working with $\sigma^2=1$).
It remains to apply Lemma~\ref{lem:suffcond}, with the choices:
\beqa
\overline{\bm{w}}_N&=&\overline{\bm{\Sigma}}^{\textnormal{(odd)}}_\beta,\qquad~~~
\underline{\bm{w}}_N=\underline{\bm{\Sigma}}^{\textnormal{(odd)}}_\beta,
\label{eq:WWdefs}\\
\overline{\bm{z}}_N&=&\overline{\bm{\Sigma}}^{\textnormal{(odd)}}_\gamma\,\mathfrak{M},\qquad
\underline{\bm{z}}_N=\underline{\bm{\Sigma}}^{\textnormal{(odd)}}_\gamma\,\mathfrak{m},
\label{eq:ZZdefs}
\eeqa
which concludes the proof for the one-lag estimator under the uniform concentration regime.

Let us switch to the analysis of the residual estimator. Using~(\ref{eq:residerrestdef}), the error corresponding to the residual estimator can be written as, for all $i,j\in \mathcal{S}$, with $i\neq j$:
\beq
\bm{e}^{\textnormal{(res)}}_{ij}=\sum_{h=1}^\infty \bm{a}^{(2h+1)}_{ij} - \sum_{h=1}^\infty \bm{a}^{(2h)}_{ij},
\eeq
and, hence, using the bounds in~(\ref{eq:offdiagsumeven}) and~(\ref{eq:offdiagsumodd}), we can write:
\beq
\bm{e}^{\textnormal{(res)}}_{ij}
\leq
(\overline{\bm{\Sigma}}^{\textnormal{(odd)}}_\beta - \underline{\bm{\Sigma}}^{\textnormal{(even)}}_\beta)
\,\bm{a}_{ij}
+
(\overline{\bm{\Sigma}}^{\textnormal{(odd)}}_\gamma\,\mathfrak{M} 
- 
\underline{\bm{\Sigma}}^{\textnormal{(even)}}_\gamma\,\mathfrak{m}
),
\label{eq:errorresrepr}
\eeq
and
\beq
\bm{e}^{\textnormal{(res)}}_{ij}
\geq
(\underline{\bm{\Sigma}}^{\textnormal{(odd)}}_\beta - \overline{\bm{\Sigma}}^{\textnormal{(even)}}_\beta)
\,\bm{a}_{ij}
+
(\underline{\bm{\Sigma}}^{\textnormal{(odd)}}_\gamma\,\mathfrak{m} 
- 
\overline{\bm{\Sigma}}^{\textnormal{(even)}}_\gamma\,\mathfrak{M}
).
\label{eq:errorresrepr2}
\eeq
Now, using the convergence results in~(\ref{eq:SigmaEvenConv}), (\ref{eq:SigmaOddConv}) and~(\ref{eq:Msimply}), simple algebraic calculations lead to:
\beqa
N p_N \, (\overline{\bm{\Sigma}}^{\textnormal{(odd)}}_\gamma\,\mathfrak{M} 
- 
\underline{\bm{\Sigma}}^{\textnormal{(even)}}_\gamma\,\mathfrak{m}
)
&\stackrel{\textnormal{p}}{\longrightarrow}&
\eta,
\nonumber\\
N p_N \, (\underline{\bm{\Sigma}}^{\textnormal{(odd)}}_\gamma\,\mathfrak{m} 
- 
\overline{\bm{\Sigma}}^{\textnormal{(even)}}_\gamma\,\mathfrak{M}
)
&\stackrel{\textnormal{p}}{\longrightarrow}&
\eta,
\eeqa
where $\eta$ corresponds to the bias for the residual estimator listed in the pertinent row of Table~(\ref{tab:Theorem1}). 
Likewise, exploiting~(\ref{eq:SigmaEvenConv}), (\ref{eq:SigmaOddConv}) and the representation of the gap $\Gamma$ of the residual estimator in Table~\ref{tab:Theorem1}, we can prove that:
\beqa
(\overline{\bm{\Sigma}}^{\textnormal{(odd)}}_\beta - \underline{\bm{\Sigma}}^{\textnormal{(even)}}_\beta)
&\stackrel{\textnormal{p}}{\longrightarrow}&
\Gamma/\kappa - 1,
\nonumber\\
(\underline{\bm{\Sigma}}^{\textnormal{(odd)}}_\beta - \overline{\bm{\Sigma}}^{\textnormal{(even)}}_\beta)
&\stackrel{\textnormal{p}}{\longrightarrow}&
\Gamma/\kappa - 1.
\eeqa
It remains to apply Lemma~\ref{lem:suffcond}, with the choices:
\beqa
\overline{\bm{w}}_N&=&\overline{\bm{\Sigma}}^{\textnormal{(odd)}}_\beta - \underline{\bm{\Sigma}}^{\textnormal{(even)}}_\beta,\nonumber\\
\underline{\bm{w}}_N&=&\underline{\bm{\Sigma}}^{\textnormal{(odd)}}_\beta - \overline{\bm{\Sigma}}^{\textnormal{(even)}}_\beta,\eeqa
and
\beqa
\overline{\bm{z}}_N&=&\overline{\bm{\Sigma}}^{\textnormal{(odd)}}_\gamma\,\mathfrak{M} 
- 
\underline{\bm{\Sigma}}^{\textnormal{(even)}}_\gamma\,\mathfrak{m},\nonumber\\
\underline{\bm{z}}_N&=&\underline{\bm{\Sigma}}^{\textnormal{(odd)}}_\gamma\,\mathfrak{m} 
- 
\overline{\bm{\Sigma}}^{\textnormal{(even)}}_\gamma\,\mathfrak{M},\nonumber\\
\eeqa
which concludes the proof of the theorem.
\end{IEEEproof}

\section{Useful Convergence Results}

\begin{lemma}[List of Convergences under Uniform Concentration]
\label{lem:listofconv}
If the connection probability fulfills~(\ref{eq:strongconc}), the convergences listed in Table~\ref{tab:Mdefs} hold true.
\end{lemma}

\begin{IEEEproof}
\begin{enumerate}
\item 
We have that:
\beq
\boxed{
\mathfrak{M}_{a,\textnormal{self}} \stackrel{\textnormal{p}}{\longrightarrow} \rho - \kappa,\qquad
\mathfrak{m}_{a,\textnormal{self}} \stackrel{\textnormal{p}}{\longrightarrow} \rho - \kappa
}
\label{eq:Maself}
\eeq
Since $\bm{a}_{ii}=\rho - \sum_{\substack{\ell\in\mathcal{N}\\ \ell\neq i}}\bm{a}_{i\ell}$, from~(\ref{eq:fundamentalassum}) we can write: 
\beq
\bm{a}_{ii} \leq \rho - \displaystyle{\frac{\kappa}{\bm{d}_{\max}}} \sum_{\substack{\ell\in\mathcal{N}\\ \ell\neq i}}\bm{g}_{i\ell}=
\rho - \kappa\,\frac{\bm{d}_i - 1}{\bm{d}_{\max}}\leq \rho - \kappa\,\frac{\bm{d}_{\min} - 1}{\bm{d}_{\max}}.
\label{eq:aiiineqappendix}
\eeq
Therefore, recalling that $\mathfrak{M}_{a,\textnormal{self}}\dfz \max_{i=1,2,\ldots,N} \bm{a}_{i i}$, we can write:
\beq
\mathfrak{M}_{a,\textnormal{self}}\leq\rho - \kappa\,\frac{\bm{d}_{\min} - 1}{\bm{d}_{\max}}.
\eeq
In view of~(\ref{eq:dmaxmin}), we have that the ratio $\frac{\bm{d}_{\min} - 1}{\bm{d}_{\max}}$ converges to $1$ in probability. 
Repeating the same reasoning with lower bounds in place of upper bounds, and with minima in place of maxima, yields the same result, and, hence, Eq.~(\ref{eq:Maself}) follows. 
\item 
We have that:
\beq
\boxed{
\mathfrak{M}_{a_2,\textnormal{self}} \stackrel{\textnormal{p}}{\longrightarrow} (\rho - \kappa)^2,\qquad
\mathfrak{m}_{a_2,\textnormal{self}} \stackrel{\textnormal{p}}{\longrightarrow} (\rho - \kappa)^2
}
\label{eq:Ma2self}
\eeq
We can write:
\beqa
\bm{a}_{ii}^{(2)}&=&\sum_{\ell\in\mathcal{N}} \bm{a}_{i\ell}\bm{a}_{\ell i}=\bm{a}_{ii}^2 + \sum_{\substack{\ell\in\mathcal{N} \\ \ell\neq i}}^N \bm{a}_{i\ell}\bm{a}_{\ell i}\nonumber\\
&\leq&
\bm{a}_{ii}^2 + \displaystyle{\frac{\kappa^2 }{\bm{d}_{\min}^2}} \sum_{\substack{\ell\in\mathcal{N} \\ \ell\neq i}}^N \bm{g}_{i\ell}
\nonumber\\
&\leq&
\bm{a}_{ii}^2 + \kappa^2 \displaystyle{\frac{\bm{d}_{\max}-1}{\bm{d}_{\min}^2}},
\eeqa
where the intermediate inequality follows by~(\ref{eq:fundamentalassum}). 
Therefore, recalling that $\mathfrak{M}_{a_2,\textnormal{self}}\dfz \max_{i\in\mathcal{N}} \bm{a}^{(2)}_{i i}$, we can write:
\beq
\mathfrak{M}_{a_2,\textnormal{self}}\leq \mathfrak{M}^2_{a,\textnormal{self}} + \kappa^2\,\frac{\bm{d}_{\max} - 1}{\bm{d}_{\min}^2},
\eeq
where the last term vanishes in probability in view of~(\ref{eq:dmaxmin}). 
Using now and~(\ref{eq:Maself}), and repeating the same reasoning with lower bounds in place of upper bounds, and with minima in place of maxima, the result in~(\ref{eq:Ma2self}) follows.
\item 
We have that:
\beq
\boxed{
\mathfrak{M}_{a,\textnormal{sum}}\stackrel{\textnormal{p}}{\longrightarrow} \rho,\qquad
\mathfrak{m}_{a,\textnormal{sum}}\stackrel{\textnormal{p}}{\longrightarrow} \rho
}
\label{eq:masum}
\eeq
The claim in~(\ref{eq:masum}) follows readily from the first relationship in~(\ref{eq:fundamentalassum}), since we can write, for any $i\neq j$: 
\beq
\sum_{\substack{\ell\in\mathcal{N} \\ \ell\neq j}}\bm{a}_{i\ell}=\rho - \bm{a}_{i j}.
\label{eq:almostrho}
\eeq
In view of~(\ref{eq:fundamentalassum}), we have that:
\beq
\bm{a}_{i j }\leq \frac{\kappa}{\bm{d}_{\min}}, 
\eeq
and, hence, $\bm{a}_{i j}$ goes to zero in probability.
\item 
We have that:
\beq
\boxed{
\mathfrak{M}^{(\mathcal{S}')}_{a,\textnormal{sum}} \stackrel{\textnormal{p}}{\longrightarrow}  \kappa(1-\xi), 
\qquad 
\mathfrak{m}^{(\mathcal{S}')}_{a,\textnormal{sum}} \stackrel{\textnormal{p}}{\longrightarrow}  \kappa(1-\xi) 
}
\label{eq:MaSprimesum}
\eeq
In view of~(\ref{eq:fundamentalassum}) we can write:
\beq
\sum_{\substack{h\in\mathcal{S}' \\  h\neq \ell,m}}\bm{a}_{h m}\leq
\frac{\kappa}{\bm{d}_{\min}}\,
\sum_{\substack{h\in\mathcal{S}' \\  h\neq \ell,m}}
\bm{g}_{hm}.
\label{eq:almostrho}
\eeq
Now we observe that the random variable:
\beq
\sum_{\substack{h\in\mathcal{S}' \\  h\neq \ell,m}}
\bm{g}_{hm}
\eeq
is a binomial random variable with number of trials equal to $S' - 2$, and success probability equal to $p_N$. 
In other words, we get the following representation:
\beq
\max_{\substack{\ell,m\in\mathcal{S}' \\ \ell\neq m}}\sum_{\substack{h\in\mathcal{S}' \\  h\neq \ell,m}}
\bm{g}_{hm}
=\bm{\mathcal{B}}_{\max}(S'-2,p_N,(S'-1)S'),
\label{eq:ourvariableasBmax}
\eeq
because maximization is carried over all pairs $\ell,m\in\mathcal{S}'$, with $\ell\neq m$.
Moreover, since in the uniform concentration regime we have:
\beq
p_N=\omega_N\frac{\log N}{N},\qquad \omega_N\stackrel{N\rightarrow\infty}{\longrightarrow} \infty,
\eeq
and since $S'/N\rightarrow 1-\xi$ as $N\rightarrow\infty$, we can regard the connection probability $p_N$ as a connection probability scaling with respect to $S'$, namely,
\beqa
p_N&=&\omega_N \frac{\log S'}{S'}\,\frac{S'}{N}\,\frac{\log N}{\log(S'/N)+\log N}\nonumber\\
&=&\omega_{S'} \frac{\log S'}{S'}\dfz p_{S'},
\label{eq:pconnSprime}
\eeqa
where:
\beq
\omega_{S'}=\omega_N\frac{S'}{N}\frac{\log N}{\log(S'/N)+\log N}\stackrel{N\rightarrow\infty}{\longrightarrow}\infty.
\eeq
This shows that the uniform concentration regime can be referred also to the scaling of the involved quantities w.r.t. $S'$ (in place of $N$). 
Accordingly, we can apply~(\ref{eq:basicmax}) with $K=(S'-1)S'$, $N=S'-2$, and $p_{S'}=\omega_{S'} \log(S')/S'$, which yields:
\beq
\displaystyle{
\frac{\bm{\mathcal{B}}_{\max}(S'-2,p_N,(S'-1)S')}{S' p_N}
}
\stackrel{\textnormal{p}}{\longrightarrow} 1.
\label{eq:dmaxoverSprime}
\eeq
It remains to rewrite~(\ref{eq:almostrho}) as:
\beqa
\sum_{\substack{h \in\mathcal{S}' \\ h\neq \ell,m}}\bm{a}_{h m}
&\leq&
\kappa\,
\underbrace{\frac{N p_N}{\bm{d}_{\min}}}_{\stackrel{\textnormal{p}}{\longrightarrow} 1}
\underbrace{\frac{S'}{N}}_{\rightarrow 1-\xi}
\nonumber\\
&\times&
\underbrace{
\frac{\bm{\mathcal{B}}_{\max}(S'-2,p_N,(S'-1)S')}{S' p_N}}_{\stackrel{\textnormal{p}}{\longrightarrow} 1}.\nonumber\\
\eeqa
Repeating the same reasoning with lower bounds in place of upper bounds, and with minima in place of maxima, yields the same result, and, hence, Eqs.~(\ref{eq:MaSprimesum}) follow.
\item 
We have that:
\beq
\boxed{
\mathfrak{M}_{c,\textnormal{self}} \stackrel{\textnormal{p}}{\longrightarrow} (\rho - \kappa)^2,
\qquad
\mathfrak{m}_{c,\textnormal{self}} \stackrel{\textnormal{p}}{\longrightarrow} (\rho - \kappa)^2
}
\label{eq:Mcself}
\eeq
This result follows readily by repeating the same steps used to prove~(\ref{eq:Ma2self}).
\item
We have that:
\beq
\boxed{
\begin{array}{lll}
\mathfrak{M}_{c,\textnormal{sum}} \stackrel{\textnormal{p}}{\longrightarrow}& \rho^2 - 2 \rho\kappa\xi + \kappa^2 \xi,
\\
\\
\mathfrak{m}_{c,\textnormal{sum}} \stackrel{\textnormal{p}}{\longrightarrow}& \rho^2 - 2 \rho\kappa\xi + \kappa^2 \xi
\end{array}
}
\label{eq:Mcsum}
\eeq
Using the definition of $\bm{C}$ in~(\ref{eq:Hmatdef}), we note that we can write:
\beqa
\sum_{\substack{h\in\mathcal{S}' \\ h\neq m}} \bm{c}_{\ell h}&=&
\sum_{\substack{h\in\mathcal{S}' \\ h\neq m}} \sum_{j\in\mathcal{N}} \bm{a}_{\ell j} \bm{a}_{j h}
\nonumber\\
&=&
\sum_{j\in\mathcal{S}} \bm{a}_{\ell j} \sum_{\substack{h\in\mathcal{S}' \\ h\neq m}} \bm{a}_{j h}+
\sum_{j\in\mathcal{S}'} \bm{a}_{\ell j} \sum_{\substack{h\in\mathcal{S}' \\ h\neq m}}  \bm{a}_{j h}.\nonumber\\
\label{eq:Mcsumdecomp}
\eeqa
Applying the same procedure used in the previous items of this section, it is readily proved that, if $j\in\mathcal{S}$ (and, hence the self-term $\bm{a}_{\ell \ell}$ is not present, because $\ell\in\mathcal{S}'$):
\beq
\max_{j\in\mathcal{S}} \sum_{\substack{h\in\mathcal{S}' \\ h\neq m}} \bm{a}_{j h}
\stackrel{\textnormal{p}}{\longrightarrow} \kappa (1-\xi),~
\min_{j\in\mathcal{S}} \sum_{\substack{h\in\mathcal{S}' \\ h\neq m}} \bm{a}_{j h}
\stackrel{\textnormal{p}}{\longrightarrow} \kappa (1-\xi),
\label{eq:Csum1}
\eeq
whereas, if $j\in\mathcal{S}'$:
\beq
\max_{j\in\mathcal{S}} \sum_{\substack{h\in\mathcal{S}' \\ h\neq m}} \bm{a}_{j h}
\stackrel{\textnormal{p}}{\longrightarrow} \rho - \kappa \xi,~
\min_{j\in\mathcal{S}} \sum_{\substack{h\in\mathcal{S}' \\ h\neq m}} \bm{a}_{j h}
\stackrel{\textnormal{p}}{\longrightarrow} \rho - \kappa \xi.
\label{eq:Csum2}
\eeq
Likewise, we can show that:
\beq
\max_{\ell\in\mathcal{S}'} \sum_{j\in\mathcal{S}} \bm{a}_{\ell j}
\stackrel{\textnormal{p}}{\longrightarrow} \kappa \xi,~
\min_{\ell\in\mathcal{S}'} \sum_{j\in\mathcal{S}} \bm{a}_{\ell j}
\stackrel{\textnormal{p}}{\longrightarrow} \kappa \xi,
\label{eq:Csum3}
\eeq
and that:
\beq
\max_{\ell\in\mathcal{S}'} \sum_{j\in\mathcal{S}'} \bm{a}_{\ell j}
\stackrel{\textnormal{p}}{\longrightarrow} \rho - \kappa \xi,~
\min_{\ell\in\mathcal{S}'} \sum_{j\in\mathcal{S}'} \bm{a}_{\ell j}
\stackrel{\textnormal{p}}{\longrightarrow} \rho - \kappa \xi.
\label{eq:Csum4}
\eeq
Plugging~(\ref{eq:Csum1})--~(\ref{eq:Csum4}) into~(\ref{eq:Mcsumdecomp}) finally yields~(\ref{eq:Mcsum}). 
\item 
We have that:
\beq
\boxed{
N p_N \,\mathfrak{M}\stackrel{\textnormal{p}}{\longrightarrow} \kappa^2 p,
\qquad
N p_N \,\mathfrak{m}\stackrel{\textnormal{p}}{\longrightarrow} \kappa^2 p 
}
\label{eq:Msimply}
\eeq
In view of~(\ref{eq:fundamentalassum}), we can write:
\beq
\sum_{\substack{\ell\in\mathcal{N} \\ \ell \neq i,j}}\bm{a}_{i\ell}\bm{a}_{\ell j}
\leq
\frac{\kappa^2}{{\bm{d}^2_{\min}}}
\sum_{\substack{\ell\in\mathcal{N} \\ \ell \neq i,j}}\bm{g}_{i\ell}\bm{g}_{\ell j}.
\label{eq:Mintermediate}
\eeq
Now we see that the quantity:
\beq
\sum_{\substack{\ell\in\mathcal{N} \\ \ell \neq i,j}}\bm{g}_{i\ell}\bm{g}_{\ell j}
\eeq
is a binomial random variable with number of trials equal to $N-2$, and success probability equal to $p_N^2$, since when $\ell\neq i,j$, the product variable $\bm{g}_{i\ell}\bm{g}_{\ell j}$ is a Bernoulli variable with success probability $p_N^2$. Therefore, we are allowed to introduce the definition:
\beq
\bm{\mathcal{B}}_{\max}(N-2,p^2_N,(N-1)N)=
\max_{\substack{i,j\in\mathcal{N} \\ i\neq j}}
\sum_{\substack{\ell\in\mathcal{N} \\ \ell \neq i,j}}\bm{g}_{i\ell}\bm{g}_{\ell j},
\eeq
which, in view of Lemma~\ref{lem:fundmaxN2lemma}, yields: 
\beq
\frac{1}{{N p_N}}\,\max_{\substack{i,j\in\mathcal{N} \\ i\neq j}}\sum_{\ell\neq i,j}\bm{g}_{i\ell}\bm{g}_{\ell j}\stackrel{\textnormal{p}}{\longrightarrow} p.
\label{eq:Max2Intermediate}
\eeq
Now, from the definition of $\mathfrak{M}$ in Table~\ref{tab:Mdefs}, line $12)$, we get: 
\beq
N p_N \mathfrak{M}
\leq
\kappa^2\,
\underbrace{
\displaystyle{\frac{\bm{\mathcal{B}}_{\max}(N-2,p^2_N,(N-1)N)}{N p_N}}
}_{\stackrel{\textnormal{p}}{\longrightarrow} p}
\,\underbrace{\frac{N^2 p_N^2}{\bm{d}^2_{\min}}}_{\stackrel{\textnormal{p}}{\longrightarrow} 1}.
\label{eq:Mintermediate}
\eeq
Repeating the same reasoning with lower bounds in place of upper bounds, and with minima in place of maxima, we get the claim in~(\ref{eq:Msimply}).
\item
We have that: 
\beq
\boxed{
\begin{array}{ll}
N p_N \, \mathfrak{M}^{(\mathcal{S}')}&\!\!\!\!\stackrel{\textnormal{p}}{\longrightarrow} \kappa^2 p (1-\xi),
\\
N p_N \, \mathfrak{m}^{(\mathcal{S}')}&\!\!\!\!\stackrel{\textnormal{p}}{\longrightarrow} \kappa^2 p (1-\xi)
\end{array}
}
\label{eq:MSprime}
\eeq
The proof for the case where $p=0$ comes from~(\ref{eq:Msimply}) because, from definitions $12)$ and $13)$ in Table~\ref{tab:Mdefs}, we see that $\mathfrak{M}^{(\mathcal{S}')}\leq \mathfrak{M}$. 
The proof for the case where $p>0$ is readily obtained by using the same arguments leading to~(\ref{eq:Msimply}).
\item 
We have that:
\beq
\boxed{
\begin{array}{ll}
N p_N\,\mathfrak{M}_{a_3,\textnormal{sum}}^{(\mathcal{S}')}&\!\!\!\stackrel{\textnormal{p}}{\longrightarrow} \kappa^3 p (1-\xi)^2,
\\
N p_N\,\mathfrak{m}_{a_3,\textnormal{sum}}^{(\mathcal{S}')}&\!\!\!\stackrel{\textnormal{p}}{\longrightarrow} \kappa^3 p (1-\xi)^2
\end{array}
}
\label{eq:Ma3sumSprime}
\eeq
We have that:
\beqa
\lefteqn{\sum_{\substack{\ell, m\in\mathcal{S}' \\ \ell\neq m}}\bm{a}_{i\ell}\bm{a}_{\ell m}\bm{a}_{m j}}\nonumber\\
&=&
\sum_{\ell\in\mathcal{S}'}\bm{a}_{i\ell}
\sum_{\substack{m\in\mathcal{S}' \\ m\neq \ell}}\bm{a}_{\ell m}\bm{a}_{m j}\nonumber\\
&\leq&
\max_{i\in\mathcal{S}}\sum_{\ell\in\mathcal{S}'}\bm{a}_{i\ell}
\max_{j\in\mathcal{S},\ell\in\mathcal{S}'}
\sum_{\substack{m\in\mathcal{S}' \\ m\neq \ell}}\bm{a}_{\ell m}\bm{a}_{m j}.
\label{eq:a3sumSprimeintermediate}
\eeqa
Reasoning as done for proving~(\ref{eq:MSprime}), we can show that:
\beq
N p_N\max_{j\in\mathcal{S},\ell\in\mathcal{S}'}\sum_{\substack{m\in\mathcal{S}' \\ m\neq \ell}}\bm{a}_{\ell m}\bm{a}_{m j}
\stackrel{\textnormal{p}}{\longrightarrow}\kappa^2 p (1-\xi).
\label{eq:maxa2sumSprime}
\eeq
Likewise, reasoning as done for proving~(\ref{eq:MaSprimesum}), we can show that:
\beq
\max_{i\in\mathcal{S}}\sum_{\ell\in\mathcal{S}'}\bm{a}_{i\ell}
\stackrel{\textnormal{p}}{\longrightarrow} \kappa(1-\xi).
\label{eq:maxasumSprime}
\eeq
Finally, using~(\ref{eq:maxa2sumSprime}) and~(\ref{eq:maxasumSprime}) into~(\ref{eq:a3sumSprimeintermediate}), repeating the same reasoning with lower bounds in place of upper bounds, and with minima in place of maxima, we get~(\ref{eq:Ma3sumSprime}).
\item
The following list of convergences is obtained by trivial variations on the previous proofs.
\beq
\boxed{
\widetilde{\mathfrak{M}}^{(\mathcal{S}')}
\stackrel{\textnormal{p}}{\longrightarrow}
\kappa^2 (1-\xi)^2,~~
\widetilde{\mathfrak{m}}^{(\mathcal{S}')}
\stackrel{\textnormal{p}}{\longrightarrow}
\kappa^2 (1-\xi)^2
}
\label{eq:MSprime2}
\eeq

\beq
\boxed{
\widetilde{\widetilde{\mathfrak{M}}}^{(\mathcal{S}')}
\stackrel{\textnormal{p}}{\longrightarrow}
\kappa^2 (1-\xi)^2,~~
\widetilde{\widetilde{\mathfrak{m}}}^{(\mathcal{S}')}
\stackrel{\textnormal{p}}{\longrightarrow}
\kappa^2 (1-\xi)^2
}
\label{eq:MSprime3}
\eeq

\beq
\boxed{
\widetilde{\mathfrak{M}}_{a,\textnormal{sum}}^{(\mathcal{S}')}
\stackrel{\textnormal{p}}{\longrightarrow}
\kappa (1-\xi),~~
\widetilde{\mathfrak{m}}_{a,\textnormal{sum}}^{(\mathcal{S}')}
\stackrel{\textnormal{p}}{\longrightarrow}
\kappa (1-\xi)
}
\label{eq:MSprimeasum}
\eeq

\beq
\boxed{
\widetilde{\widetilde{\mathfrak{M}}}_{a,\textnormal{sum}}^{(\mathcal{S}')}
\stackrel{\textnormal{p}}{\longrightarrow}
\kappa^3 (1-\xi)^2,~~
\widetilde{\widetilde{\mathfrak{m}}}_{a,\textnormal{sum}}^{(\mathcal{S}')}
\stackrel{\textnormal{p}}{\longrightarrow}
\kappa^3 (1-\xi)^2
}
\label{eq:MSprimeasum2}
\eeq

\end{enumerate}
\end{IEEEproof}

\section{Useful Lemma}
\begin{lemma}
\label{lem:lemmaseries}
Let $0<\alpha<1$, $0<\rho_{\ell}<1$ and $\beta_{\ell}\in\mathbb{R}$ for all $\ell=1,2,\ldots,L$, and introduce the following recursion: 
\beq
f_{k+1}=\alpha f_k +\sum_{\ell=1}^L \beta_{\ell} \rho_{\ell}^k.
\label{eq:frec0}
\eeq
Then, $f_k$ is equal to:
\beq
\left(f_0 +\sum_{\ell=1}^L \frac{\beta_{\ell}}{\alpha-\rho_{\ell}}\right)\alpha^{k} - \sum_{\ell=1}^L \frac{\beta_{\ell}}{\alpha-\rho_{\ell}} \rho_{\ell}^k,
\label{eq:frec}
\eeq
and, hence, can be cast in the following form:
\beq
f_k=\sum_{\ell=0}^{L} \widetilde{\beta}_{\ell} \widetilde{\rho}_{\ell}^k,
\label{eq:fsol}
\eeq
with obvious choices for $\widetilde{\beta}_{\ell}$ and $\widetilde{\rho}_{\ell}$.
\end{lemma}
\begin{IEEEproof}
Exploiting~(\ref{eq:frec0}), we can write:
\beqa
f_1&=&\alpha f_0 + \sum_{\ell=1}^L \beta_{\ell}, \nonumber\\
f_2&=&\alpha^2 f_0 + \alpha\sum_{\ell=1}^L \beta_{\ell} +\sum_{\ell=1}^L \beta_{\ell}\rho_{\ell},\nonumber\\
&\vdots&\nonumber\\
f_k&=&\alpha^k f_0 + \sum_{\ell=1}^L \beta_{\ell}\sum_{j=0}^{k-1}\alpha^{k-1-j}\rho_{\ell}^j.
\eeqa
The last equation can be manipulated as follows:
\beqa
f_k&=&\alpha^k f_0 + \sum_{\ell=1}^L \beta_{\ell}\alpha^{k-1}\sum_{j=0}^{k-1}\left(\frac{\rho_{\ell}}{\alpha}\right)^j\nonumber\\
&=&\alpha^k f_0 + \sum_{\ell=1}^L \beta_{\ell}\alpha^{k-1}\frac{1 - (\rho_{\ell}/\alpha)^k}{1-\rho_{\ell}/\alpha}
\nonumber\\
&=&\alpha^k f_0 + \sum_{\ell=1}^L \beta_{\ell}\frac{\alpha^k - \rho_{\ell}^k}{\alpha-\rho_{\ell}},
\eeqa
which corresponds to~(\ref{eq:frec}). 
\end{IEEEproof}

\section{Sample Consistency}
\label{app:samplecons}
We start with a useful lemma that characterizes the output of the clustering algorithm proposed in Sec.~\ref{subsec:clualgo}, when its input is constituted by two well-separated classes.
\begin{lemma}[Consistency of the clustering algorithm]
\label{lem:myclu}
Let 
\beq
v_1\leq v_2\leq\ldots\leq v_L
\eeq
be the set of points to be clustered in two classes. 
Let $\nu_0<\nu_1$ be two real numbers such that, for some $k\in\{1,2,\ldots,L\}$ and for all $i\in\{1,2,\ldots,L\}$:
\beq
|v_i-\nu_0|<\epsilon~~\textnormal{if}~~i\leq k,~~~~~~
|v_i-\nu_1|<\epsilon~~\textnormal{if}~~i>k.
\eeq
Then, for sufficiently small $\epsilon$, the optimal partition of the clustering algorithm is that corresponding to $j^{\star}=k$, namely:
\beqa
\mathcal{C}_0(k)&=&\{v_i: |v_i-\nu_0| <\epsilon\},\nonumber\\
\mathcal{C}_1(k)&=&\{v_i: |v_i-\nu_1|<\epsilon\}.
\label{eq:C01jstar}
\eeqa
\end{lemma}
\begin{IEEEproof}
It suffices to prove the lemma for the choice $\nu_0=0$. 
We first show that~(\ref{eq:C01jstar}) is an admissible configuration. To this aim, we need to show that the midpoint between the centroids separates the clusters. The centroids of the clusters in~(\ref{eq:C01jstar}) fulfill in fact the following bounds:
\beq
-\epsilon<c_0(k)<\epsilon,\qquad \nu_1-\epsilon<c_1(k)<\nu_1+\epsilon,
\label{eq:centstarbounds}
\eeq
and, hence, their midpoint fulfills the condition:
\beq
\frac{\nu_1}{2}-\epsilon<\frac{c_0(k)+c_1(k)}{2}<\frac{\nu_1}{2}+\epsilon.
\eeq
For sufficiently small $\epsilon$, this configuration is clearly admissible, because the midpoint between the two centroids separates the two clusters. Moreover, in view of~(\ref{eq:centstarbounds}) the distance between the centroids fulfills the bound:
\beq
c_1(k)-c_0(k)>\nu_1-2\epsilon.
\label{eq:maxdistcentbound}
\eeq
Next we show that another admissible configuration (if any) distinct from~(\ref{eq:C01jstar}) exhibits a smaller inter-cluster distance. 
Let us consider a configuration $\mathcal{C}_0(j)$ with $j>k$. 
Since now the point $v_j>\nu_1-\epsilon$ belongs to $\mathcal{C}_0(j)$, the centroid midpoint that separates the two clusters must be not smaller than $v_j$, yielding:
\beq
\nu_1-\epsilon<\frac{c_0(j)+c_1(j)}{2}\Rightarrow 
-c_0(j)<c_1(j) - 2(\nu_1-\epsilon).
\label{eq:cent0j}
\eeq
On the other hand, we have clearly:
\beq
c_1(j)<\nu_1+\epsilon,
\label{eq:cent1j}
\eeq
which used along with~(\ref{eq:cent0j}) yields:
\beq
c_1(j)-c_0(j)<2 c_1(j) - 2(\nu_1-\epsilon)<4\epsilon.
\eeq
Thus, for sufficiently small $\epsilon$, we see that any $j>k$ would produce a distance between the centroids that is dominated by the distance in~(\ref{eq:maxdistcentbound}), and since the situation is similar for $j<k$, the proof is complete.
\end{IEEEproof}
\begin{IEEEproof}[Proof of Theorem~\ref{theor:samplecons}] 
In the following, for $i,j\in\mathcal{S}$, the $(i,j)$-th entry of the estimated matrix $\widehat{\bm{A}}_{\mathcal{S},n}$ is denoted by $\widehat{\bm{a}}_{ij}(n)$. Likewise, the $(i,j)$-th entry of the limiting estimated matrix $\widehat{\bm{A}_{\mathcal{S}}}$ is denoted by $\widehat{\bm{a}}_{ij}$. 
Let us introduce the following events (we recall that bold notation refers to random variables):
\beqa
\mathcal{E}_0(n,N)&=&\left\{
\max_{\substack{i,j\in\mathcal{S}: \bm{a}_{ij}=0 \\ i\neq j}}
\left|
s_N\widehat{\bm{a}}_{ij}(n)-\eta
\right|<\epsilon\right\},\nonumber\\
\mathcal{E}_1(n,N)&=&\left\{\max_{\substack{i,j\in\mathcal{S}: \bm{a}_{ij}>0 \\ i\neq j}}
\left|
s_N\widehat{\bm{a}}_{ij}(n)-\eta-\Gamma
\right|<\epsilon
\right\}.\nonumber\\
\eeqa
In view of Lemma~\ref{lem:myclu}, assuming a sufficiently small $\epsilon$, the clustering algorithm proposed in Sec.~\ref{subsec:clualgo} reconstructs the exact graph whenever $\mathcal{E}_0(n,N)\cap\mathcal{E}_1(n,N)$ occurs. 
Accordingly, the theorem will be proved if we show that for some $n(N)$:
\beq
\lim_{N\rightarrow\infty}\P[\mathcal{E}_0(n(N),N)\cap\mathcal{E}_1(n(N),N)]=1.
\label{eq:claiminevents}
\eeq
To this aim, let us introduce the events:
\beqa
\mathcal{E}_0(N)&=&\left\{
\max_{\substack{i,j\in\mathcal{S}: \bm{a}_{ij}=0 \\ i\neq j}}
\left|
s_N\widehat{\bm{a}}_{ij}-\eta
\right|<\epsilon/2\right\},\nonumber\\
\mathcal{E}_1(N)&=&
\left\{\max_{\substack{i,j\in\mathcal{S}: \bm{a}_{ij}>0 \\ i\neq j}}
\left|
s_N\widehat{\bm{a}}_{ij}-\eta-\Gamma
\right|<\epsilon/2
\right\},\nonumber\\
\eeqa
and
\beq
\mathcal{F}(n,N)=\left\{
s_N\,\|\widehat{\bm{A}}_{\mathcal{S},n}-\widehat{\bm{A}}_{\mathcal{S}}\|_{\max}<\epsilon/2\right\}.
\label{eq:fnN}
\eeq
By triangle inequality, we have the implication:
\beq
\mathcal{E}_0(N)\cap\mathcal{E}_1(N)\cap\mathcal{F}(n,N)
\Rightarrow
\mathcal{E}_0(n,N)\cap\mathcal{E}_1(n,N),
\eeq
yielding:
\beq
\P\left[\mathcal{E}_0(n,N)\cap\mathcal{E}_1(n,N)\right]\geq
\P\left[\mathcal{E}_0(N)\cap\mathcal{E}_1(N)\cap\mathcal{F}(n,N)\right].
\label{eq:eventineq01}
\eeq
Since by assumption the limiting matrix estimator $\widehat{\bm{A}}_{\mathcal{S}}$ achieves universal local structural consistency we know that:
\beq
\lim_{N\rightarrow\infty} \P\left[
\mathcal{E}_0(N)\cap\mathcal{E}_1(N)
\right]=1.
\label{eq:ulsevent}
\eeq
Consider now a sequence $\epsilon_N>0$ that vanishes as $N\rightarrow\infty$. 
Since $\widehat{\bm{A}}_{\mathcal{S},n}$ is in the class defined by~(\ref{eq:asystableclass}), for any fixed $N$, there exists $n(N)$ such that, for all $n\geq n(N)$ we have: 
\beq
\P\left[
\mathcal{F}(n,N)
\right]>1-\epsilon_N,
\label{eq:asystableevent}
\eeq
or, by the sandwich theorem:
\beq
\lim_{N\rightarrow\infty} \P\left[
\mathcal{F}(n(N),N)
\right]=1.
\label{eq:asystableevent}
\eeq
Using now~(\ref{eq:ulsevent}) and~(\ref{eq:asystableevent}) in~(\ref{eq:eventineq01}) implies~(\ref{eq:claiminevents}) and, hence, the claim of the theorem.
\end{IEEEproof}

\section{Sample Complexity}
\label{app:samplecomple}
In the following we will make repeated use of the following inequality, holding for any two matrices $A$, $B$:
\beq
\|A B\|_{\max}
\leq
\min(\|A\|_{\max}\|B\|_1,\|A\|_{\infty}\|B\|_{\max}).
\eeq

Given a regular diffusion matrix fulfilling Assumption~\ref{assum:regdiffmat} with parameters $\rho$ and $\kappa$, let $\zeta=\rho-\kappa$. We introduce, for a small $\delta$, with $0<\delta<1-\zeta^2$, the auxiliary quantities:
\beq
\varphi_1\dfz\frac{(1-\zeta^2-\delta)(1-\rho)}{16\sqrt{2}\sigma^2},~~
\varphi_2\dfz\frac{(1-\rho^2)^2 (1-\rho)}{16\sqrt{2}\sigma^2}.\label{eq:varphidef}
\eeq
We notice that when $\delta$ is sufficiently small, $\varphi_1>\varphi_2$.
Next we introduce some auxiliary functions. 
Let
\beq
{\sf b}_{n}(x)\!=\!\!\left\{\!
\begin{array}{ll}
\!\!\!\min\left\{
1\,,\, S^2 \left(e^{-n/2}
+
e^{-\left[\sqrt{n} x -\sqrt{2}\right]^2}\right)
\!\right\},\!\!\!\!&\sqrt{n}x>\sqrt{2}
\\
\\
\!\!\!1,&\sqrt{n}x\leq \sqrt{2}
\end{array}
\right.
\label{eq:anphidef}
\eeq
where $x$ is a positive quantity, $n$ is the sample size and $S$ is the number of probed nodes. 
It is immediate to verify that ${\sf b}_{n}(x)$ is non-increasing in $x$, and that:
\beq
\lim_{n\rightarrow\infty} {\sf b}_n(x)=0~~\textnormal{for all }x>0.
\label{eq:anphiconv}
\eeq
The second auxiliary function is:
\beq
{\sf b}_N(\delta)=\P\left[
\frac{\min_{i\in\mathcal{N}} [\bm{R}_0]_{ii}}
{(\max_{i\in\mathcal{N}} [\bm{R}_0]_{ii})^2}
<\frac{1-\zeta^2-\delta}{\sigma^2}
\right].
\label{eq:bndeldefapp}
\eeq
We have that:
\beq
\lim_{N\rightarrow\infty} {\sf b}_N(\delta)=0.
\label{eq:minmaxR0conv}
\eeq
The convergence in~(\ref{eq:minmaxR0conv}) follows by observing that:
\beq
[\bm{R}_0]_{ii}=\sigma^2[(I-\bm{A}^2)^{-1}]_{ii}=
\sigma^2\left(1+\sum_{h=1}^{\infty}\bm{a}^{(2h)}_{ii}\right),
\label{eq:R0selfappexpan}
\eeq
which, in light of~(\ref{eq:SigmaEvenConv}), implies~(\ref{eq:minmaxR0conv}).

The following lemma characterizes the rate of convergence of the sample correlation estimators. 
This lemma adapts Lemmas $1$ and $2$ in~\cite{HanLuLiuJMLR} to exploit the peculiarities of our setting.
\begin{lemma}[Convergence rate of the sample correlation estimators]
\label{lem:concorr}
Let us consider the VAR model in~(\ref{eq:VARmodel}) with Gaussian input source and initial state distributed according to the stationary distribution. Let the combination matrix $\bm{A}$ fulfill Assumption~\ref{assum:regdiffmat} with parameters $\rho$ and $\kappa$, and let the underlying graph fulfill the Erd\H{o}s-R\'enyi model under the uniform concentration regime. 
Then, for all $\epsilon>0$ we have:
\beqa
\lefteqn{
\P\left[\|[\widehat{\bm{R}}_{0,n}]_{\mathcal{S}}-[\bm{R}_0]_{\mathcal{S}}\|_{\max}>\epsilon\right]}
\nonumber\\
&\leq&
3\Big[
{\sf b}_{n}(\epsilon\,\varphi_1)
+
{\sf b}_{n}(\epsilon\,\varphi_2)
\,{\sf b}_N(\delta)\Big],
\label{eq:boundR0app}
\eeqa
and
\beqa
\lefteqn{\P\left[\|[\widehat{\bm{R}}_{1,n}]_{\mathcal{S}}-[\bm{R}_1]_{\mathcal{S}}\|_{\max}>\epsilon\right]}
\nonumber\\
&\leq&
4\Big[
{\sf b}_{n-1}(\epsilon\,\varphi_1)
+
{\sf b}_{n-1}(\epsilon\,\varphi_2)
\,{\sf b}_N(\delta)\Big].
\label{eq:boundR1app}
\eeqa
\end{lemma}
\begin{IEEEproof}
Let us start by examining the probability in~(\ref{eq:boundR0app}) for a fixed realization of the combination matrix $\bm{A}=A$, i.e., we consider the quantity, for $i,j\in\mathcal{N}$:
\beq
\P[|\widehat{\bm{R}}_{0,n}-R_0|_{ij}|>\epsilon],~~~R_0=\sigma^2(I-A^2)^{-1},
\label{eq:corrboundeterministic}
\eeq 
where the conditioning on the event $\bm{A}=A$ is removed since the input source $\{\bm{x}_n\}$ in~(\ref{eq:VARmodel}) is independent from the random generation of the combination matrix. 
An upper bound on the probability appearing in~(\ref{eq:corrboundeterministic}) is derived in~\cite{HanLuLiuJMLR}. 
We now modify the proof of Lemma~$1$ in~\cite{HanLuLiuJMLR} to exploit the peculiarities of our model. 
To this end, let us consider, for $m,m'=1,2,\ldots,n$, the following known relationship that is obtained exploiting~(\ref{eq:VARmodel}):\footnote{The specific representation $R_{m-m'}=A^{|m'-m|} R_0$ holds due to the symmetry of $A$.}
\beq
R_{m-m'}=\E[\bm{y}_m\bm{y}_{m'}^{\top}]=A^{|m-m'|} R_0,
\eeq
In~\cite{HanLuLiuJMLR} the following inequality is used:
\beq
\|R_{m-m'}\|_{\max}\leq \|R_{m-m'}\|_2\leq \rho^{|m-m'|} \|R_0\|_2.
\label{eq:lemma1otherpaper}
\eeq 
In our case we can replace~(\ref{eq:lemma1otherpaper}) by the following inequality, which exploits additional constraints on $A$:
\beqa
\|R_{m-m'}\|_{\max}&\leq& \|A^{|m-m'|}\|_{\infty} \|R_0\|_{\max}\nonumber\\
&=&\rho^{|m-m'|} \max_{i\in\mathcal{N}} [R_0]_{ii},
\label{eq:betterboundapp}
\eeqa
where the last equality comes from noticing that: $i)$ the matrix $A/\rho$ is doubly stochastic; and $ii)$ the off-diagonal entries of the correlation matrix are upper bounded as $|[R_0]_{ij}|\leq\sqrt{[R_0]_{ii}[R_0]_{jj}}$ by Cauchy-Schwarz inequality. 
Let us now introduce the following definitions:
\beq
\varphi(R_0)=\frac{(1-\rho)}{16\sqrt{2}}\frac{\min_{i\in\mathcal{N}} [R_0]_{ii}}{(\max_{i\in\mathcal{N}} [R_0]_{ii})^2}.
\label{eq:variousdefincludingphi}
\eeq
Applying~(\ref{eq:betterboundapp}) in the proof of Lemma~$1$ in~\cite{HanLuLiuJMLR} (with every other detail of the proof being unaltered) we get, for $\sqrt{n}\epsilon\varphi(R_0)>\sqrt{2}$ and all $i,j\in\mathcal{N}$ with $i\neq j$: 
\beqa
\lefteqn{\P\left[|[\widehat{\bm{R}}_{0,n}]_{ij}-[R_0]_{ij}|>\epsilon\right]}
\nonumber\\
&\leq&
3 \left(e^{-n/2}+
e^{-\left[\sqrt{n} \epsilon\,\varphi(R_0) -\sqrt{2}\right]^2}\right),
\label{eq:boundR0app2}
\eeqa
which, using the union bound over the set of probed nodes $\mathcal{S}$ and the definition in~(\ref{eq:anphidef}) yields, for all $n$: 
\beq
\P\left[\|[\widehat{\bm{R}}_{0,n}]_{\mathcal{S}}-[R_0]_{\mathcal{S}}\|_{\max}>\epsilon\right]\leq
3 \,{\sf b}_{n}(\epsilon\,\varphi(R_0)).
\label{eq:boundR0app2}
\eeq
We now complete the proof by taking into account the randomness of the combination matrix $\bm{A}$. 
To this end, let $\mathcal{R}$ the set of all possible realizations of $\bm{R}_0$ (we recall that $\bm{R}_0$ depends on $\bm{A}$, which can take only a finite number of values since the possible realizations of the graph are finite). 
Let us introduce the following set:
\beq
\mathcal{T}\dfz\left\{R_0\in\mathcal{R}: \frac{\min_{i\in\mathcal{N}} [R_0]_{ii}}{(\max_{i\in\mathcal{N}} [R_0]_{ii})^2}\geq
\frac{1-\zeta^2-\delta}{\sigma^2}\right\}.
\eeq
By the law of total probability we can write:
\beqa
\lefteqn{\P\Big[\|[\widehat{\bm{R}}_{0,n}]_{\mathcal{S}}-[\bm{R}_0]_{\mathcal{S}}\|_{\max}>\epsilon\Big]}
\nonumber\\
&=&
\sum_{R_0\in\mathcal{R}} 
\P\Big[\|[\widehat{\bm{R}}_{0,n}]_{\mathcal{S}}-[R_0]_{\mathcal{S}}\|_{\max}>\epsilon\Big]\;
\P\Big[\bm{R}_0=R_0\Big]
\nonumber\\
&=&
\sum_{R_0\in\mathcal{T}} 
\P\Big[\|[\widehat{\bm{R}}_{0,n}]_{\mathcal{S}}-[R_0]_{\mathcal{S}}\|_{\max}>\epsilon\Big]\;
\P\Big[\bm{R}_0=R_0\Big]
\nonumber\\
&+&
\sum_{R_0\in\mathcal{R}\setminus\mathcal{T}} 
\P\Big[\|[\widehat{\bm{R}}_{0,n}]_{\mathcal{S}}-[R_0]_{\mathcal{S}}\|_{\max}>\epsilon\Big]\;
\P\Big[\bm{R}_0=R_0\Big]
\nonumber\\
&\leq&
3\sum_{R_0\in\mathcal{T}} 
{\sf b}_{n}(\epsilon\,\varphi(R_0))
\;
\P\Big[\bm{R}_0=R_0\Big]
\label{eq:totalprobapp31}
\\
&+&
3\sum_{R_0\in\mathcal{R}\setminus\mathcal{T}} 
{\sf b}_{n}(\epsilon\,\varphi(R_0))
\;
\P\Big[\bm{R}_0=R_0\Big],
\label{eq:totalprobapp32}
\eeqa
where we have used~(\ref{eq:boundR0app2}).
Recalling now the definition of $\varphi_1$ in~(\ref{eq:varphidef}) and of $\varphi(R_0)$ in~(\ref{eq:variousdefincludingphi}), we see that when $R_0\in\mathcal{T}$ we have $\varphi(R_0)>\varphi_1$, yielding, since ${\sf b}_n(x)$ is non-increasing in $x$:
\beq
{\sf b}_{n}(\epsilon\,\varphi(R_0))\leq {\sf b}_{n}(\epsilon\, \varphi_1).
\label{eq:boundphi31}
\eeq
Moreover, from~(\ref{eq:R0selfappexpan}) we see that:
\beq
\sigma^2 \leq [R_0]_{ii} \leq \frac{\sigma^2}{1-\rho^2},
\eeq
yielding:
\beq
\frac{\min_{i\in\mathcal{N}} [\bm{R}_0]_{ii}}
{(\max_{i\in\mathcal{N}} [\bm{R}_0]_{ii})^2}\geq \frac{(1-\rho^2)^2}{\sigma^2},
\eeq
and in view of the definition of $\varphi_2$ in~(\ref{eq:varphidef}) we conclude that $\varphi(R_0)\geq\varphi_2$ for all $R_0$, implying then:
\beq
{\sf b}_{n}(\epsilon\,\varphi(R_0))\leq {\sf b}_{n}(\epsilon\,\varphi_2).
\label{eq:boundphi32}
\eeq
Applying~(\ref{eq:boundphi31}) in~(\ref{eq:totalprobapp31}), and~(\ref{eq:boundphi32}) in~(\ref{eq:totalprobapp32}), and further using the definition of ${\sf b}_N(\delta)$ in~(\ref{eq:bndeldefapp}) we finally get~(\ref{eq:boundR0app}).
In order to obtain~(\ref{eq:boundR1app}), we must repeat the same steps shown above to the proof of Lemma~$2$ in~\cite{HanLuLiuJMLR}.
\end{IEEEproof}

Let us comment briefly on the structure of the bounds in~(\ref{eq:boundR0app}) and~(\ref{eq:boundR1app}). 
For clarity of presentation, we focus on~(\ref{eq:boundR0app}), since similar considerations would apply to~(\ref{eq:boundR1app}). 
We notice that the bound in~(\ref{eq:boundphi32}) could be used for all possible realizations of $\bm{R}_0$, while we used it only in~(\ref{eq:totalprobapp32}). This is because in~(\ref{eq:totalprobapp31}) we managed to obtain a better bound exploiting the degree-concentration properties exhibited by our system. 
In fact, examining the bound in~(\ref{eq:boundR0app}), we see that the term ${\sf b}_{n}(\epsilon\,\varphi_2)$ multiplies a quantity ${\sf b}_N(\delta)$ that is independent on the sample size $n$, and vanishes as the network size goes to infinity. 
The other term, ${\sf b}_{n}(\epsilon\,\varphi_1)$, goes to zero with the number of samples faster than ${\sf b}_{n}(\epsilon\,\varphi_2)$, since, exploiting the peculiarities of our model we obtained a higher (i.e., better) constant $\varphi_1>\varphi_2$. 
As a possible extension, one could try to refine also the term ${\sf b}_{n}(\epsilon\,\varphi_2)$ by exploring alternative bounding techniques --- see, e.g.,~\cite{BentoIbrahimiMontanari,RaoKipnisJavidiEldarGoldsmithCDC2016}.

\begin{corollary}[Scaling law useful for sample complexity]
\label{cor:bigomega}
Assume the same conditions used in Lemma~\ref{lem:concorr}. 
If
\beq
n(N)=\Omega\left(\, (Np_N)^2\log S \, \right),
\label{eq:samplecomplawapp}
\eeq
then we have that:
\beqa
\!\!\!\!\!\!\!\!\!\!\!\!\!\!\!&&\lim_{N\rightarrow\infty}
\P\left[\|[\widehat{\bm{R}}_{0,n(N)}]_{\mathcal{S}}-[\bm{R}_0]_{\mathcal{S}}\|_{\max}>\frac{\epsilon}{N p_N}\right]=0,
\label{eq:corollaryappR0}\\
\!\!\!\!\!\!\!\!\!\!\!\!\!\!\!&&\lim_{N\rightarrow\infty}
\P\left[\|[\widehat{\bm{R}}_{1,n(N)}]_{\mathcal{S}}-[\bm{R}_1]_{\mathcal{S}}\|_{\max}>\frac{\epsilon}{N p_N}\right]=0.
\label{eq:corollaryappR1}
\eeqa
\end{corollary}
\begin{IEEEproof}
We will prove the claim with reference to~(\ref{eq:corollaryappR0}), with the proof being identical for~(\ref{eq:corollaryappR1}). 
Examining~(\ref{eq:boundR0app}), we see that in order to get~(\ref{eq:corollaryappR0}) it suffices to guarantee that:
\beq
\lim_{N\rightarrow\infty} b_{n(N),S}\left(\frac{\epsilon\,\varphi_1}{N p_N}\right)=0, 
\label{eq:middleclaimapp}
\eeq
since the second term in~(\ref{eq:boundR0app}) goes to zero automatically as $N\rightarrow\infty$ due to the presence of ${\sf b}_N(\delta)$, whatever law $n(N)$ is chosen. 
Furthermore, in order to prove~(\ref{eq:middleclaimapp}), it suffices to examine only the second exponential term in~(\ref{eq:anphidef}), since the first term obviously vanishes provided that $n\rightarrow\infty$. 
Accordingly, the second exponential term in~(\ref{eq:anphidef}) can be written as:
\beq
\exp\left\{-\left(\frac{\sqrt{n} \,\epsilon\,\varphi_1}{N p_N} - \sqrt{2}\right)^2 + 2\log S\right\},
\eeq
whose exponent can be rewritten, after some straightforward algebra, as:
\beq
-\left[
\left(
\sqrt{\frac{n\epsilon^2\,\varphi_1^2}{2 (N p_N)^2\log S}}
-
\frac{1}{N p_N\sqrt{\log S}}
\right)^2
- 1
\right]\log S^2.
\eeq
Now, since $S/N\rightarrow \xi>0$ as $N\rightarrow\infty$, we see that $S\rightarrow\infty$. 
Therefore, the desired claim in~(\ref{eq:middleclaimapp}) is obtained if the quantity under brackets becomes asymptotically negative. Noticing that the term $(N p_N\sqrt{\log S})^{-1}$ vanishes as $N\rightarrow\infty$, we see that this condition is verified if we have, for $\epsilon'>0$:
\beq
\frac{n\epsilon^2\,\varphi_1^2}{2 (N p_N)^2\log S}>1 + \epsilon'
\eeq
from some $N$ onward, which corresponds to~(\ref{eq:samplecomplawapp}).
\end{IEEEproof}

The next two lemmas are auxiliary to the proof of Theorem~\ref{theor:samplecomplexity}.
\begin{lemma}
\label{lem:gennormprop}
Let ${\sf A}$, ${\sf R}_0$, ${\sf R}_1$, $\widehat{{\sf A}}$, $\widehat{{\sf R}}_0$ and $\widehat{{\sf R}}_1$ be square matrices of equal size, with:
\beq
\|{\sf A}\|_{\infty}\leq 1,~~{\sf R}_1={\sf A}{\sf R}_0.
\label{eq:Anormboundapp}
\eeq 
Assume that $\widehat{{\sf R}}_0$ is invertible and let 
\beq
\widehat{{\sf A}}=\widehat{{\sf R}}_1\widehat{{\sf R}}_0^{-1}.
\label{eq:hatAapp}
\eeq
Then we have that:
\beq
\|\widehat{{\sf A}}-{\sf A}\|_{\max}\leq 
\|\widehat{{\sf R}}_0^{-1}\|_1\,
\left(
\|\widehat{{\sf R}}_0-{\sf R}_0\|_{\max}
+
\|\widehat{{\sf R}}_1-{\sf R}_1\|_{\max}
\right).
\label{eq:normerr1}
\eeq
\end{lemma}
\begin{IEEEproof}
From~(\ref{eq:hatAapp}) we can write:
\beqa
\|\widehat{{\sf A}}-{\sf A}\|_{\max}
&=&
\|(\widehat{{\sf R}}_1-{\sf A}\widehat{{\sf R}}_0)\widehat{{\sf R}}_0^{-1}\|_{\max}
\nonumber\\
&\leq&
\|\widehat{{\sf R}}_0^{-1}\|_{1}\,
\|\widehat{{\sf R}}_1-{\sf A}\widehat{{\sf R}}_0\|_{\max}
\nonumber\\
&=&
\|\widehat{{\sf R}}_0^{-1}\|_{1}\,
\|\widehat{{\sf R}}_1-{\sf R}_1 \nonumber\\
&+& \underbrace{{\sf R}_1 - {\sf A}{\sf R}_0}_{=0} + {\sf A}{\sf R}_0 - {\sf A}\widehat{{\sf R}}_0\|_{\max}
\nonumber\\
&\leq&
\|\widehat{{\sf R}}_0^{-1}\|_{1}\,
\|\widehat{{\sf R}}_1-{\sf R}_1\|_{\max}\nonumber\\
&+&
\|\widehat{{\sf R}}_0^{-1}\|_{1}\,
\|{\sf A}\|_{\infty} \|\widehat{{\sf R}}_0-{\sf R}_0\|_{\max},
\eeqa
and the result in~(\ref{eq:normerr1}) follows in view of~(\ref{eq:Anormboundapp}).
\end{IEEEproof}

\begin{lemma}
\label{lem:gennormprop2}
Let ${\sf A}$, ${\sf R}_0$, ${\sf R}_1$, $\widehat{{\sf A}}$, $\widehat{{\sf R}}_0$ and $\widehat{{\sf R}}_1$ be square matrices of equal size, with:
\beq
\|{\sf A}\|_{\infty}\leq 1,~~{\sf R}_1={\sf A}{\sf R}_0.
\label{eq:Anormboundapp2}
\eeq 
Let the $i$-th row of $\widehat{{\sf A}}$ be a solution to the optimization problem:
\beq
\min_{x} 
\left\|
x \, \widehat{{\sf R}}_0 - [\widehat{{\sf R}}_1]_i
\right\|_{\infty}~~\textnormal{s.t. }\|x\|_1\leq 1,
\label{eq:optmaxnorm}
\eeq
where $x$ is a row vector and $[\widehat{{\sf R}}_1]_i$ is the $i$-th row of $\widehat{{\sf R}}_1$.

If ${\sf R}_0$ is invertible, then we have that:
\beq
\|\widehat{{\sf A}}-{\sf A}\|_{\max}
\leq
2\|{\sf R}_0^{-1}\|_{1}\left(
\|\widehat{{\sf R}}_0-{\sf R}_0\|_{\max}
+
\|\widehat{{\sf R}}_1-{\sf R}_1\|_{\max}
\right).
\label{eq:normerr2}
\eeq
\end{lemma}
\begin{IEEEproof}
We recall that for a vector the $\|\cdot\|_{\infty}$ norm amounts to the maximum absolute value of its entries, whereas for matrices it is the maximum absolute row sum. Accordingly, we see that the matrix ${\sf A}$ introduced in the statement of the lemma is a candidate solution to~(\ref{eq:optmaxnorm}) since by assumption $\|{\sf A}\|_{\infty}\leq 1$. 
As a result, a solution to the optimization problem in~(\ref{eq:optmaxnorm}) fulfills:
\beqa
\|\widehat{{\sf A}}\widehat{{\sf R}}_0-\widehat{{\sf R}}_1\|_{\max}
&\leq&
\|{\sf A}\widehat{{\sf R}}_0-\widehat{{\sf R}}_1\|_{\max}\nonumber\\
&=&
\|{\sf A}\widehat{{\sf R}}_0
-{\sf A}{\sf R}_0
+
{\sf R}_1
-\widehat{{\sf R}}_1\|_{\max}\nonumber\\
&\leq&
\|{\sf A}\|_{\infty}\,\|
\widehat{{\sf R}}_0
-{\sf R}_0
\|_{\max}
+
\|
{\sf R}_1
-\widehat{{\sf R}}_1\|_{\max}.\nonumber\\
\label{eq:demifinalbound}
\eeqa
On the other hand, we can write:
\beqa
\|\widehat{{\sf A}}-{\sf A}\|_{\max}
&=&
\|(\widehat{{\sf A}}{\sf R}_0-{\sf R}_1){\sf R}_0^{-1}\|_{\max}
\nonumber\\
&\leq&
\|{\sf R}_0^{-1}\|_{1}\,
\|\widehat{{\sf A}}{\sf R}_0-{\sf R}_1\|_{\max}\nonumber\\
&=&
\|{\sf R}_0^{-1}\|_{1}\,
\|\widehat{{\sf A}}{\sf R}_0
-\widehat{{\sf A}}\widehat{{\sf R}}_0 \nonumber\\
&+& \widehat{{\sf A}}\widehat{{\sf R}}_0
-\widehat{{\sf R}}_1+\widehat{{\sf R}}_1-{\sf R}_1\|_{\max}\nonumber\\
&\leq&
\|{\sf R}_0^{-1}\|_{1}\,
\|\widehat{{\sf A}}\|_{\infty}\,\|\widehat{{\sf R}}_0-{\sf R}_0\|_{\max}
\nonumber\\
&+&
\|{\sf R}_0^{-1}\|_{1}\,
\|\widehat{{\sf A}}\widehat{{\sf R}}_0-\widehat{{\sf R}}_1\|_{\max}
\nonumber\\
&+&
\|{\sf R}_0^{-1}\|_{1}\,
\|\widehat{{\sf R}}_1-{\sf R}_1\|_{\max},
\eeqa
and the claim in~(\ref{eq:normerr2}) follows from~(\ref{eq:demifinalbound}).
\end{IEEEproof}
\begin{IEEEproof}[Proof of Theorem~\ref{theor:samplecomplexity}]
We must evaluate the sample complexity of the estimators in~(\ref{eq:sampleGra}), ~(\ref{eq:sampleonelag}),~(\ref{eq:sampleres}) and~(\ref{eq:regularsamGra}).  
Examining the proof of Theorem~\ref{theor:samplecons}, and in particular Eq.~(\ref{eq:fnN}), we see that the sample complexity of $\widehat{\bm{A}}_{\mathcal{S},n}$ can be determined by finding a law $n(N)$ that ensures, for any arbitrarily small $\epsilon>0$:
\beq
\lim_{N\rightarrow\infty}
\P\left[
\|\widehat{\bm{A}}_{\mathcal{S},n(N)}-\widehat{\bm{A}}_{\mathcal{S}}\|_{\max}>\frac{\epsilon}{N p_N}
\right]=0.
\label{eq:relevantprobtoshow}
\eeq 
For what concerns the estimators in~(\ref{eq:sampleonelag}) and~(\ref{eq:sampleres}), the convergence in~(\ref{eq:relevantprobtoshow}) can be shown by examining the following probabilities, for $i=0,1$:
\beq
\P\left[\|[\widehat{\bm{R}}_{i,n}]_{\mathcal{S}}-[\bm{R}_i]_{\mathcal{S}}\|_{\max}>\frac{\epsilon}{N p_N}\right],
\label{eq:corrconvapp}
\eeq
and the claim of the theorem comes from application of Corollary~\ref{cor:bigomega}.
Let us consider then the regularized sample Granger estimator defined by~(\ref{eq:regularsamGra}). 
First we show that the {\em limiting} Granger estimator obeys the following bound:
\beq
\|\widehat{\bm{A}}^{\textnormal{(Gra)}}_{\mathcal{S}}\|_{\infty}\leq 1.
\label{eq:Granormboundapp}
\eeq
In view of~(\ref{eq:errGraMatform}) we can write:
\beq
\widehat{\bm{A}}^{\textnormal{(Gra)}}_{\mathcal{S}}=\bm{A}_{\mathcal{S}}+\bm{E}^{\textnormal{(Gra)}}
=
\bm{A}_{\mathcal{S}}+\bm{A}_{\mathcal{S},\mathcal{S}'}\underbrace{\bm{H}[\bm{A}^2]_{\mathcal{S}'\mathcal{S}}}_{\bm{F}},
\label{eq:normboundedby2}
\eeq
where we remark that the matrix $\bm{F}$ is nonnegative by construction. 
Exploiting~(\ref{eq:normboundedby2}) we have:
\beq
\sum_{j\in\mathcal{S}}
\widehat{\bm{a}}^{\textnormal{(Gra)}}_{ij}=
\sum_{j\in\mathcal{S}}
\bm{a}_{ij}+
\sum_{\ell\in\mathcal{S}'}\bm{a}_{i\ell}
\sum_{j\in\mathcal{S}}\bm{f}_{\ell j}
\leq 1-\rho\leq 1,
\eeq
where the last inequality follows since from Eq.~$(75)$ in~\cite{tomo} we know that:
\beq
\sum_{j\in\mathcal{S}} \bm{f}_{\ell j}\leq 1.
\label{eq:boundF}
\eeq
According to~(\ref{eq:Granormboundapp}), we can now use~(\ref{eq:normerr2}) and write:
\beqa
\lefteqn{
\|\widehat{\bm{A}}^{\textnormal{(reGra)}}_{\mathcal{S},n}-\widehat{\bm{A}}^{\textnormal{(Gra)}}_{\mathcal{S}}\|_{\max}}
\nonumber\\
&\leq&
2\|([\bm{R}_0]_{\mathcal{S}})^{-1}\|_{1}\left(
\|[\widehat{\bm{R}}_{0,n}]_{\mathcal{S}}-[\bm{R}_0]_{\mathcal{S}}\|_{\max}\right.
\nonumber\\
&+&
\left.\|[\widehat{\bm{R}}_{1,n}]_{\mathcal{S}}-[\bm{R}_1]_{\mathcal{S}}\|_{\max}
\right).
\label{eq:regulGrabound}
\eeqa
We conclude that the scaling law in~(\ref{eq:samplecomplaw}) will be sufficient for the regularized Granger estimator if we show that $\|([\bm{R}_0]_{\mathcal{S}})^{-1}\|_{1}$ is bounded by some constant. 
To this end, let 
\beq
\bm{Z}=\frac{I-\bm{A}^2}{\sigma^2}
\Leftrightarrow
\bm{R}_0=\bm{Z}^{-1}.
\eeq 
From one of the block matrix representations for the inverse matrix we have that~\cite{horn2012matrix}:
\beqa
\left([\bm{R}_0]_{\mathcal{S}}\right)^{-1}&=&
\bm{Z}_{\mathcal{S}} - \bm{Z}_{\mathcal{S}\mathcal{S}'}(\bm{Z}_{\mathcal{S}'})^{-1}\bm{Z}_{\mathcal{S}'\mathcal{S}}\nonumber\\
&=&
\frac{I_{\mathcal{S}} - [\bm{A}^2]_{\mathcal{S}} - [\bm{A}^2]_{\mathcal{S}\mathcal{S}'} \bm{H} [\bm{A}^2]_{\mathcal{S}'\mathcal{S}}}
{\sigma^2}
\nonumber\\
&=&
\frac{I_{\mathcal{S}} - [\bm{A}^2]_{\mathcal{S}} - [\bm{A}^2]_{\mathcal{S}\mathcal{S}'} \bm{F}}{\sigma^2}.
\label{eq:MILrep}
\eeqa
Using~(\ref{eq:boundF}) we see that the diagonal entries of $\left([\bm{R}_0]_{\mathcal{S}}\right)^{-1}$ are positive, whereas the off-diagonal entries are nonpositive since $\bm{F}$ and $\bm{A}$ are nonnegative matrices. Recalling that $\|\bm{A}^2\|_{\infty}=\rho^2$, from~(\ref{eq:MILrep}) we get: 
\beq
\left([\bm{R}_0]_{\mathcal{S}}\right)^{-1}\leq \frac{1+\rho^2}{\sigma^2}.
\label{eq:invsubmat}
\eeq
Let us finally prove that in the dense case where $p_N\rightarrow p>0$, the convergence in~(\ref{eq:corrconvapp}) with the scaling law in~(\ref{eq:samplecomplaw}) holds for the (non-regularized) Granger estimator in~(\ref{eq:sampleGra}). 
To this end, we expand $[\widehat{\bm{R}}_{0,n}]_{\mathcal{S}}$ as:
\beqa
[\widehat{\bm{R}}_{0,n}]_{\mathcal{S}}&=&
[\bm{R}_0]_{\mathcal{S}}+
\left(
[\widehat{\bm{R}}_{0,n}]_{\mathcal{S}}-[\bm{R}_0]_{\mathcal{S}}
\right)\nonumber\\
&=&
[\bm{R}_0]_{\mathcal{S}}\Big(
I_{\mathcal{S}}+\underbrace{
([\bm{R}_0]_{\mathcal{S}})^{-1}\,([\widehat{\bm{R}}_{0,n}]_{\mathcal{S}}-[\bm{R}_0]_{\mathcal{S}})}_{\bm{{\sf D}}}
\Big).\nonumber\\
\label{eq:usefulmatrepapp}
\eeqa
Now we observe that:
\beqa
\|\bm{{\sf D}}\|_1
&\leq& \|([\bm{R}_0]_{\mathcal{S}})^{-1}\|_1\,\,
\|([\widehat{\bm{R}}_{0,n}]_{\mathcal{S}}-[\bm{R}_0]_{\mathcal{S}})\|_1\nonumber\\
&\leq&
\frac{1+\rho^2}{\sigma^2}
\,S \,\,
\|([\widehat{\bm{R}}_{0,n}]_{\mathcal{S}}-[\bm{R}_0]_{\mathcal{S}})\|_{\max}\nonumber\\
&=&
\frac{1+\rho^2}{\sigma^2 p_N}\,\frac{S}{N}\,
\left(N p_N \,\,\|([\widehat{\bm{R}}_{0,n}]_{\mathcal{S}}-[\bm{R}_0]_{\mathcal{S}})\|_{\max}\right)\nonumber\\
&\leq&
\left(\frac{(1+\rho^2)\xi}{\sigma^2 p}+\epsilon\right)\,\epsilon=\epsilon',\nonumber\\
\eeqa
where: $i)$ in the second-to-last inequality we used~(\ref{eq:invsubmat}) and the fact that for an $S\times S$ matrix the $\|\cdot\|_1$ norm is upper bounded by $S$ times the $\|\cdot\|_{\max}$ norm; $ii)$ we used the fact that $S/N\rightarrow \xi$ and $p_N\rightarrow p>0$ as $N\rightarrow\infty$; and $iii)$ in view of~(\ref{eq:corrconvapp}), the last inequality holds (and, hence, $[\widehat{\bm{R}}_{0,n}]_{\mathcal{S}}$ is invertible) with probability converging to $1$ as $N\rightarrow\infty$.
Finally, using in~(\ref{eq:usefulmatrepapp}) a known result about the inversion of perturbed matrices we get~\cite{horn2012matrix}:
\beqa
\|([\widehat{\bm{R}}_{0,n}]_{\mathcal{S}})^{-1}\|_1
&\leq& 
\|([\bm{R}_{0}]_{\mathcal{S}})^{-1}\|_1\,\,
\|(I_{\mathcal{S}}+\bm{{\sf D}})^{-1}\|_1\nonumber\\
&\leq&
\frac{(1+\rho^2)(1-\|\bm{{\sf D}}\|_1)^{-1}}{\sigma^2}\nonumber\\
&\leq& \frac{1+\rho^2}{\sigma^2(1-\epsilon')},
\eeqa
and the proof of the theorem is now complete in view of Lemma~\ref{lem:gennormprop}.
\end{IEEEproof}

\end{appendices}


\end{document}